\documentclass{article}
\usepackage{geometry}

\usepackage{amsthm}
\newtheorem{theorem}{Theorem}[section]
\newtheorem{corollary}{Corollary}[theorem]
\newtheorem{lemma}[theorem]{Lemma}

\newtheorem{remark}{Remark}

\usepackage{bm}
\usepackage{amsmath}
\usepackage{amssymb}
\usepackage{stmaryrd}
\usepackage{mismath}

\usepackage[ruled,linesnumbered]{algorithm2e}

\usepackage{graphicx}
\graphicspath{{./figures/}}
\usepackage{subcaption}

\usepackage{cite}
\usepackage{csquotes}

\usepackage{tabu}
\usepackage{booktabs}
\usepackage{diagbox}
\usepackage{multirow}

\usepackage{xcolor}

\usepackage[shortlabels]{enumitem}

\usepackage{caption}
\usepackage{float}

\usepackage{authblk}
\usepackage{indentfirst}

\usepackage{todonotes}
\usepackage{siunitx}
\usepackage{hyperref}
\usepackage{cleveref}

\newcounter{stepnum}

\begin{document}
\title{A Hybrid Discontinuous Galerkin - Neural Network Method for Solving Hyperbolic Conservation Laws with Temporal Progressive Learning}

\author[1]{Yan Shen}
\author[1,2]{Jingrun Chen}
\author[2, \thanks{Corresponding author:wukekever@ustc.edu.cn}]{Keke Wu}

\affil[1]{\small School of Mathematical Sciences, University of Science and Technology of China, Hefei, 230026, People's Republic of China}
\affil[2]{\small School of Mathematical Sciences and Suzhou Institute for Advanced Research, University of Science and Technology of China, Suzhou 215127, People’s Republic of China}

\maketitle

\begin{abstract}
For hyperbolic conservation laws, traditional methods and physics-informed neural networks (PINNs) often encounter difficulties in capturing sharp discontinuities and maintaining temporal consistency.
To address these challenges, we introduce a hybrid computational framework by coupling discontinuous Galerkin (DG) discretizations with a temporally progressive neural network architecture. Our method incorporates a structure-preserving weak-form loss—combining DG residuals and Rankine–Hugoniot jump conditions—with a causality-respecting progressive training strategy. The proposed framework trains neural networks sequentially across temporally decomposed subintervals, leveraging pseudo-label supervision to ensure temporal coherence and solution continuity. This approach mitigates error accumulation and enhances the model’s capacity to resolve shock waves and steep gradients without explicit limiters. Besides, a theoretical analysis establishes error bounds for the proposed framework, demonstrating convergence toward the physical solution under mesh refinement and regularized training. Numerical experiments on Burgers and Euler equations show that our method consistently outperforms standard PINNs, PINNs-WE, and first-order DG schemes in both accuracy and robustness, particularly in capturing shocks and steep gradients. These results highlight the promise of combining classical discretization techniques with machine learning to develop robust and accurate solvers for nonlinear hyperbolic systems.
\end{abstract}

\section{Introduction}

Hyperbolic conservation laws govern a wide range of physical phenomena, including fluid flow, wave propagation, and gas dynamics. These systems generally show strong nonlinear characteristics, resulting in shock waves and contact discontinuities. These features have been challenging to accurately simulate using numerical methods.

Classical numerical methods such as finite difference \cite{leveque2002finite} and finite element \cite{brenner2008mathematical} schemes are effective in many smooth regions but often suffer from numerical diffusion or spurious oscillations near discontinuities. To improve robustness, flux-based solvers such as the Godunov \cite{godunov1959difference} and HLLC \cite{toro2013riemann} schemes have been widely adopted. More advanced methods—including ENO, WENO \cite{harten1997uniformly, jiang1996efficient}, and discontinuous Galerkin (DG) schemes \cite{cockburn2003discontinuous}—can provide higher accuracy, but still risk generating nonphysical oscillations near steep gradients. While slope limiters \cite{sweby1984high} can mitigate these artifacts, they introduce added complexity and computational overhead.

In recent years, machine learning, especially deep neural networks (DNNs), has emerged as a compelling alternative for solving partial differential equations (PDEs). Neural networks operate independently of structured grids and are well-suited to approximating high-dimensional functions, which makes them attractive for addressing the curse of dimensionality. Their nonlinear expressiveness and adaptability offer potential advantages in modeling intricate solution structures, such as shocks and rarefaction waves.

A range of techniques have been developed to exploit these capabilities. Magiera et al. \cite{magiera2020constraint} introduced constraint-aware networks for solving Riemann problems in a data-driven manner. Mao et al. \cite{mao2020physics} applied physics-informed neural networks (PINNs) \cite{RAISSI2019686} with randomized clustering to enhance shock resolution in high-speed flows \cite{cai2021physics}. Other efforts have explored shock detection without labeled data. Despite some progress, there are still challenges. As noted in \cite{wang2022pinn}, DNNs exhibit a spectral bias toward low-frequency modes, and vanilla PINNs relying on strong-form residuals frequently encounter difficulties with training stability and resolving sharp discontinuities.

To address these limitations, several improvements have been proposed. Patel et al. \cite{patel2022} designed control-volume PINNs inspired by finite volume schemes, integrating entropy inequalities and TVD constraints. Michoski et al. \cite{michoski2020} introduced artificial viscosity via modified equations to stabilize solutions, despite the cost of conservation fidelity and increased parameter sensitivity. The computation of second-order derivatives in these settings is also resource-intensive. Papados et al. \cite{papados2021} proposed W-PINNs-DE, emphasizing initial condition training and spatial extension, though at the cost of resolution near shocks. Liu et al. \cite{liu2024discontinuity} introduced Physics-Informed Neural Networks with Equation Weight (PINNs-WE), which adaptively adjusts loss function weights to better capture discontinuities. On the optimization side, Lu et al. \cite{lu2021deepxde} introduced residual-based adaptive sampling, while others employed gradient-aware weighting to guide training \cite{ferrer2024, mao2023physics}. More recently, Zhou et al. \cite{zhou2024} presented the Relaxation Neural Network framework, which augments PINNs with auxiliary systems but introduces additional modeling constraints and complexity.

A fundamental limitation across many of these approaches is the neglect of temporal causality. Hyperbolic systems evolve along characteristic curves at finite speeds, and classical solvers respect this structure through sequential time-stepping. In contrast, standard PINNs simultaneously approximate the solution over the entire spacetime domain, without enforcing causality. As a result, they may overfit to later-time data before learning the initial state, leading to unphysical behavior and accumulated error.

In this work, we propose a hybrid framework that integrates DG discretization with neural networks through a temporally progressive training strategy. The time domain is partitioned into subintervals, and a neural network is trained for each interval in sequence, using pseudo-labels from the preceding stage as supervision. This structure conforms to the causal dynamics of hyperbolic systems and substantially reduces temporal error accumulation.
To enhance stability and accuracy, a physics-guided loss function is designed, incorporating DG weak residuals, Rankine–Hugoniot jump conditions, and initial and boundary constraints. This formulation facilitates sharp shock representation without artificial viscosity or explicit limiters.

Our contributions are summarized as follows:
\begin{itemize}
  \item \textbf{Hybrid DG-Neural Network Framework:} We combine the strength of DG discretizations with the flexibility of neural networks to approximate solutions with strong nonlinearities.
  \item \textbf{Progressive Temporal Training:} We introduce a sequential training scheme that respects causality and improves long-time stability.
  \item \textbf{Physics-Guided Loss Design:} We formulate a composite loss incorporating DG residuals, jump conditions, and pseudo-label supervision to enhance accuracy near discontinuities.
  \item \textbf{Theoretical Guarantees:} We provide an error analysis that establishes bounds on the total prediction error, proving convergence to the physical solution under mesh refinement and regularized training.
  \item \textbf{Numerical Validation:} Comprehensive experiments on Burgers and Euler equations demonstrate superior accuracy and robustness over standard PINNs, PINNs-WE, and first-order DG methods, particularly in solving shocks and steep gradients.
\end{itemize}

The rest of this paper is organized as follows. Section \ref{sec:preliminaries} reviews the DG method, the challenges of using PINNs with hyperbolic systems, and the potential of the Progressive Neural Networks (PNNs) to solve the time-dependent problems. Section \ref{sec:methods} presents the proposed hybrid architecture and loss formulation. Theoretical analysis is provided in Section \ref{sec:theory}, followed by numerical experiments in Section \ref{sec:experiments}. Conclusions and directions for future research are given in Section \ref{sec:conclusion}.

\section{Preliminaries}
\label{sec:preliminaries}
\subsection{Discontinuous Galerkin Method}

Consider hyperbolic conservation laws
\begin{equation}\label{con law}
    \partial_t\mathbf{u}(\mathbf{x},t) + \nabla_\mathbf{x} \cdot \mathbf{f}(\mathbf{u}(\mathbf{x},t)) = 0, \quad \mathbf{x} = (x_1, \cdots, x_d) \in \Omega \subset \mathbb{R}^d,\quad t\in\mathcal{T}
\end{equation}
with corresponding initial and boundary conditions\cite{leveque1992numerical}.

For clarity, our analysis concentrates on the scalar field $u(x,t)$ defined over a one-dimensional spatial domain $\Omega = [a, a+L]$, discretized into $N$ computational cells $I_j = [x_{j-\frac{1}{2}}, x_{j+\frac{1}{2}}]$, where $j=1, \ldots, N$. The temporal domain is given by $\mathcal{T} = [0,T]$. We construct the finite element space
\begin{equation}
\label{eqn: test-space}
V_h = \{v(x): v|_{I_j} \in P^q(I_j), \quad j=1, \ldots, N\},
\end{equation}
where $P^q(I_j)$ denotes the space of polynomials with degree $\le q$ restricted to $I_j$.

Multiplying \eqref{con law} by a test function $v \in V_h$ and integrating by parts over each cell yields the semi-discrete DG formulation:
\begin{equation}\label{local form}
    \int_{I_j} u_t v \, \mathrm{d}x - \int_{I_j} f(u) v_x \, \mathrm{d}x + f_{j+\frac{1}{2}} v_{j+\frac{1}{2}}^- - f_{j-\frac{1}{2}} v_{j-\frac{1}{2}}^+ = 0,
\end{equation}
where the interface fluxes $f_{j+\frac{1}{2}}$ are replaced by numerical fluxes $\hat{f}_{j+\frac{1}{2}} = \hat{f}(u_{j+\frac{1}{2}}^-, u_{j+\frac{1}{2}}^+)$, and similar to $\hat{f}_{j-\frac{1}{2}}$, satisfying consistency, monotonicity, and Lipschitz continuity for stability \cite{leveque1992numerical}. We adopt the Lax–Friedrichs numerical flux in the calculation in the following text:
\begin{equation}\label{LF}
    f_{j+\frac{1}{2}}^{\mathrm{LF}} = \frac{1}{2} \left[ f(u_{j+\frac{1}{2}}^-) + f(u_{j+\frac{1}{2}}^+) - \alpha (u_{j+\frac{1}{2}}^+ - u_{j+\frac{1}{2}}^-) \right], \quad \alpha = \max_u |f'(u)|.
\end{equation}

Let $u_h \in V_h$ represent the numerical approximation of the solution $u$. The weak formulation is given by
\begin{equation}\label{weak form}
\int_{I_j} (u_h)_t v \mathrm{d}x - \int_{I_j} f(u_h) v_x \mathrm{d}x + \hat{f}_{j+\frac{1}{2}} v_{j+\frac{1}{2}}^- - \hat{f}_{j-\frac{1}{2}} v_{j-\frac{1}{2}}^+ = 0, \quad \forall v \in V_h.
\end{equation}

By projecting $u_h$ onto the local basis $\{\phi_i(x)\}_{i=0}^q$ of $V_h$, we obtain
\begin{equation}
(u_h)_i = \langle u_h, \phi_i \rangle_{I_j}, \quad 0 \le i \le q,
\end{equation}
where $\langle \cdot, \cdot \rangle$ represents the standard $L^2$ inner product, and $\{\phi_i(x)\}$ are canonical basis functions.

For numerical computations, we employ orthonormal Legendre polynomials on each cell:
\begin{equation}
\phi_0^{(j)}(x) = \frac{1}{\sqrt{h_j}}, \quad
\phi_1^{(j)}(x) = \frac{2\sqrt{3}}{h_j^{3/2}} (x - x_j), \quad \ldots
\end{equation}
with $h_j = x_{j+\frac{1}{2}} - x_{j-\frac{1}{2}}$. This basis choice remains consistent throughout our formulation. The approximate solution is expressed as
\begin{equation}
\tilde{u}_h = \sum_{i=0}^q (u_h)_i \phi_i.
\end{equation}

Substituting the approximate solution $\tilde{u}_h$ into the weak formulation \eqref{weak form} and choosing $v = \phi_m$ for $0 \le m \le q$ yields
\begin{equation}
\int{I_j} \left( \sum_{i=0}^q (u_h)_i \phi_i \right)_t \phi_m \mathrm{d}x +
\hat{f}_{j+\frac{1}{2}} (\phi_m)_{j+\frac{1}{2}}^- - \hat{f}_{j-\frac{1}{2}} (\phi_m)_{j-\frac{1}{2}}^+ = 0.
\end{equation}

Applying a forward Euler temporal discretization with time step $\Delta t$, we obtain the discrete formulation:
\begin{equation}
\label{eq:euler-forward}
\begin{aligned}
&\sum_{i=0}^q (u_h)_i^{n+1} \int{I_j} \phi_i \phi_m \mathrm{d}x - \sum_{i=0}^q (u_h)_i^{n} \int_{I_j} \phi_i \phi_m \mathrm{d}x \\
=& \Delta t \left[ \int_{I_j} f\left( \sum_{i=0}^q (u_h)_i^{n} \phi_i \right) \phi_m' \mathrm{d}x - \hat{f}_{j+\frac{1}{2}}^{n} (\phi_m)_{j+\frac{1}{2}}^- + \hat{f}_{j-\frac{1}{2}}^{n} (\phi_m)_{j-\frac{1}{2}}^+ \right].
\end{aligned}
\end{equation}

Introducing the coefficient vector $\tilde{\tilde{u}}_h = [(u_h)_0, (u_h)_1, \dots, (u_h)_q]^T$, source term $\mathbf{b}$, and mass matrix $\mathcal{M} \in \mathbb{R}^{(q+1)\times (q+1)}$ with entries
\begin{equation}
m_{rs} = \int_{I_j} \phi_r \phi_s , \mathrm{d}x,
\end{equation}
leads to the compact DG scheme formulation:
\begin{equation}\label{DG scheme}
\mathcal{M} \tilde{\tilde{u}}_h^{n+1} = \mathcal{M} \tilde{\tilde{u}}_h^{n} + \Delta t \cdot \mathbf{b}.
\end{equation}

\subsection{Physics-Informed Neural Networks (PINNs)}

Physics-Informed Neural Networks (PINNs) provide a robust framework to solve PDEs by embedding the governing equations and boundary conditions into the training loss \cite{RAISSI2019686}. Unlike completely data-driven methods, PINNs enforce physical laws as soft constraints, enabling effective generalization with limited data.

To approximate solutions of the one-dimensional hyperbolic conservation law \eqref{con law}, we represent the state by a neural network $u_\theta(x,t)$ mapping spatio-temporal inputs to predicted values. Training minimizes a composite loss that enforces the PDE residual, initial, and boundary conditions.

The physics-informed loss penalizes the PDE residual:
\begin{equation}
\label{eq:continuous_loss}
\begin{aligned}
    \mathcal{L}_{\mathrm{int}}(\theta) &= \frac{1}{|\mathcal{T}\times\Omega|}\int_0^T \int_{\Omega} \left| \frac{\partial u_\theta}{\partial t} + \frac{\partial}{\partial x} f(u_\theta) \right|^2 \mathrm{d}x \, \mathrm{d}t\\
    &\approx \frac{1}{N_{\mathrm{int}}} \sum_{i=1}^{N_{\mathrm{int}}} \left| \frac{\partial}{\partial t}u_\theta(x_{\mathrm{int}}^i, t_{\mathrm{int}}^i) + \frac{\partial}{\partial x} f(u_\theta(x_{\mathrm{int}}^i, t_{\mathrm{int}}^i)) \right|^2,
\end{aligned}
\end{equation}
where 
$\{ (x_{\mathrm{int}}^i, t_{\mathrm{int}}^i) \}_{i=1}^{{N_\mathrm{int}}}$ are the sampled collocation points in the domain $[0, T] \time \Omega$,
and derivatives are computed via automatic differentiation.

The boundary and initial conditions are enforced through additional loss terms:
\begin{equation}
\mathcal{L}_{\mathrm{bdy}}(\theta) = \frac{1}{|\mathcal{T}\times\partial\Omega|}\int_{\mathcal{T}\times\partial\Omega} |u_\theta - u_{\mathrm{bdy}}|^2 \mathrm{d}s\mathrm{d}t, \quad
\mathcal{L}_0(\theta) = \frac{1}{|\Omega|}\int{\Omega} |u_\theta(x,0) - u_0(x)|^2 \mathrm{d}x,
\end{equation}
which are approximated using collocation points $\{(x_{\mathrm{bdy}}^i,t_{\mathrm{bdy}}^i)\}_{i=1}^{N{\mathrm{bdy}}}$ on the boundary and $\{x_0^i\}_{i=1}^{N_0}$ at $t=0$:
\begin{equation}
\mathcal{L}_{\mathrm{bdy}}(\theta) \approx \frac{1}{N_{\mathrm{bdy}}}\sum_{i=1}^{N_{\mathrm{bdy}}} |u_\theta(x_{\mathrm{bdy}}^i,t_{\mathrm{bdy}}^i) - u_{\mathrm{bdy}}^i|^2, \quad
\mathcal{L}_0(\theta) \approx \frac{1}{N_0}\sum_{i=1}^{N_0} |u_\theta(x_0^i,0) - u_0^i|^2.
\end{equation}

The total loss combines these components with weights $\lambda_{\mathrm{int}}, \lambda_{\mathrm{bdy}}, \lambda_0$:
\begin{equation}
    \mathcal{L}(\theta) = \lambda_{\mathrm{int}} \mathcal{L}_{\mathrm{int}} + \lambda_{\mathrm{bdy}} \mathcal{L}_{\mathrm{bdy}} + \lambda_0 \mathcal{L}_0.
\end{equation}

Although vanilla PINNs with uniform loss weights ($\lambda_{\mathrm{int}}=\lambda_{\mathrm{bdy}}=\lambda_0=1$) perform well for simple hyperbolic equations given suitable parameterization, they face significant challenges when applied to complex hyperbolic conservation laws. Specifically, their accuracy deteriorates and errors accumulate over time \cite{wang2022pinn, zhang2024physicsinformedneuralnetworks}.

To illustrate this limitation, we implement a conventional PINN architecture (4 hidden layers, 64 neurons per layer) for the Sod shock tube problem \eqref{Sod}. The network is trained using a two-phase approach: initial optimization via Adam (30,000 iterations at $10^{-3}$ learning rate) followed by L-BFGS refinement ($10^{-2}$ learning rate).

\begin{equation}
\begin{cases}
   & \frac{\partial}{\partial t} \begin{bmatrix} 
        \rho \\ \rho u \\ E 
    \end{bmatrix} 
    + \nabla \cdot \begin{bmatrix}
        \rho u \\ \rho u^2 + p \\ u(E + p)
    \end{bmatrix}
= 0,\\
&(\rho(x,0),u(x,0),p(x,0))=\left\{
\begin{aligned}
    &(1.0, 0.0, 1.0),  & x < 1.0 \\
    &(0.125, 0.0, 0.1),  & x > 1.0
\end{aligned}\right.
\end{cases}
\label{Sod}
\end{equation}

\begin{figure}[!htbp]
    \centering
    \includegraphics[width=0.8\linewidth]{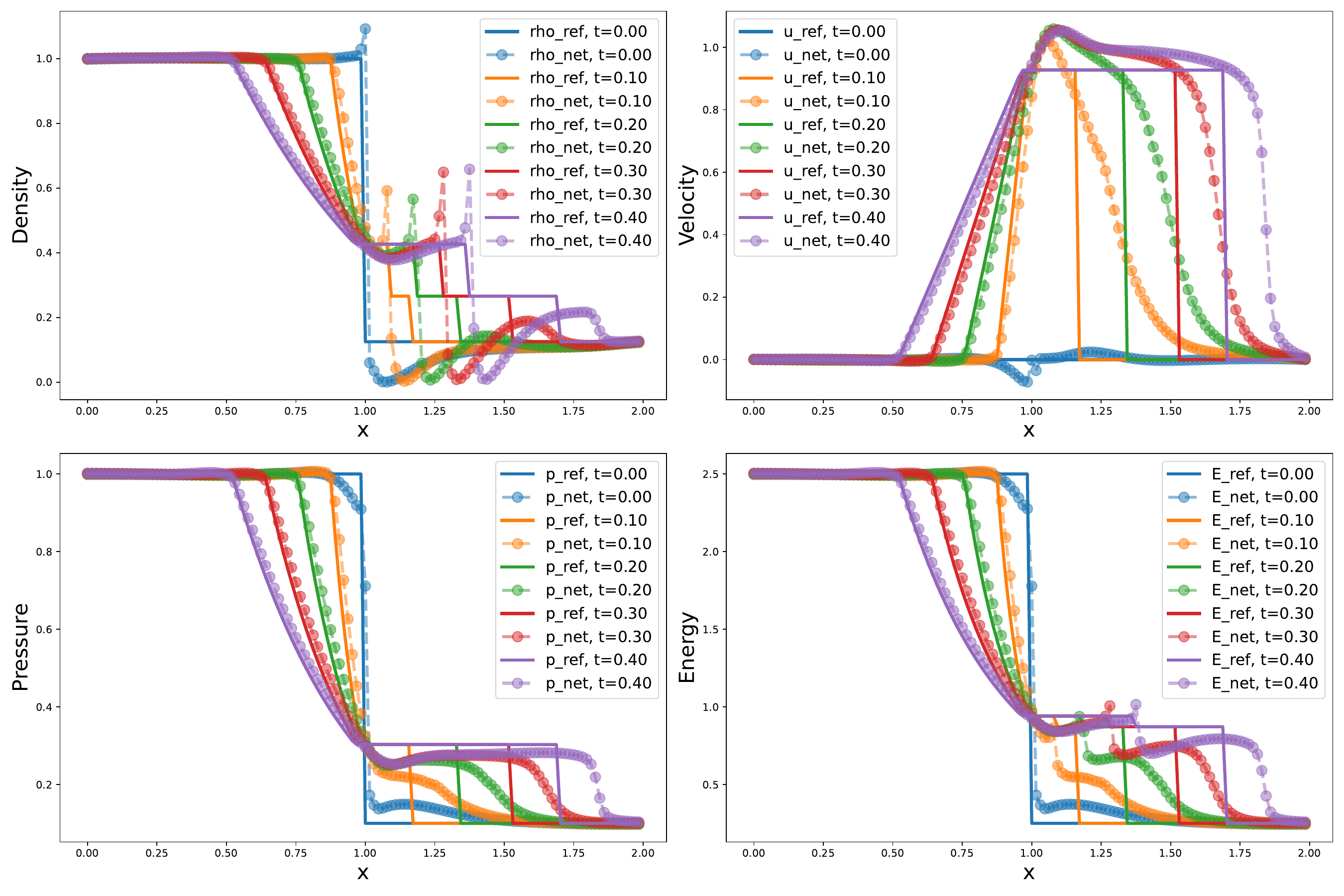}
    \caption{Vanilla PINN prediction for the 1D Sod shock tube. Solid lines: exact solution; scatter points: PINN predictions over time.}
    \label{class_pinn}
\end{figure}

As shown in Figure \ref{class_pinn}, while automatic differentiation offers implementation convenience, its application to hyperbolic conservation laws faces stability and convergence challenges due to solution discontinuities.
The development of hybrid frameworks that integrate classical numerical schemes with neural architectures has been spurred by this insight. This suggests a promising direction for future research.

\subsection{Progressive Neural Networks}

As stated in \cite{rusu2016progressive}, progressive neural networks (PNNs) are an advanced model architecture that demonstrates remarkable knowledge retention capability and facilitates continuous learning across multiple tasks.
By implementing lateral connections, these networks efficiently leverage features acquired from previous tasks, thereby mitigating catastrophic forgetting. A progressive network is built on a foundational column, which is a neural network made up of $L$ layers with hidden activations $h_i$ in $R_i$, where $n_i$ is the number of units at layer $i \le L$. The initial parameters $\theta^{(1)}$ are trained to convergence.  When transitioning to subsequent tasks, the network expands by adding new columns rather than modifying existing parameters. 
For instance, when addressing a second task, the original parameters $\theta^{(1)}$ are frozen (fixed and not participating in training) while a new column with parameters $\theta^{(2)}$ is initialized randomly. 
This new column incorporates lateral connections, allowing layer $h_i^{(2)}$ to receive inputs from both its preceding layer $h_{i-1}^{(2)}$ and the corresponding layer in the first column $h_{i-1}^{(1)}$.
This architectural pattern extends naturally to accommodate $K$ tasks. The activation function for any layer $i$ in the $k$-th column is given by:
\begin{equation}
    h_i^{(k)} = \sigma (\underbrace{W_i^{(k)}h_{i-1}^{(k)}}_{\text{learnable}} + \underbrace{\sum_{j<k} U_i^{(k:j)}h_{i-1}^{(j)}}_{\text{frozen}}).
\end{equation}
Here, $W_i^{(k)}$ denotes the weight matrix for layer $i$ of column $k$, $h_i^{(k)}$ represents the hidden layer activations at layer $i$ in column $k$, and $U_i^{(k:j)}$ represents the lateral connection matrices linking column $k$ to all previous columns $j$ (where $j < k$). 
The function $\sigma$ is an element-wise non-linear activation function that enables the network to learn complex, non-linear patterns and relationships in the data.

\begin{remark}
The key difference between PNNs and transfer learning is knowledge retention. Transfer learning risks forgetting during fine-tuning, whereas PNNs freeze prior parameters and use lateral connections to reuse features. This design prevents catastrophic forgetting and promotes effective, continual learning.
\end{remark}

\begin{figure}[!htbp]
    \centering
    \includegraphics[width=0.6\linewidth]{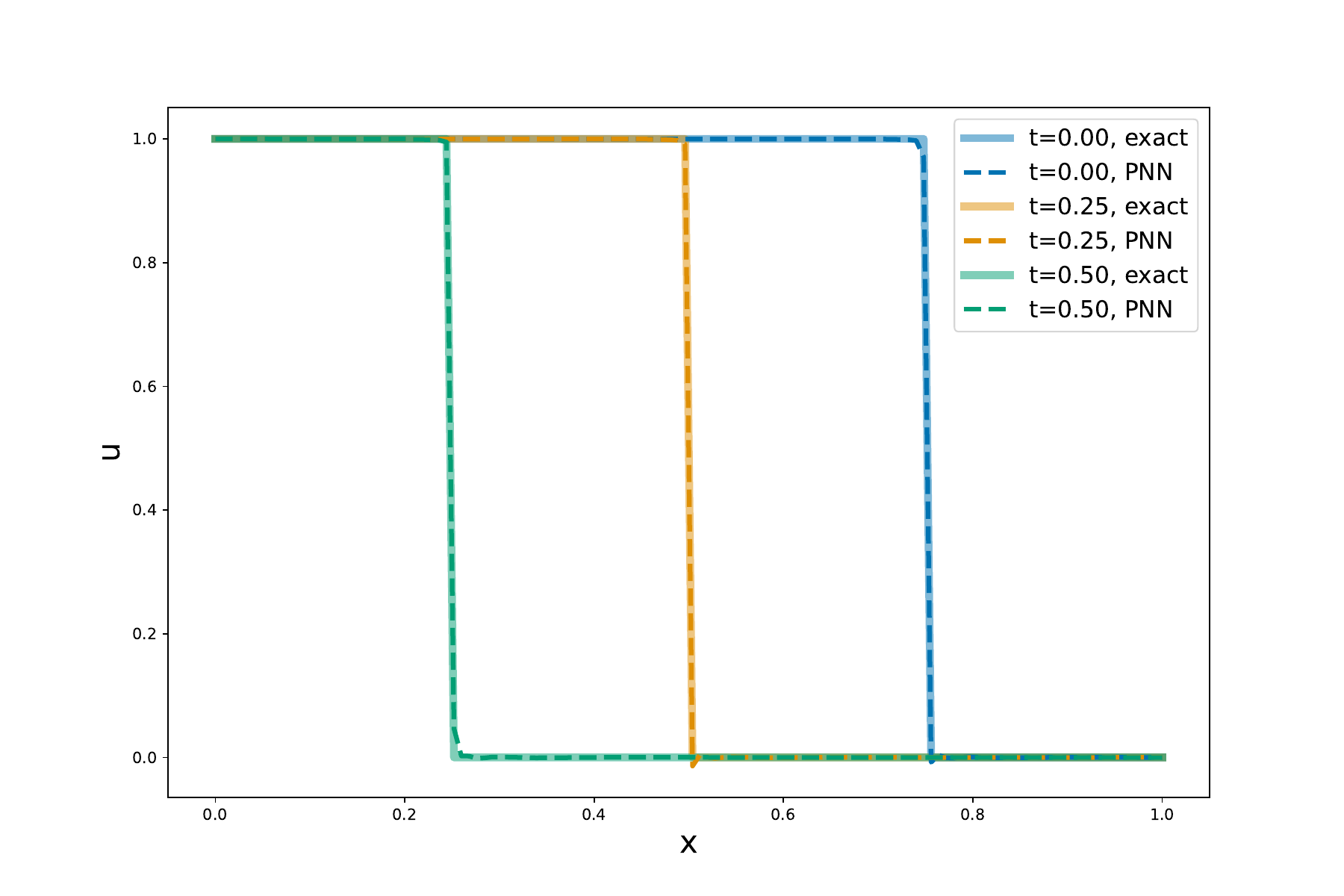}
    \caption{Exact vs. progressive neural network solutions for the linear advection equation across time intervals. Dashed lines: PNN predictions for sequential tasks.}
    \label{fig:pnn}
\end{figure}

For time-dependent PDEs, PNNs can align each column with a time subinterval. This partitions the temporal domain, enabling incremental learning of solution dynamics while stabilizing long-time simulations. We illustrate this on the linear advection equation:
\begin{equation}
    \left\{\begin{aligned}
        u_t - u_x &= 0, \\
        u(x,0) &= \begin{cases}
            1.0, & x < 0.75, \\
            0.0, & x > 0.75,
        \end{cases}
    \end{aligned}\right.
\end{equation}
which is decomposed into three sequential learning tasks. Training uses a two-stage optimization: 20,000 Adam iterations ($1\times10^{-3}$ learning rate), followed by 1,000 L-BFGS steps ($1\times10^{-2}$ learning rate). As Figure \ref{fig:pnn} shows, PNN predictions closely match exact solutions with errors near $10^{-3}$, demonstrating the approach’s promise for structured temporal learning in PDEs.

\section{Methodology}
\label{sec:methods}
\subsection{Overview of the Progressive Causal Learning Framework}
To address the challenges of maintaining temporal causality and accurately capturing discontinuities in hyperbolic conservation laws, we introduce a hybrid framework that integrates discontinuous Galerkin (DG) discretization, the spectral characteristics of neural networks, and progressive network architectures. This section details the architecture and loss design of our model.

Neural networks exhibit a spectral bias toward learning low-frequency components, as observed by Xu et al. \cite{JohnXu2020frequency}. This property allows them to act as natural smoothers, suppressing high-frequency numerical oscillations commonly encountered in hyperbolic problems. Leveraging this implicit regularization, our method avoids explicit limiters—often required in high-order DG schemes—while maintaining physical consistency and reducing computational complexity.

We adopt a progressive training strategy in the temporal domain to respect the causal structure of hyperbolic PDEs, in which information propagates along characteristics with finite speed. Instead of solving the whole spatio-temporal domain simultaneously, we divide the time interval into subdomains and associate a separate neural network with each. The training of these networks is done sequentially: outputs and latent representations from earlier intervals are used as pseudo-labels for guiding training in the next interval. This design reinforces temporal coherence and effectively alleviates the propagation of cumulative errors.

\subsection{Physics-Guided Loss Functions}

Building on the framework outlined above, we model the solution using a neural network $u_{\theta}(x, t)$ that approximates the solution at each spatio-temporal location. To ensure both physical fidelity and numerical stability, we construct a composite loss function incorporating physical constraints, numerical structure, and temporal regularization.

\paragraph{Initial and Boundary Loss}
To enforce the prescribed initial and boundary conditions, we define the corresponding loss terms based on integral formulations and their discrete approximations.

The boundary loss $\mathcal{L}_{\mathrm{bdy}}$ is defined as
\begin{equation}
    \label{eq:bdy_cont}
    \begin{aligned}
        \mathcal{L}_{\mathrm{bdy}}(\theta) &=\frac{1}{|\mathcal{T}\times\partial\Omega|} \int_0^T \int_{\partial \Omega} \left| u_{\theta}(x,t) - u_{\mathrm{bdy}}(x,t) \right|^2 \, \mathrm{d}S \, \mathrm{d}t \\
        &\approx \frac{1}{N_{\mathrm{bdy}}} \sum_{i=1}^{N_{\mathrm{bdy}}}\left| u_{\theta}(x_{\mathrm{bdy}}^{i}, t_{\mathrm{bdy}}^{i}) - u_{\mathrm{bdy}}(x_{\mathrm{bdy}}^{i}, t_{\mathrm{bdy}}^{i}) \right|^2,
    \end{aligned}
\end{equation}
where $\{(x_{\mathrm{bdy}}^i, t_{\mathrm{bdy}}^i)\}_{i=1}^{N_{\mathrm{bdy}}}$ are collocation points sampled along the spatial boundary $\partial \Omega$ over the time interval $\mathcal{T}$.

The initial condition loss $\mathcal{L}_{\mathrm{IC}}$ is defined as
\begin{equation}
    \label{eq:ic_cont}
    \mathcal{L}_{\mathrm{IC}}(\theta) =\frac{1}{|\Omega|} \int_{\Omega} \left| u_{\theta}(x,0) - u_0(x) \right|^2 \, \mathrm{d}x 
    \approx \frac{1}{N_{\mathrm{IC}}} \sum_{i=1}^{N_{\mathrm{IC}}}
    \left| u_{\theta}(x_{\mathrm{IC}}^{i}, 0) - u_0(x_{\mathrm{IC}}^{i}) \right|^2,
\end{equation}
with $\{x_{\mathrm{IC}}^i\}_{i=1}^{N_{\mathrm{IC}}}$ denoting samples in the spatial domain at initial time $t=0$.

\paragraph{DG-inspired Structural Loss}
To embed the conservative structure of the governing equations into the network training, we propose a physics-guided loss term that rigorously enforces conservation laws. This loss penalizes the residual derived from a second-order Runge–Kutta (RK2) discretization. Formally, the loss function is defined as follows (with $\hat{\mathcal{T}} = [0, T - \Delta t]$):
\begin{equation} \label{eq:DG-loss-RK2}
\begin{aligned}
\mathcal{L}_{\mathrm{DG}}(\theta) 
&= \frac{1}{|\hat{\mathcal{T}}\times\Omega|}\int_0^{T - \Delta t} \int_\Omega d(x,t) \cdot \left( \frac{\partial u_\theta}{\partial t}(x,t) + \nabla \cdot f(u_\theta(x,t)) \right)^2 \, \mathrm{d}x \, \mathrm{d}t \\
&\approx \frac{1}{N_{\mathrm{DG}}} \sum_{i=1}^{N_{\mathrm{DG}}} \left[ d\left(x_{\mathrm{DG}}^i, t_{\mathrm{DG}}^i\right) \cdot \left( \frac{u_\theta^{n+1}(x_{\mathrm{DG}}^i, t_{\mathrm{DG}}^i) - u_\theta^n(x_{\mathrm{DG}}^i, t_{\mathrm{DG}}^i)}{\Delta t} + \frac{1}{2} \left( k_1^i + k_2^i \right) \right) \right]^2.
\end{aligned}
\end{equation}
where $\{(x_{\mathrm{DG}}^i, t_{\mathrm{DG}}^i)\}_{i=1}^{N_{\mathrm{DG}}}$ denotes a set of collocation points randomly sampled from the spatio-temporal domain. The residual terms $k_1^i$ and $k_2^i$ are constructed using a two-stage RK2 method with Lax-Friedrichs flux as \eqref{LF}, defined as:
\begin{equation} \label{eq:rk2-k1k2-discrete}
k_1^i = - \frac{1}{\Delta x} \left( f^{\mathrm{LF}}_{i+\frac{1}{2}} - f^{\mathrm{LF}}_{i-\frac{1}{2}} \right), \quad
k_2^i = - \frac{1}{\Delta x} \left( f^{\mathrm{LF},*}_{i+\frac{1}{2}} - f^{\mathrm{LF},*}_{i-\frac{1}{2}} \right).
\end{equation}
The starred quantities $f^{\mathrm{LF},*}$ correspond to fluxes evaluated using the intermediate RK2 state $u_\theta^{*,i} = u_\theta^n(x_{\mathrm{DG}}^i, t_{\mathrm{DG}}^i) + \Delta t \cdot k_1^i$.

Using finite difference approximations, we estimate the local solution jump $\llbracket u \rrbracket$, where $\llbracket \cdot \rrbracket$ denotes the jump operator quantifying differences across shock interfaces. For any scalar or vector quantity $u$, the jump is defined as $\llbracket u \rrbracket = u_l - u_r$, with $u_l$ and $u_r$ representing the left and right limits of $u$ at the interface, respectively.
This enables adaptive weighting of the residual loss in regions with sharp features. The weighting function $d(x,t)$ identifies abrupt variations or steep gradients and defines a binary weight:
\begin{equation}
    \label{eq:jump-weight}
    d(x,t) = 
    \begin{cases}
        0.1, & \llbracket u \rrbracket > \epsilon, \\
        1,   & \llbracket u \rrbracket \leq \epsilon.
    \end{cases}
\end{equation}

Unlike the approach of Liu et al. \cite{liu2024discontinuity}, which applies smaller weights solely to regions with sharp gradient declines, our method accounts for both rapid increases and decreases in solution gradients. This ensures more robust identification of singular structures, improving training stability and solution fidelity near discontinuities.

\begin{remark}
A key advantage of this formulation is eliminating explicit slope limiters or higher-order reconstructions (e.g., WENO), typically required to suppress spurious oscillations near discontinuities. This simplification stems from neural networks' inherent spectral bias—their tendency to learn smooth, low-frequency functions, implicitly regularizing solutions even for sharp gradients.
\end{remark}
\begin{remark}
After sampling spatio-temporal points $(x_{\mathrm{DG}}^i, t_{\mathrm{DG}}^i)$ within the domain, each $x_{\mathrm{DG}}^i$ is mapped to its corresponding computational cell for precise flux evaluation. Our experiments reveal resolution-dependent sampling behavior: for high resolutions (128/256 cells), half-grid point placement enhances accuracy and stability, while for coarse grids (32/64 cells), fully random in-cell sampling produces more robust training. This empirical, resolution-dependent approach serves as an effective practical prior, leaving more systematic adaptive collocation methods for future investigation.
\end{remark}
\begin{remark}
Alternatively, the DG-inspired loss can be reformulated through explicit RK2 time integration of the predicted solution. The network output $u_\theta(x,t)$ is advanced in time to generate the numerical approximation $u_h(x,t+\Delta t)$ using the same two-stage RK2 scheme similar to \ref{eq:euler-forward}. This formulation expresses the residual loss \eqref{eq:DG-loss-RK2} as the forward error between the network prediction $u_\theta(x,t+\Delta t)$ and RK2-evolved state $u_h(x,t+\Delta t)$, scaled by $\Delta t^{-1}$:
\begin{equation}
\label{eq:RK2-direct-loss}
\mathcal{L}_{\mathrm{RK2}}(\theta) = \frac{1}{|\hat{\mathcal{T}}\times\Omega|}\int_0^{T-\Delta t}\int_\Omega d(x,t) \cdot \left( \frac{u_\theta(x,t+\Delta t) - u_h(x,t+\Delta t)}{\Delta t} \right)^2 \, \mathrm{d}x \,\mathrm{d}t.
\end{equation}
While structurally equivalent to \eqref{eq:DG-loss-RK2}, this formulation offers distinct physical interpretation by casting the residual as a one-step prediction error. We employ this representation in the framework diagram for consistency with our implementation.
\end{remark}

\paragraph{Rankine-Hugoniot Loss}

To strengthen conservation law enforcement near discontinuities, we incorporate the Rankine-Hugoniot (RH) jump condition directly into the training objective. The condition is enforced at $N_{\mathrm{RH}}$ sampled points $\{(x_{\mathrm{RH}}^i, t_{\mathrm{RH}}^i)\}_{i=1}^{N{\mathrm{RH}}}$ throughout the computational domain, with the corresponding loss defined as:
\begin{equation}
    \label{eq:RH-loss}
    \begin{aligned}
    \mathcal{L}_{\mathrm{RH}}(\theta) 
    &= \frac{1}{|\mathcal{T}\times\Omega|}\int_0^{T} \int_{\Omega} \left| s \cdot \llbracket u_{\theta}(x,t) \rrbracket - \llbracket f(u_{\theta})(x,t) \rrbracket \right|^2 \, \delta_{\Gamma}(x,t) \, \mathrm{d}x \, \mathrm{d}t \\
    &\approx \frac{1}{N_{\mathrm{RH}}} \sum_{i=1}^{N_{\mathrm{RH}}} \left| s \cdot \llbracket u_{\theta}(x_{\mathrm{RH}}^{i}, t_{\mathrm{RH}}^{i}) \rrbracket - \llbracket f(u_{\theta})(x_{\mathrm{RH}}^{i}, t_{\mathrm{RH}}^{i}) \rrbracket \right|^2,
    \end{aligned}
\end{equation}
where $s$ represents the discontinuity propagation speed. This formulation ensures the learned solution satisfies physically admissible jump conditions at discontinuities, significantly improving both accuracy and stability of the neural network approximation near singular features.

\paragraph{Supervised Pseudo-Label Loss}

To maintain temporal coherence and reduce error propagation across time steps, we implement a supervised loss using pseudo-labels \cite{Lee2013PseudoLabel}. The time domain $\mathcal{T}$ is divided into $M$ uniform subintervals $[t_{0},t_{1}], [t_{1},t_{2}], \cdots, [t_{M-1},t_{M}]$, where $t_{0}=0$ and $t_{M}=T$. For each $k$-th time interval $[t_{k-1}, t_{k}]$, we define the corresponding neural network approximation as $u_{\theta^{(k)}}$ with the supervised loss:
\begin{equation}
    \mathcal{L}_{\mathrm{tmp}}^{(k)}(\theta) = \frac{1}{N_{\mathrm{sup}}}\sum_{i=1}
    ^{N_{\mathrm{sup}}}\vert u_{\theta^{(k)}}(x_{\mathrm{sup}}^{k,i},t_{\mathrm{sup}}^{k,i}
    )-u_{\theta^{(k-1)}}(x_{\mathrm{sup}}^{k,i},t_{\mathrm{sup}}^{k,i})\vert^{2},
\end{equation}
where $\{x_{\mathrm{sup}}^{k,i},t_{\mathrm{sup}}^{k,i}\}_{i=1}^{N_{\mathrm{sup}}}$
are uniformly sampled points in the domain $[a, a+L] \times [t_{0}, t_{k-1}]$. To
further capture this discrepancy in the frequency domain, we apply the Fast Fourier
Transform (FFT) to the prediction differences between adjacent neural networks,
i.e.,
$\{u_{\theta^{(k)}}(x_{\mathrm{sup}}^{k,i},t_{\mathrm{sup}}^{k,i})\}_{i=1}^{N_{\mathrm{sup}}}$
and
$\{u_{\theta^{(k-1)}}(x_{\mathrm{sup}}^{k,i}, t_{\mathrm{sup}}^{k,i})\}_{i=1}^{N_{\mathrm{sup}}}$,
leading to the frequency-domain loss:
\begin{equation}
    \mathcal{L}_{\mathrm{freq}}^{(k)}(\theta) = \frac{1}{N_{\mathrm{sup}}}\Vert  \mathcal{F}(\{u_{\theta^{(k)}}
    (x_{\mathrm{sup}}^{k,i},t_{\mathrm{sup}}^{k,i})\}_{i=1}^{N_{\mathrm{sup}}}) -
    \mathcal{F}(\{u_{\theta^{(k-1)}}(x_{\mathrm{sup}}^{k,i},t_{\mathrm{sup}}^{k,i}
    )\}_{i=1}^{N_{\mathrm{sup}}}) \Vert_{2}^{2},
\end{equation}
where $\mathcal{F}$ denotes the FFT operator. This frequency-domain loss helps keep adjacent neural networks consistent across different frequency components. This is particularly important for capturing both low and high-frequency features of the solution. The total supervised loss for the $k$-th time block is as follows:
\begin{equation}
    \mathcal{L}_{\mathrm{sup}}^{(k)}(\theta) = \omega \mathcal{L}_{\mathrm{tmp}}^{(k)}
    (\theta ) + (1 - \omega) \mathcal{L}_{\mathrm{freq}}^{(k)}(\theta),
\end{equation}
where $\omega$ is a weighting parameter balancing the temporal and frequency domain losses. Generally, we set $\omega=1/2$ to give equal weight to both losses.

\begin{remark}
This supervised loss mechanism, combined with pseudo-labeling, enables our model to maintain solution stability while accurately resolving both low- and high-frequency solution components. The pseudo-labels generated by the preceding network $u_{\theta^{(k-1)}}$ serve as approximate solution data containing inherent noise, which automatically provides beneficial regularization during training.
This methodology shares conceptual similarities with noise injection techniques, where controlled perturbations improve model robustness. While these pseudo-labels do not constitute exact ground truth, they offer valuable supervisory signals that direct the learning process while maintaining greater flexibility than purely supervised approaches. The controlled uncertainty in training data fosters more robust and generalizable feature learning beyond simple pattern memorization - a crucial advantage for conservation law problems.
\end{remark}

\subsection{Network Architecture and Training Strategy}
Our framework employs a specialized neural network architecture for solving time-dependent PDEs. The temporal domain is decomposed into discrete subintervals, each handled by a dedicated neural network trained sequentially. This design intrinsically enforces temporal causality by eliminating backward error propagation - a fundamental limitation in conventional PINNs. While existing approaches like causal reweighting \cite{WANG2024116813} impose temporal directionality through loss scaling, our method structurally embeds causality via progressive supervision and modular network expansion.

\begin{figure}[!htbp]
    \centering
    \includegraphics[width=0.8\linewidth]{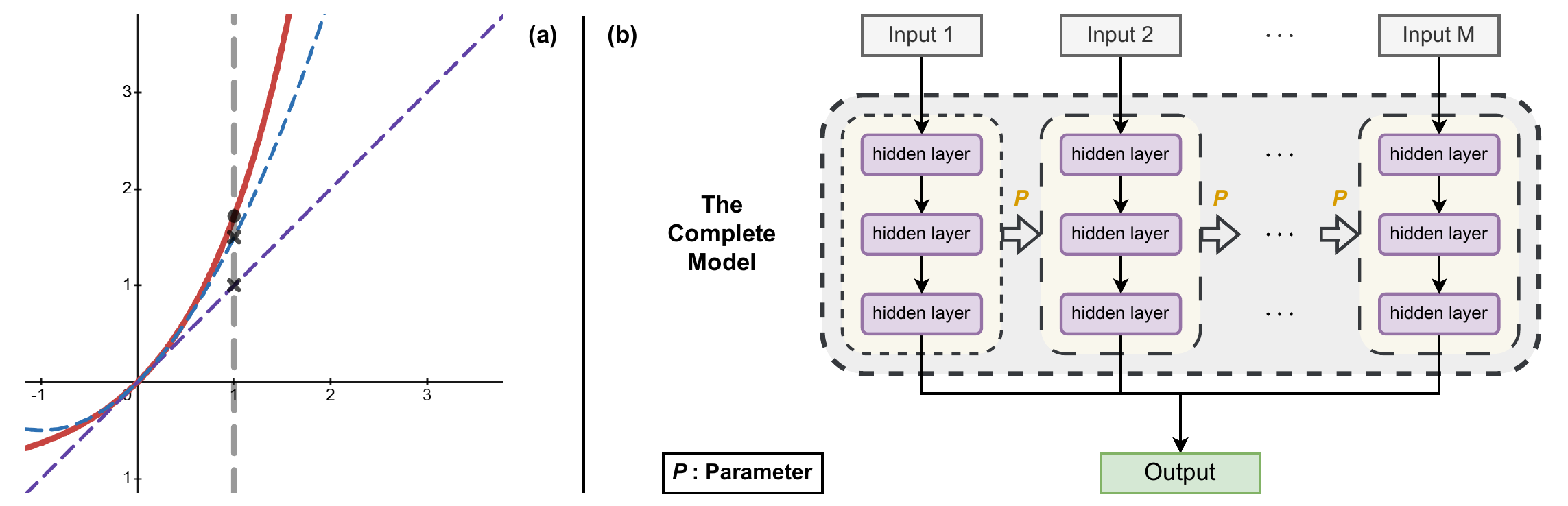}
    \caption{
    (a)  Conceptual analogy: Progressive temporal learning enhances solution accuracy similarly to how higher-order Taylor terms improve series approximation.
    (b) Network architecture: Progressive temporal expansion through sequential training, parameter freezing, and task-wise extension.
    }
    \label{fig:taylor}
\end{figure}

To facilitate the stability and efficiency, a two-phase optimization strategy is utilized. The first phase employs the Adam optimizer for initial training, followed by a fine-tuning phase with the L-BFGS optimizer. This hierarchical optimization aligns with our progressive refinement philosophy. The underlying mechanism can be conceptualized through a Taylor series analogy (Figure \ref{fig:taylor}): just as higher-order terms provide increasingly accurate function approximations, our model enhances prediction fidelity through gradual temporal expansion of the network architecture.

Our design is inspired by recent studies of temporal implicit PINNs, which highlight the impact of network structure on learning processes. Empirically, we observe that moderate temporal decomposition can improve model performance by localizing learning and reducing optimization difficulty. However, excessive partitioning may lead to instability and reduced generalization ability as error accumulation across subintervals becomes significant. These findings emphasize the need for balanced architectural choices that improve both expressiveness and robustness.

\begin{figure}[!htbp]
    \centering
    \includegraphics[width=0.8\linewidth]{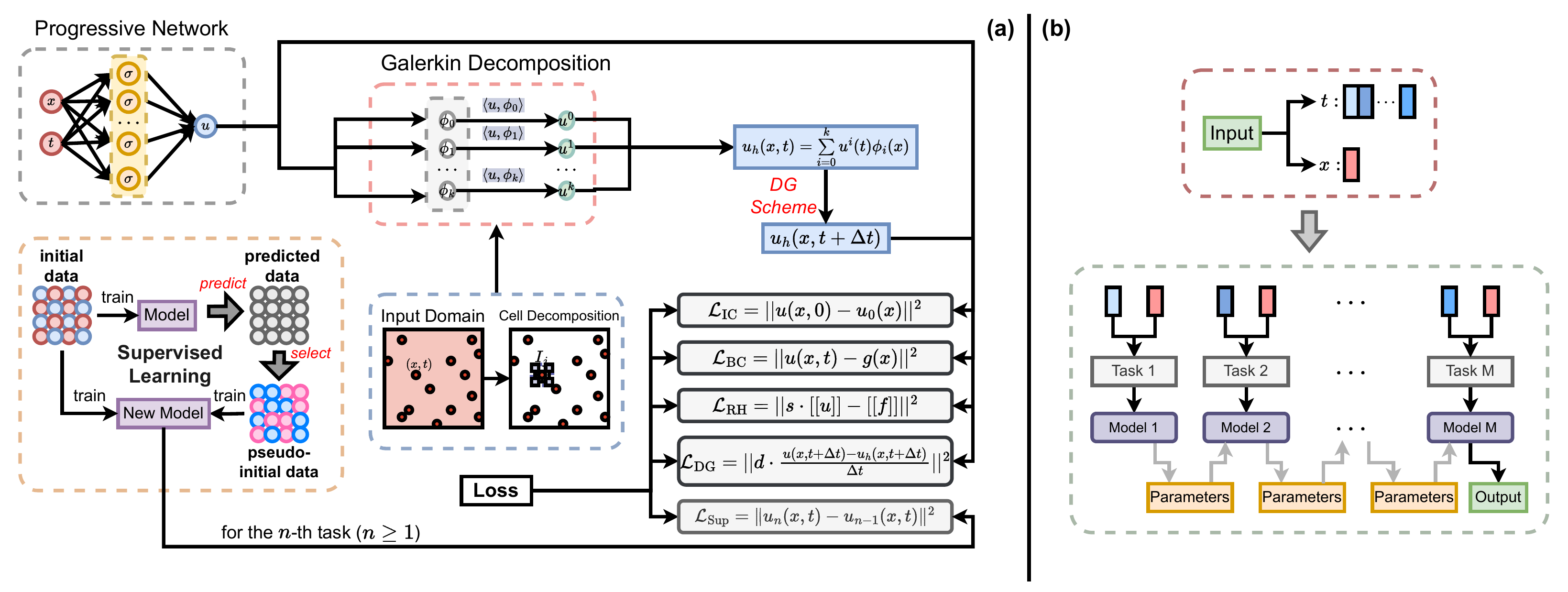}
    \caption{
    Overview of the training framework: 
    (a) Schematic of the total loss function incorporating physics-informed and supervised components. 
    (b) Illustration of the progressive network structure with task-wise inheritance and temporal expansion.
    }
    \label{framework}
\end{figure}

An overview of the proposed training pipeline, including the architecture of the progressive network and the composition of the total loss, is provided in Figure \ref{framework}. The complete objective function integrates several physical information and supervision components, denoted as
\begin{equation}
    \mathcal{L}= \omega_{\mathrm{IC}}\mathcal{L}_{\mathrm{IC}} + \omega_{\mathrm{bdy}}\mathcal{L}_{\mathrm{bdy}} + \omega_{\mathrm{DG}}\mathcal{L}_{\mathrm{DG}} + \omega_{\mathrm{RH}}\mathcal{L}_{\mathrm{RH}} + \omega_{\mathrm{sup}}\mathcal{L}_{\mathrm{sup}}.
\end{equation}

\subsection{Algorithmic Framework}

We now outline the complete algorithmic framework of the progressive learning strategy used to solve hyperbolic conservation laws over time $\mathcal{T}$. The model evolves by dividing the time domain into multiple subintervals, each of which is handled by specialized neural network components trained sequentially. This staged approach allows for flexible expansion of the model and progressive refinement over time.

\begin{algorithm}[H]
\caption{Progressive Neural Solver for Hyperbolic Conservation Laws}
\label{alg:model}
\KwIn{Number of time partitions $M$; spatial domain $x$; temporal domain $t \in \mathcal{T}$}
\KwOut{A trained neural network for multi-stage PDE approximation}

Partition $\mathcal{T}$ into $M$ uniform subintervals: $\mathcal{T}^{(k)} = [t_{k-1}, t_k],\,1 \leq k \leq M$ \;
Initialize base model parameters $\theta^{(1)}$ \;

\For{$k = 1$ \KwTo $M$}{
    \tcc{Temporal expansion with parameter freezing}
    Freeze parameters $\{\theta^{(1)}, \dots, \theta^{(k-1)}\}$ \;
    Append new trainable subnetwork to obtain $\theta^{(k)}$ \;

    Construct training dataset $\mathcal{D}^{(k)}$ on $\mathcal{T}^{(k)}$ \;

    \While{training not converged}{
        Evaluate physics-informed loss $\mathcal{L}_{\text{DG}}^{(k)}$ on $\mathcal{D}^{(k)}$ \;
        Evaluate initial/interface condition loss $\mathcal{L}_{\text{IC}}^{(k)}$ \;
        Evaluate boundary condition loss $\mathcal{L}_{\text{bdy}}^{(k)}$ \;
        Evaluate Rankine-Hugoniot loss $\mathcal{L}_{\text{RH}}^{(k)}$ \;

        \If{$k > 1$}{
            Construct pseudo-labeled dataset from previous task predictions \;
            Evaluate pseudo-supervised loss $\mathcal{L}_{\text{sup}}^{(k)}$ \;
        }

        Aggregate total loss: \\
        \quad $\mathcal{L}^{(k)} = \mathcal{L}_{\text{DG}}^{(k)} + \mathcal{L}_{\text{IC}}^{(k)} + \mathcal{L}_{\text{bdy}}^{(k)} + \mathcal{L}_{\text{RH}}^{(k)} + \mathrm{1}_{\{k>1\}} \cdot \mathcal{L}_{\text{sup}}^{(k)}$ \;

        Update trainable parameters: \\
        \quad $\theta^{(k)} \leftarrow \theta^{(k)} - \eta \nabla_{\theta^{(k)}} \mathcal{L}^{(k)}$ \;
    }
    Store interface predictions for initialization of task $k+1$ \;
}
\end{algorithm}

\section{Theoretical Analysis}
\label{sec:theory}
We analyze the convergence of our model in two steps: (i) the neural network (NN) approximation to the discontinuous Galerkin (DG) solution, and (ii) the DG approximation to the exact solution. The details of proofs will be provided in the appendix \ref{app:proof}.

\subsection{NN Approximation to DG Solution}

\begin{theorem}[NN-to-DG Error Bound]
\label{thm:nn_dg}
Let $\mathcal{V}_h$ be the piecewise polynomial DG space on a mesh of size $h$ and polynomial degree $p \ge 0$. 
For a given time subinterval $\mathcal{T}^{(k)} = [t_{k-1}, t_k]$, assume:
\begin{enumerate}[(i)]
    \item (\textbf{Representability}) There exists $\delta_{\mathrm{approx}}^{(k)} \ge 0$ such that
    $$
        \inf_{\theta \in \Theta_{\mathrm{NN}}} \| u_{\theta} - u_{\mathrm{DG}}^{(k)} \|_{L^2(\mathcal{T}^{(k)} \times \Omega)} \le \delta_{\mathrm{approx}}^{(k)},
    $$
    where $\Theta_{\mathrm{NN}}$ denotes the admissible parameter set of the NN subnetwork for task $k$.
    \item (\textbf{Optimization}) The training procedure produces parameters $\theta^{(k)}$ satisfying
    $$
        \mathcal{L}^{(k)}(\theta^{(k)}) - \min_{\theta} \mathcal{L}^{(k)}(\theta) \le \epsilon_{\mathrm{opt}}^{(k)},
    $$
    with $\mathcal{L}^{(k)}$ the total loss on $\mathcal{T}^{(k)}$ and $\epsilon_{\mathrm{opt}}^{(k)}$ the optimization error.
    \item (\textbf{Coercivity of loss}) The loss $\mathcal{L}^{(k)}$ is coercive with respect to the $L^2$-error, i.e.,
    $$
        \| u_{\theta} - u_{\mathrm{DG}}^{(k)} \|_{L^2(\mathcal{T}^{(k)} \times \Omega)}
        \le C_{\mathrm{coer}}^{(k)} \, \mathcal{L}^{(k)}(\theta)^{1/2},
    $$
    for some constant $C_{\mathrm{coer}}^{(k)}$ depending only on the DG discretization and the PDE coefficients.
\end{enumerate}
Then the trained NN prediction $u_{\mathrm{NN}}^{(k)} := u_{\theta^{(k)}}$ satisfies
\begin{equation}
\label{eq:nn_dg_bound}
    \| u_{\mathrm{NN}}^{(k)} - u_{\mathrm{DG}}^{(k)} \|_{L^2(\mathcal{T}^{(k)} \times \Omega)}
    \le C_{\mathrm{coer}}^{(k)} \, \epsilon_{\mathrm{opt}}^{(k)\,1/2} + \delta_{\mathrm{approx}}^{(k)}.
\end{equation}
\end{theorem}

\subsection{DG Approximation to Exact Solution}

\begin{theorem}[DG-to-Exact Error Bound]
\label{thm:dg_exact}
Let $u_{\mathrm{exact}}$ be the entropy solution of a hyperbolic conservation law on $\Omega \times [0,T]$ with piecewise smooth initial data. 
Let $u_{\mathrm{DG}}$ be the DG solution with polynomial degree $p \ge 0$ and mesh size $h$.
Assume a stable numerical flux and a CFL condition ensuring $L^2$-stability. Then:
\begin{enumerate}[(1)]
    \item (\textbf{Smooth regions}) If $u_{\mathrm{exact}}$ is smooth on a subdomain $D \subset \Omega \times [0,T]$, there exists a constant $C_{\mathrm{smooth}}$ such that
    $$
        \| u_{\mathrm{DG}} - u_{\mathrm{exact}} \|_{L^2(D)} \le C_{\mathrm{smooth}} \, h^{p+1}.
    $$
    \item (\textbf{Discontinuous regions}) In the vicinity of shocks or contact discontinuities, there exists a constant $C_{\mathrm{disc}}$ such that
    $$
        \| u_{\mathrm{DG}} - u_{\mathrm{exact}} \|_{L^2(\Omega \times [0,T])} \le C_{\mathrm{disc}} \, h^{1/2},
    $$
    where the Rankine-Hugoniot loss $\mathcal{L}_{\mathrm{RH}}$ enforces correct shock speeds up to an additional term $\delta_{\mathrm{RH}}$:
    $$
        |\lambda_{\mathrm{DG}} - \lambda_{\mathrm{exact}}| \le \delta_{\mathrm{RH}}.
    $$
\end{enumerate}
The constants $C_{\mathrm{smooth}}$ and $C_{\mathrm{disc}}$ are independent of $h$.
\end{theorem}

\subsection{Global Convergence with Progressive Training}

\begin{corollary}[Total Error Bound]
\label{cor:total_error}
Let $S_k \ge 1$ be the stability factor of the DG evolution operator on $\mathcal{T}^{(k)}$, defined such that an $L^2$-perturbation of size $\eta$ at $t_{k-1}$ produces at most $S_k \eta$ at $t_k$. 
Under Theorems \ref{thm:nn_dg} and \ref{thm:dg_exact}, the overall error of the progressive NN solver satisfies
\begin{equation}
\label{eq:total_error_bound}
\begin{aligned}
    \| u_{\mathrm{NN}} - u_{\mathrm{exact}} \|_{L^2(\Omega \times [0,T])} 
    &\le \sum_{k=1}^M \Big( \big( \prod_{j=k+1}^M S_j \big) \,
        \big[ C_{\mathrm{coer}}^{(k)} \epsilon_{\mathrm{opt}}^{(k)\,1/2} + \delta_{\mathrm{approx}}^{(k)} \big] \Big) \\
    &\quad + C_{\mathrm{smooth}} \, h^{p+1} + C_{\mathrm{disc}} \, h^{1/2} + \delta_{\mathrm{RH}}.
\end{aligned}
\end{equation}
If $S_k \equiv 1$ (ideal $L^2$-stability), the product term disappears and the error accumulates linearly in $M$.
\end{corollary}
The above theoretical results rigorously characterize the error behavior of our progressive neural solver for hyperbolic conservation laws.  
Specifically, Theorem \ref{thm:nn_dg} quantifies the approximation quality of the neural network with respect to the underlying DG solution, explicitly separating the contributions from network expressivity and optimization accuracy.  
Theorem \ref{thm:dg_exact} recalls classical DG convergence rates toward the exact entropy solution, accounting for both smooth solution regions and discontinuities, while highlighting the role of the Rankine-Hugoniot loss in ensuring physically consistent shock capturing.  
Finally, Corollary \ref{cor:total_error} integrates these results, emphasizing the importance of stability factors governing temporal error propagation across successive subintervals in the progressive training scheme.  

Collectively, these bounds demonstrate that the overall error can be effectively controlled by refining the mesh size, improving network training, and properly managing error propagation through the progressive architecture.  
This theoretical foundation justifies the efficacy of our approach and motivates the forthcoming numerical experiments designed to validate these convergence properties and illustrate the solver’s capability in resolving complex solution features such as shocks and rarefactions.

\section{Numerical Experiments}
\label{sec:experiments}
In this section, we conduct numerical experiments to evaluate the proposed hybrid framework in both scalar and system-level hyperbolic conservation laws. The experiments are designed to assess accuracy, stability, and shock-capturing capability under varying initial and mesh conditions.

\subsection{Experimental Design}
For reproducibility, we fix the random seed to 7 across all experiments. Our unified neural architecture handles both Burgers and Euler equations through an expandable framework. The base network consists of 4 hidden layers (64 neurons each), with 4 additional layers appended during each progressive training phase. This design enables incremental temporal learning while maintaining stability of previously acquired solutions.

Key implementation details:
\begin{itemize}
    \item Loss weighting: $\omega_{\mathrm{IC}} = 10$ for initial conditions, 1 for other terms.
    \item Activation: Standard $\tanh$ (or $\tanh(x)+1$ when positivity constraints apply).
    \item Spatial resolutions tested: 32, 64, 128, and 256 cells.
    \item Progressive tasks: 1 to 4 phases (studying temporal decomposition effects).
\end{itemize}
We avoid sigmoid activations due to their inferior accuracy from gradient vanishing issues. The progressive expansion systematically increases model capacity while preserving learned dynamics through layer-wise freezing and expansion.

\subsection{Burgers Equation}
We begin with the scalar nonlinear Burgers equation:
\begin{equation}
    \begin{cases}
        u_t + \left(\frac{1}{2}u^2\right)_x = 0, \\
        u(x,0) = u_0(x).
    \end{cases}
    \label{burgerseq}
\end{equation}

The model is tested on three representative initial conditions:
\begin{itemize}
    \item \textbf{Case 1 (Gaussian initial condition):} A smooth single-pulse profile defined as
    \begin{equation}
        u_0(x) = e^{-x^2}, \quad x \in [-2,2].
    \end{equation}

    \item \textbf{Case 2 (Riemann initial condition):} A discontinuous step function given by
    \begin{equation}
        u_0(x) = 
        \begin{cases}
            1, & 0 < x < \frac{1}{2}, \\
            0, & \frac{1}{2} < x < 1.
        \end{cases}
    \end{equation}

    \item \textbf{Case 3 (Sine wave initial condition):} A smooth periodic profile expressed as
    \begin{equation}
        u_0(x) = \sin(\pi x), \quad x \in [0,2].
    \end{equation}
\end{itemize}

With the spatial resolution fixed at 256 cells, we conducted a systematic evaluation of different task configurations. The comparative results in Table \ref{tab:burgers-task} reveal that training without temporal decomposition (single-task approach) produces substantially larger errors compared to progressive training strategies. The model accuracy shows consistent improvement when increasing the number of tasks from one to three, demonstrating the effectiveness of progressive temporal decomposition in stabilizing the optimization process and more accurately capturing solution dynamics. However, extending to four tasks leads to marginally increased errors in some test cases, suggesting that excessive temporal segmentation may reach a point of diminishing returns where additional decomposition could potentially introduce redundant parameters or compromise the model's generalization capability.

\begin{figure}[!htbp]
  \centering

  \begin{subfigure}[b]{0.24\textwidth}
    \includegraphics[width=\textwidth]{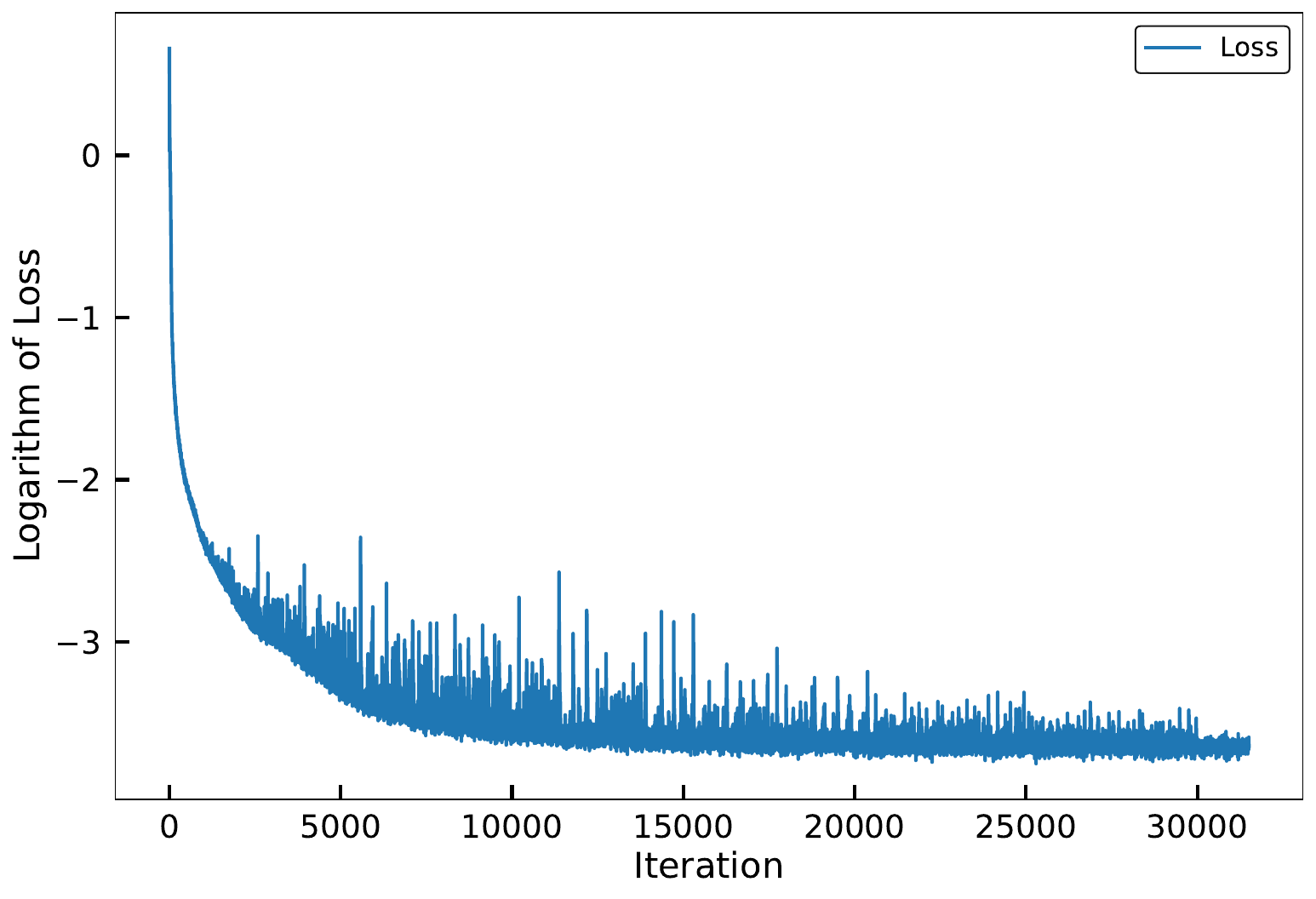}
    \caption{1 task: Loss}
  \end{subfigure}
  \begin{subfigure}[b]{0.24\textwidth}
    \includegraphics[width=\textwidth]{error/loss/gauss_1.pdf}
    \caption{2 tasks: Loss}
  \end{subfigure}
  \begin{subfigure}[b]{0.24\textwidth}
    \includegraphics[width=\textwidth]{error/loss/gauss_1.pdf}
    \caption{3 tasks: Loss}
  \end{subfigure}
  \begin{subfigure}[b]{0.24\textwidth}
    \includegraphics[width=\textwidth]{error/loss/gauss_1.pdf}
    \caption{4 tasks: Loss}
  \end{subfigure}

  \begin{subfigure}[b]{0.24\textwidth}
    \includegraphics[width=\textwidth]{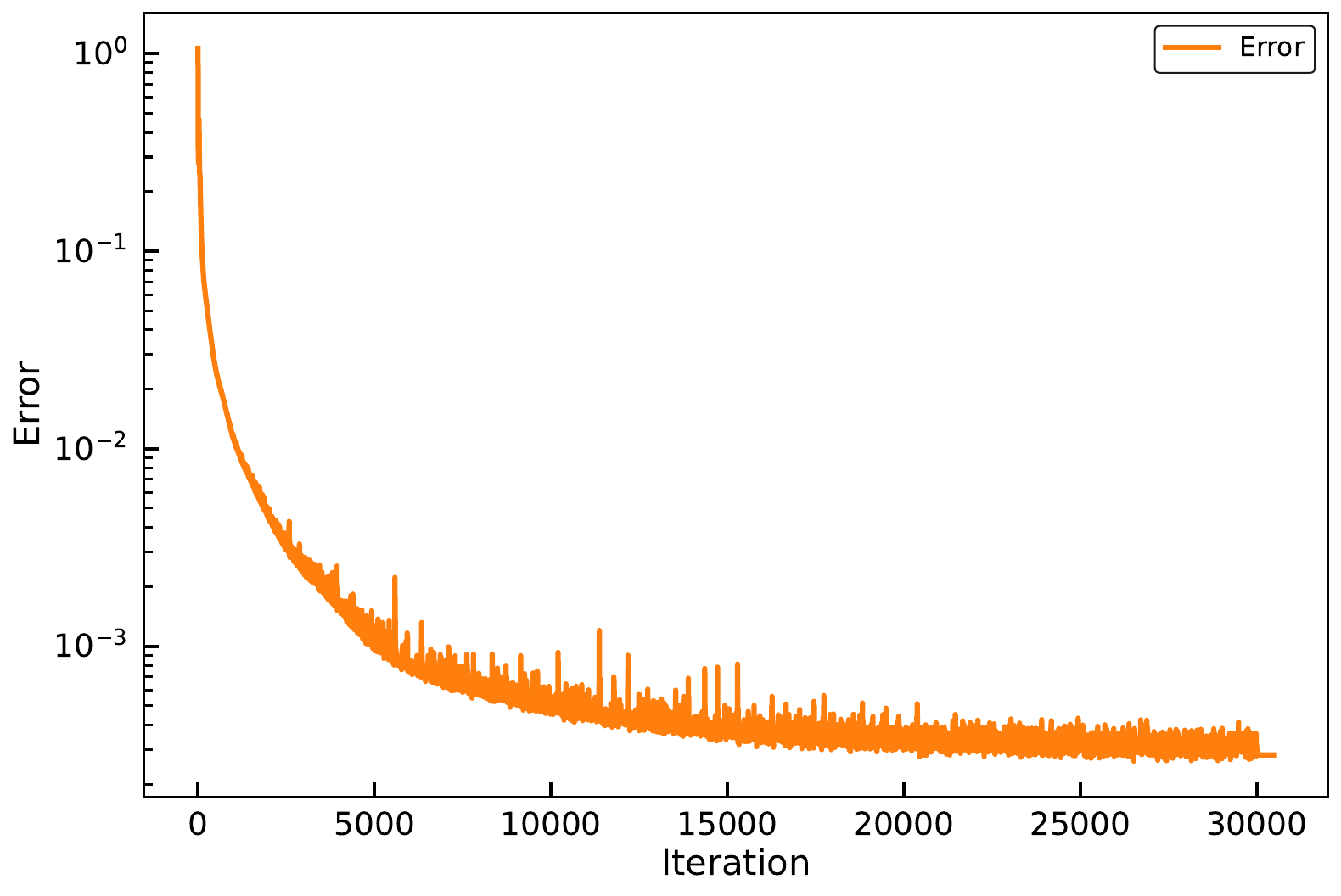}
    \caption{1 task: Error}
  \end{subfigure}
  \begin{subfigure}[b]{0.24\textwidth}
    \includegraphics[width=\textwidth]{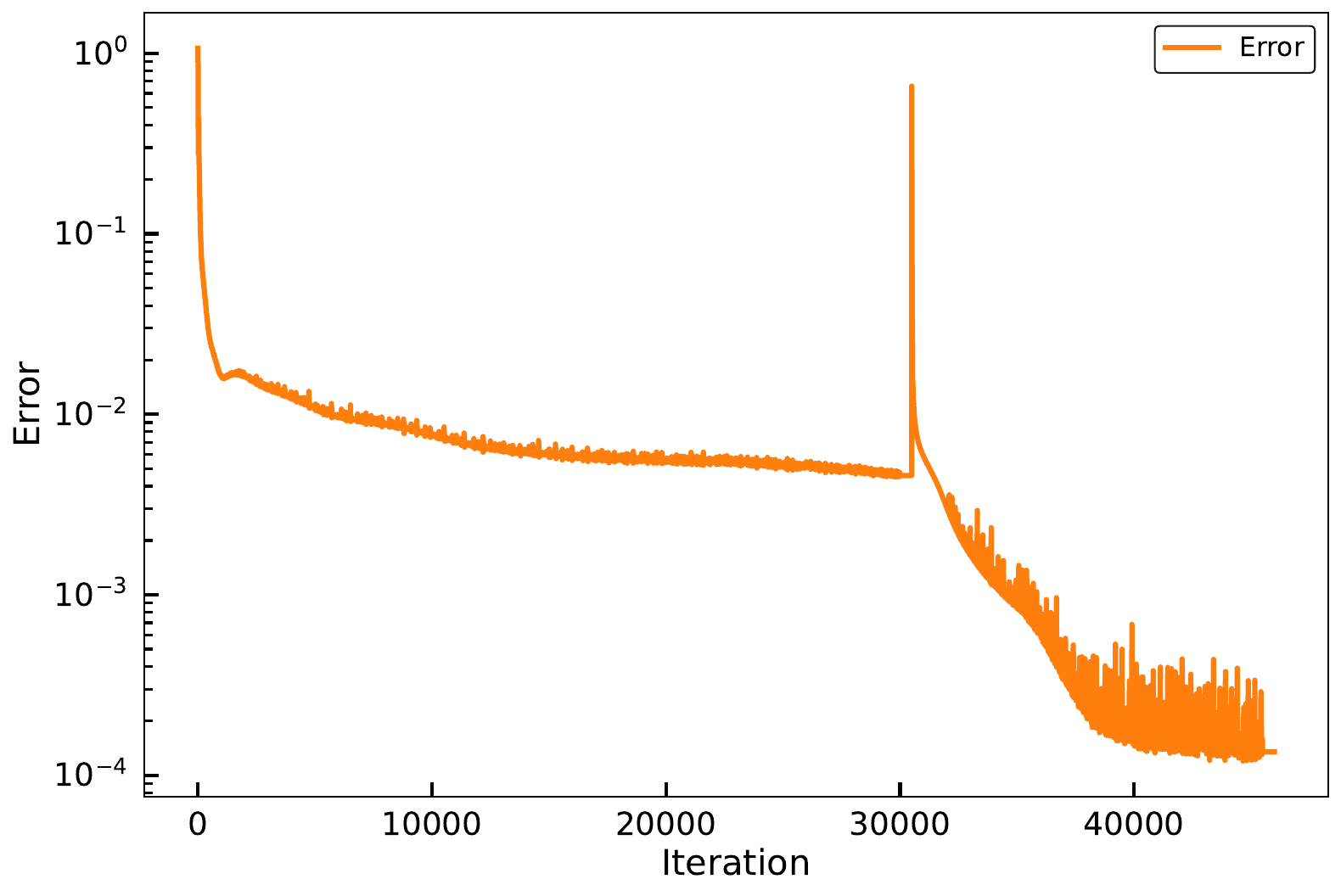}
    \caption{2 tasks: Error}
  \end{subfigure}
  \begin{subfigure}[b]{0.24\textwidth}
    \includegraphics[width=\textwidth]{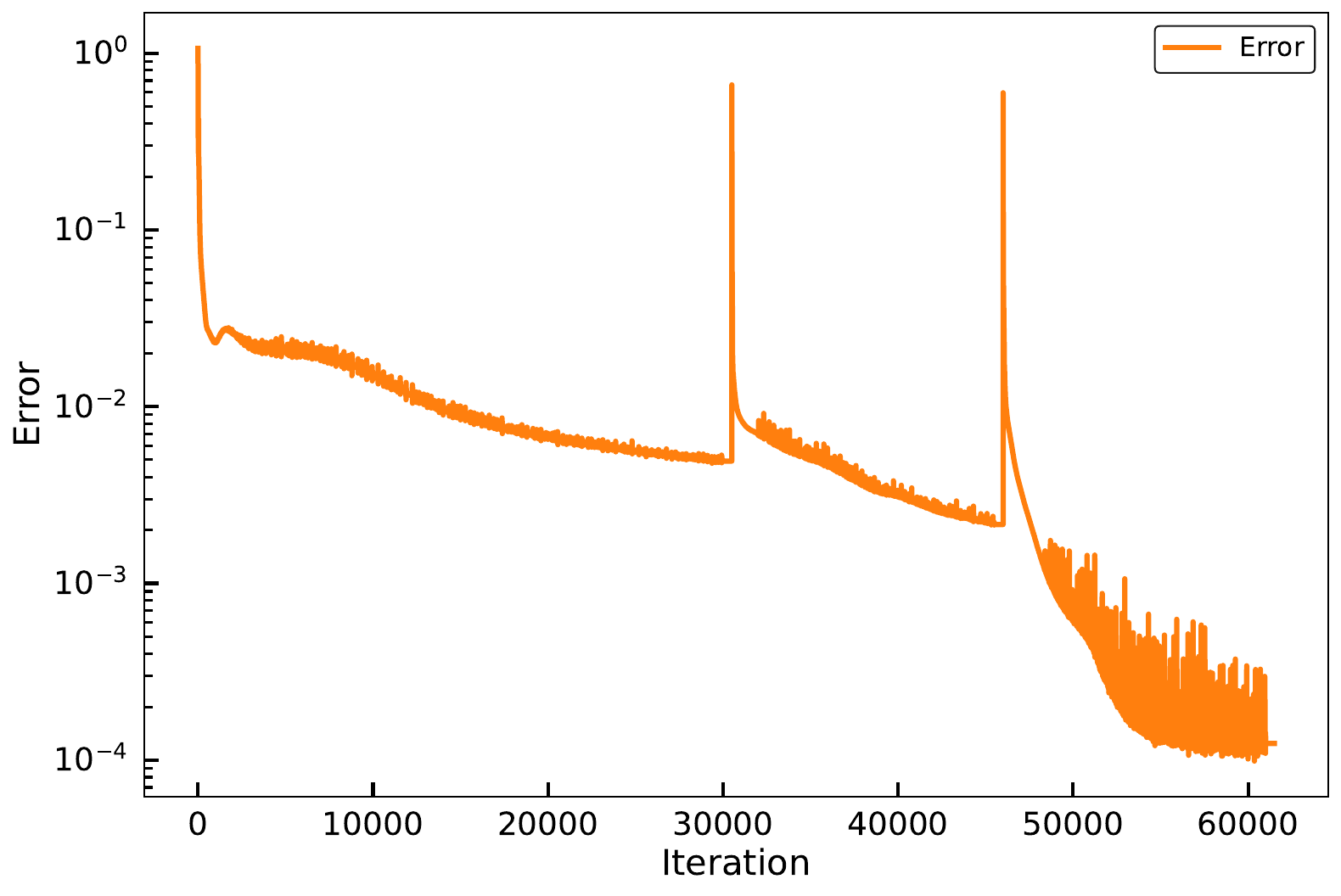}
    \caption{3 tasks: Error}
  \end{subfigure}
  \begin{subfigure}[b]{0.24\textwidth}
    \includegraphics[width=\textwidth]{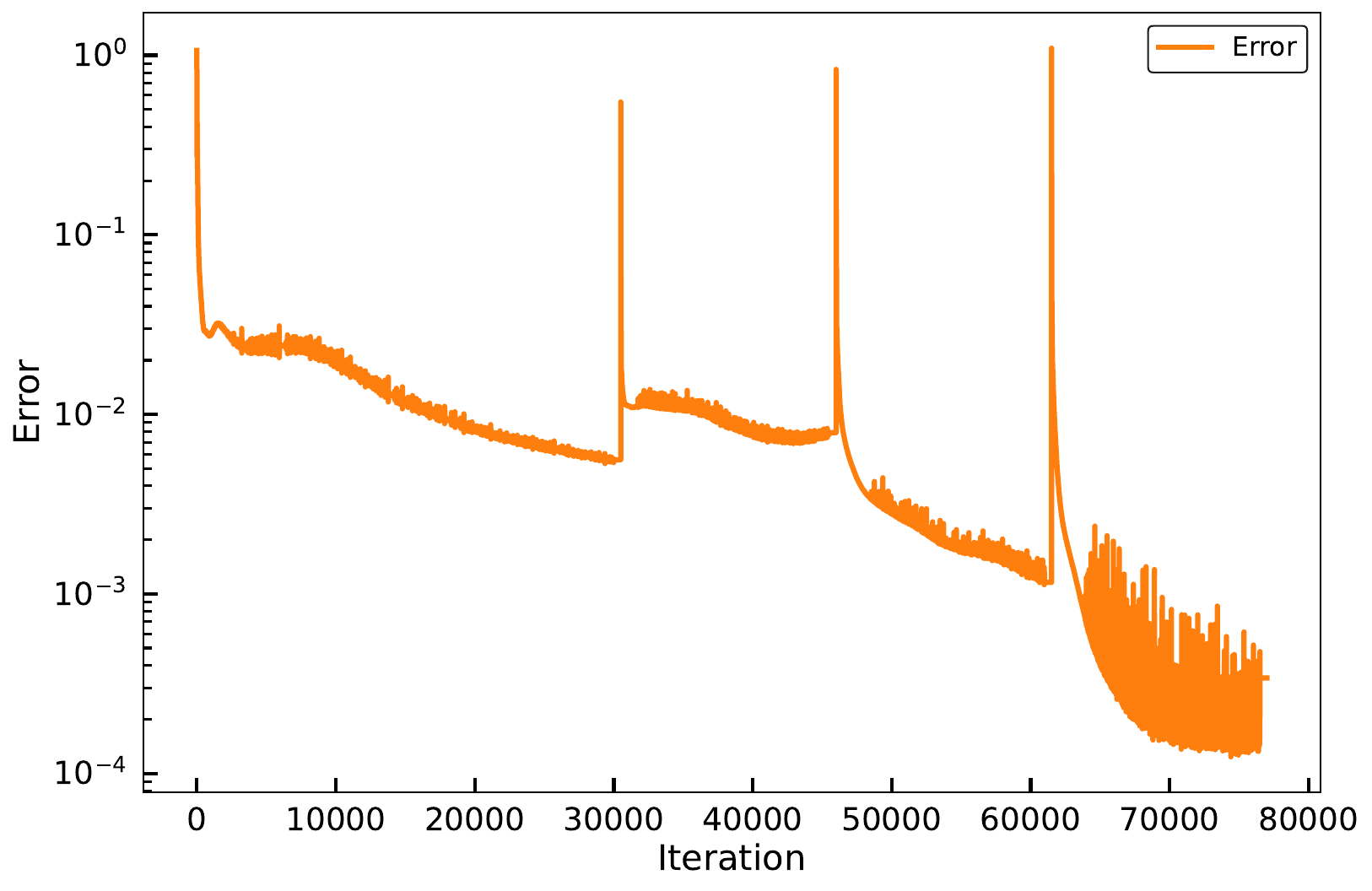}
    \caption{4 tasks: Error}
  \end{subfigure}

  \caption{Evolution of the training loss (top) and prediction error (bottom) across sequential tasks for the Gaussian initial condition. The results demonstrate improved performance through progressive temporal decomposition.}
  \label{fig:loss_error_tasks}
\end{figure}

\begin{table}[htbp]
  \centering
  \caption{Relative $L^2$ errors of the Burgers equation for different initial conditions under varying numbers of temporal tasks (256 cells).}
  \label{tab:burgers-task}
  \sisetup{scientific-notation = true}
  \begin{tabular}{lcccc}
    \toprule
    \textbf{Initial Condition} & \textbf{1 Task} & \textbf{2 Tasks} & \textbf{3 Tasks} & \textbf{4 Tasks} \\
    \midrule
    Gaussian & \num{9.89e-3} & \num{6.26e-3} & \textbf{\num{6.16e-3}} & \num{1.52e-2} \\
    Riemann  & \num{2.90e-2} & \textbf{\num{2.33e-2}} & \num{2.38e-2} & \num{3.03e-2} \\
    Sine     & \num{1.45e-1} & \textbf{\num{6.41e-2}} & \num{6.58e-2} & \num{6.45e-2} \\
    \bottomrule
  \end{tabular}
\end{table}
To validate these findings, we examine the training dynamics through loss and prediction error trajectories, illustrated in Figure \ref{fig:loss_error_tasks} using the Gaussian initial condition case. The visualization reveals two key patterns: first, a consistent error reduction across successive training tasks, confirming the effectiveness of our progressive learning approach; second, an evolving loss landscape between tasks, attributable to both network expansion and the incorporation of pseudo-label supervision. These observations quantitatively demonstrate how our method achieves progressively refined solutions through temporal decomposition.

\begin{figure}[!htbp]
  \centering
  \begin{subfigure}[b]{\textwidth}
    \centering
    \begin{minipage}[b]{0.31\textwidth}
      \centering
      \includegraphics[width=\linewidth]{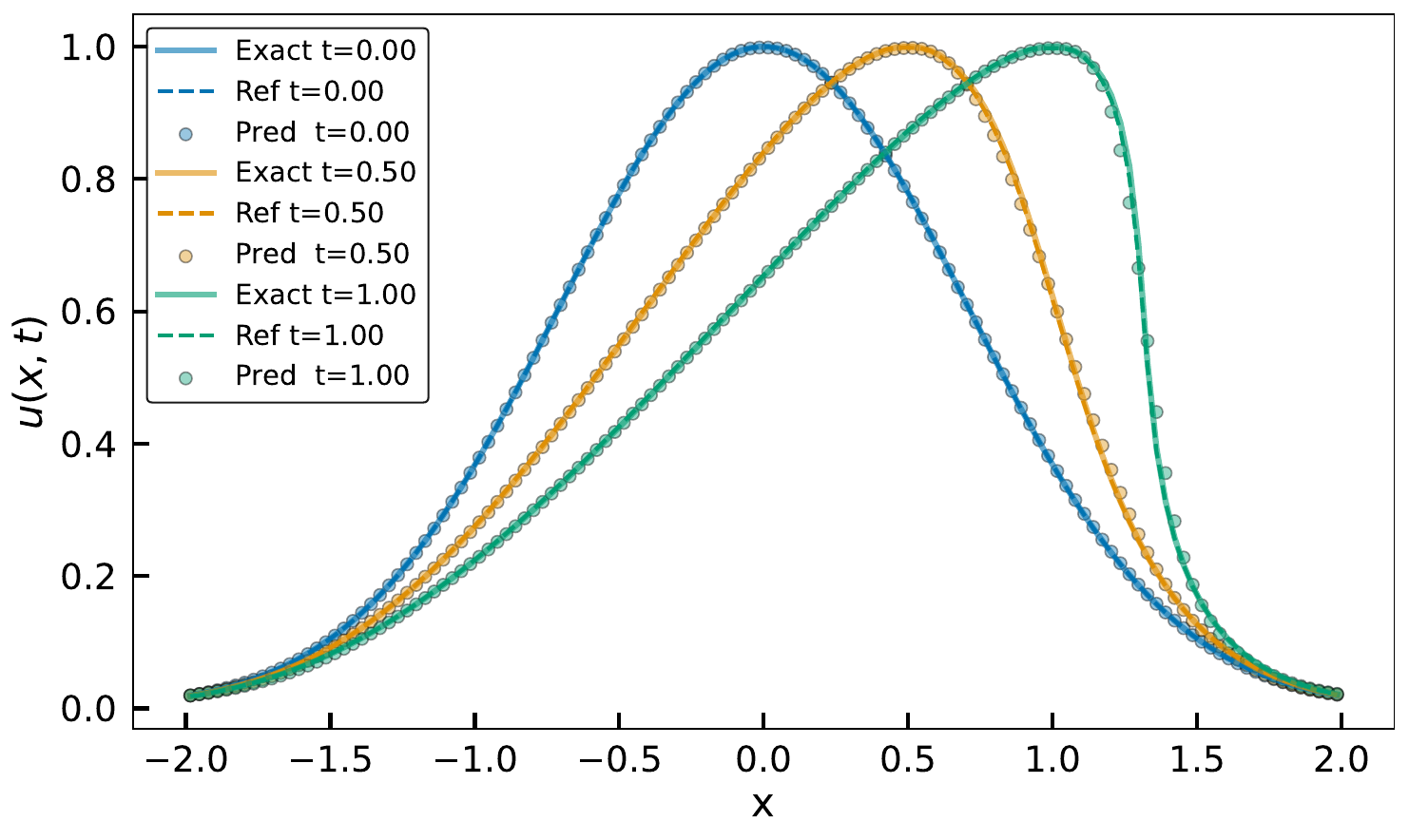}
      \subcaption{Gauss: 1 task}
    \end{minipage}
    \hspace{0.02\textwidth}
    \begin{minipage}[b]{0.31\textwidth}
      \centering
      \includegraphics[width=\linewidth]{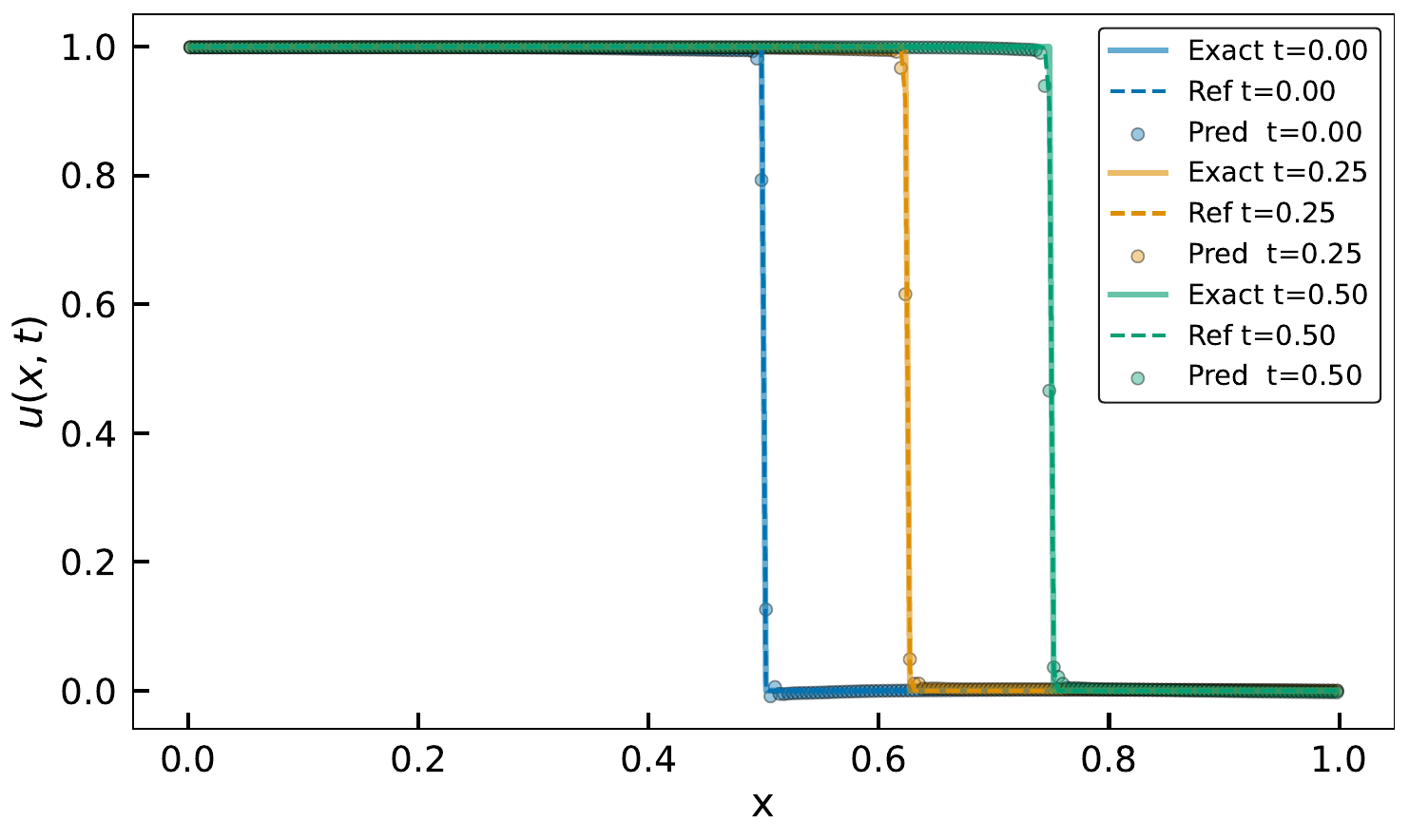}
      \subcaption{Riemann: 1 task}
    \end{minipage}
    \hspace{0.02\textwidth}
    \begin{minipage}[b]{0.31\textwidth}
      \centering
      \includegraphics[width=\linewidth]{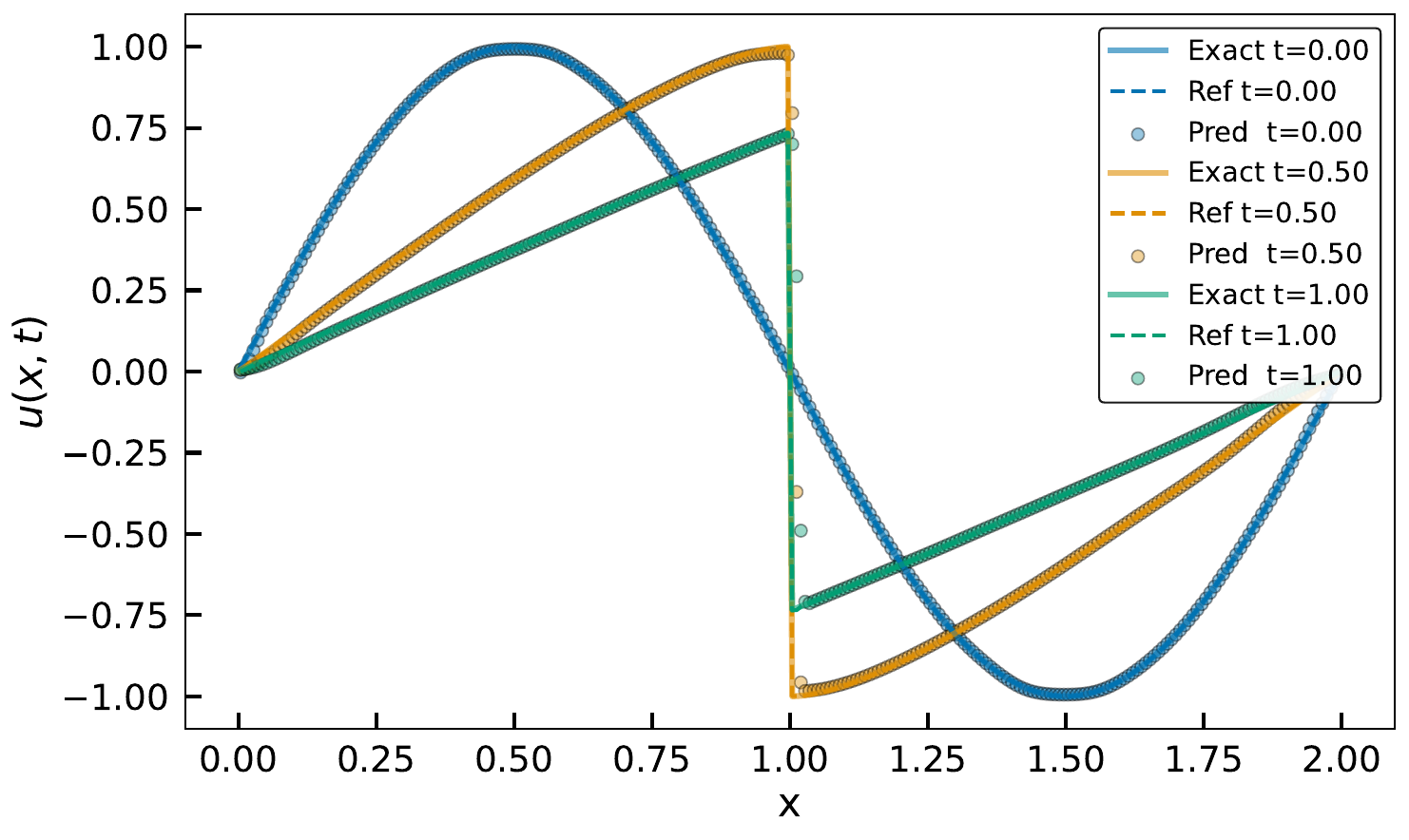}
      \subcaption{Sine: 1 task}
    \end{minipage}
  \end{subfigure}

  \begin{subfigure}[b]{\textwidth}
    \centering
    \begin{minipage}[b]{0.31\textwidth}
      \centering
      \includegraphics[width=\linewidth]{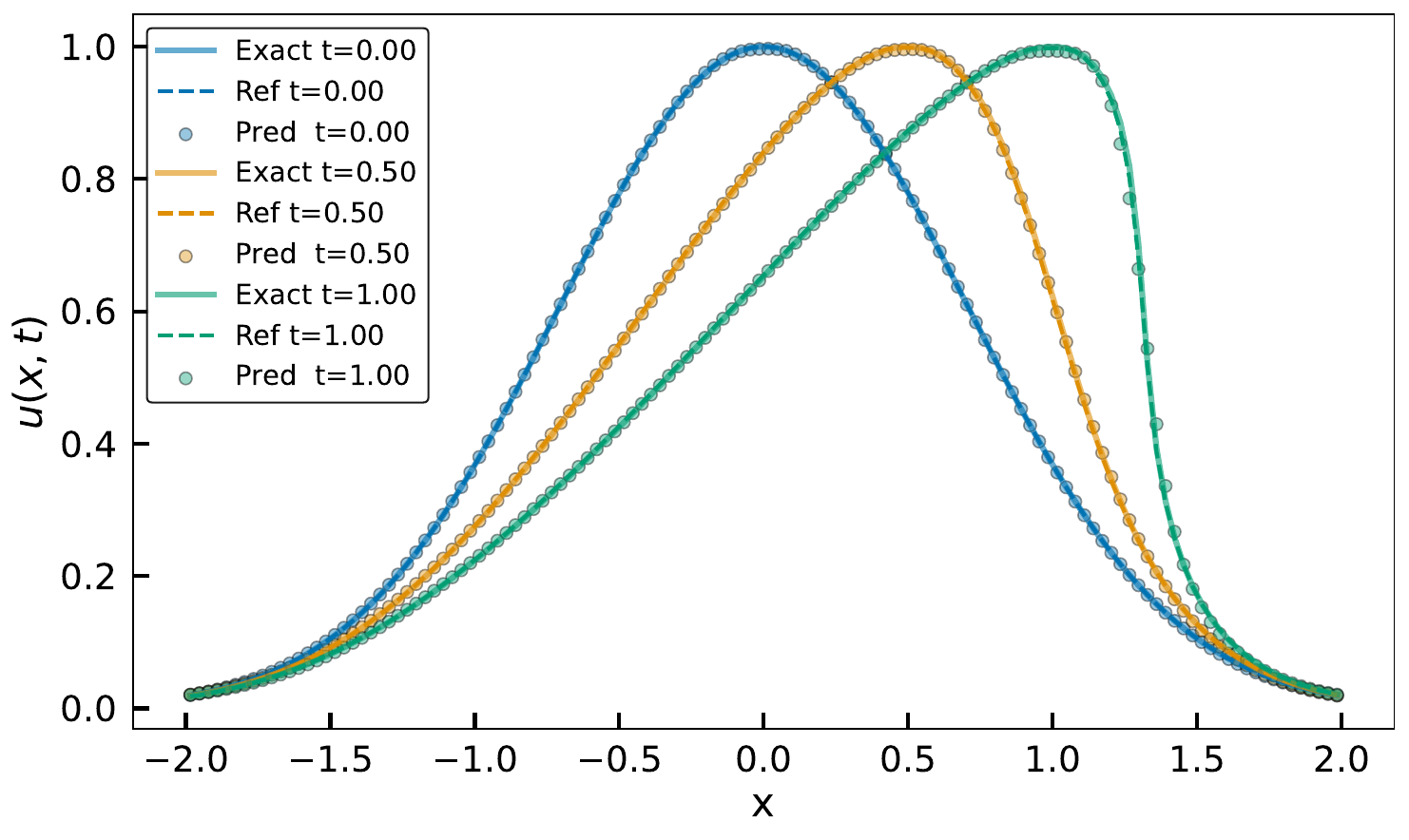}
      \subcaption{Gauss: 2 tasks}
    \end{minipage}
    \hspace{0.02\textwidth}
    \begin{minipage}[b]{0.31\textwidth}
      \centering
      \includegraphics[width=\linewidth]{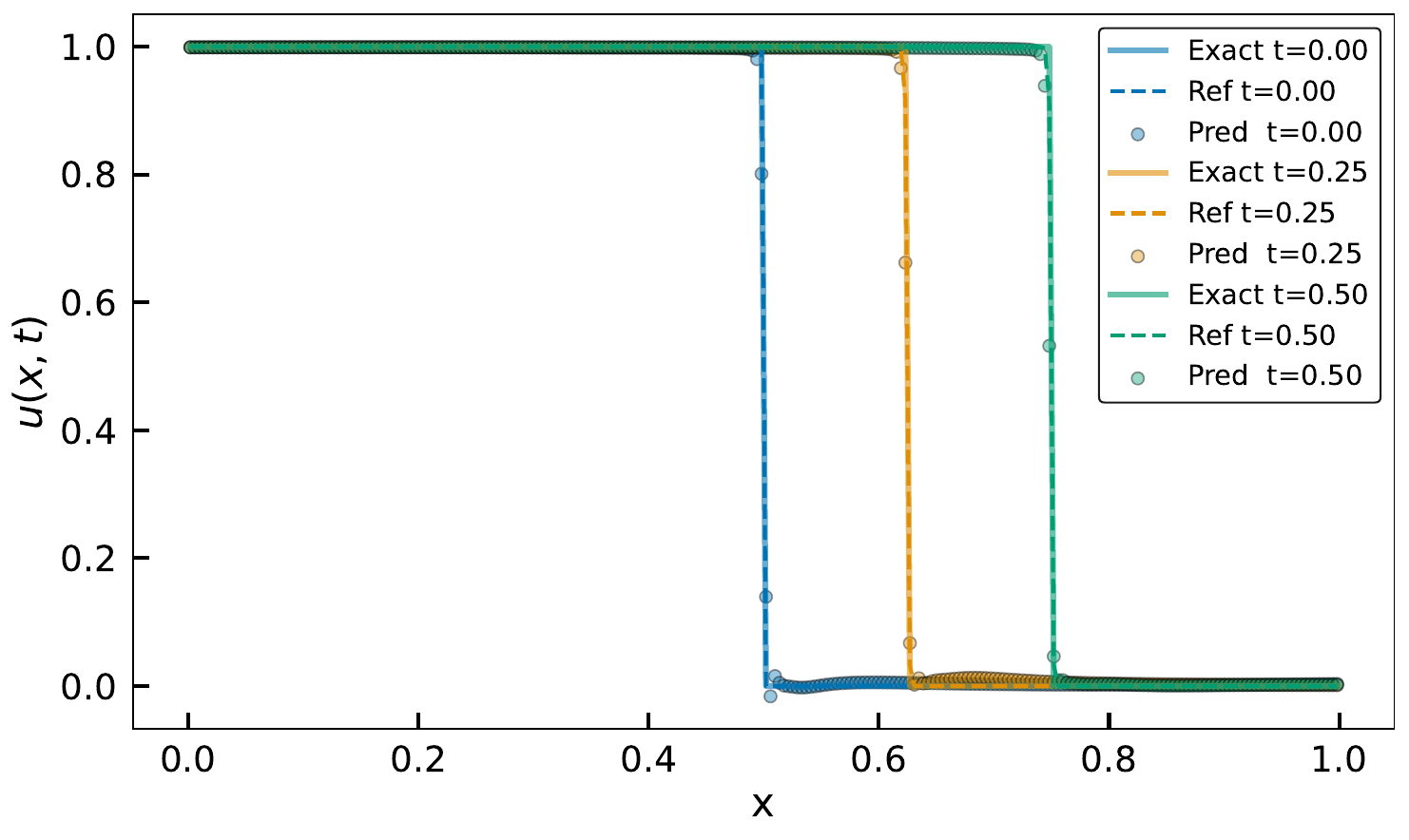}
      \subcaption{Riemann: 2 tasks}
    \end{minipage}
    \hspace{0.02\textwidth}
    \begin{minipage}[b]{0.31\textwidth}
      \centering
      \includegraphics[width=\linewidth]{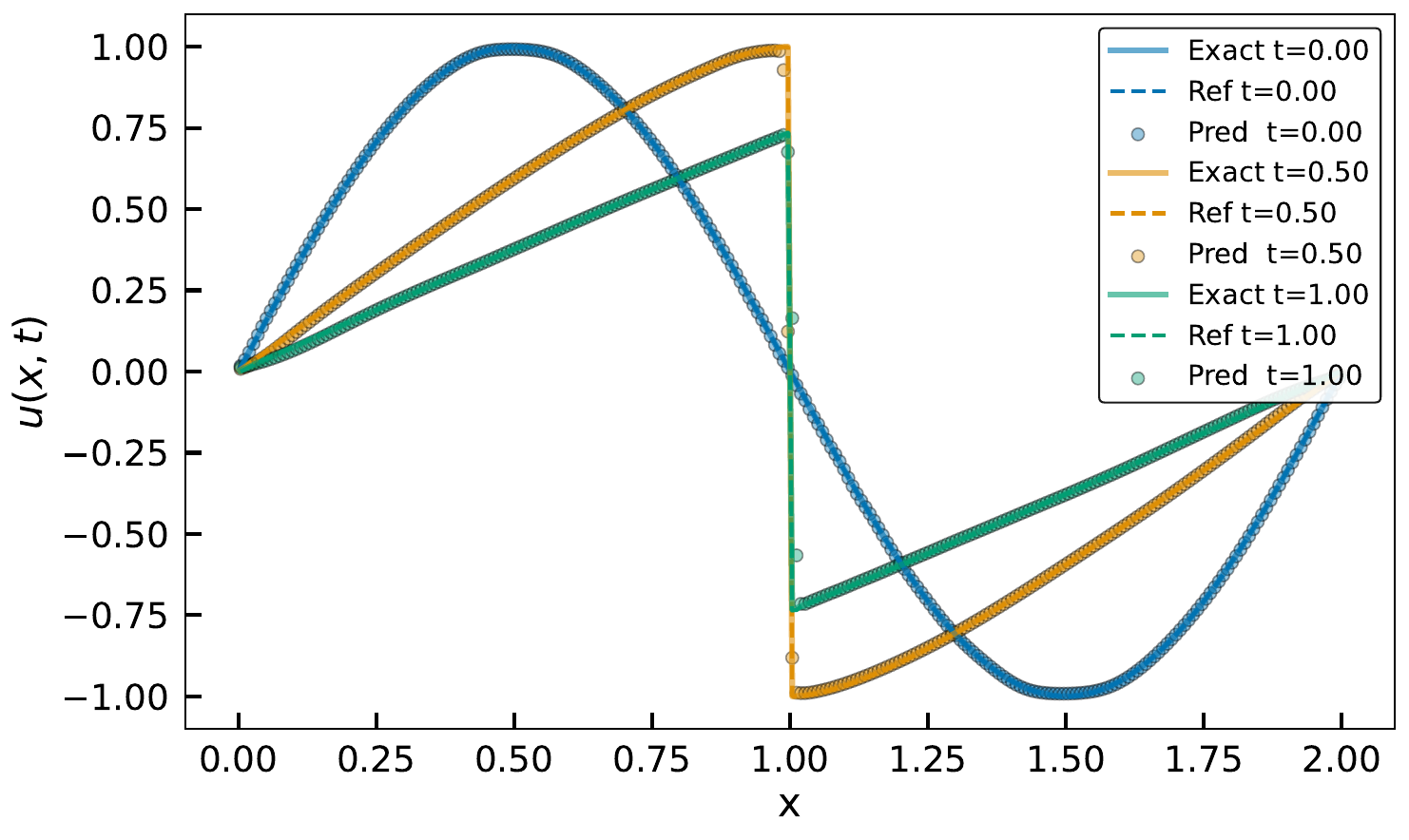}
      \subcaption{Sine: 2 tasks}
    \end{minipage}
  \end{subfigure}

  \begin{subfigure}[b]{\textwidth}
    \centering
    \begin{minipage}[b]{0.31\textwidth}
      \centering
      \includegraphics[width=\linewidth]{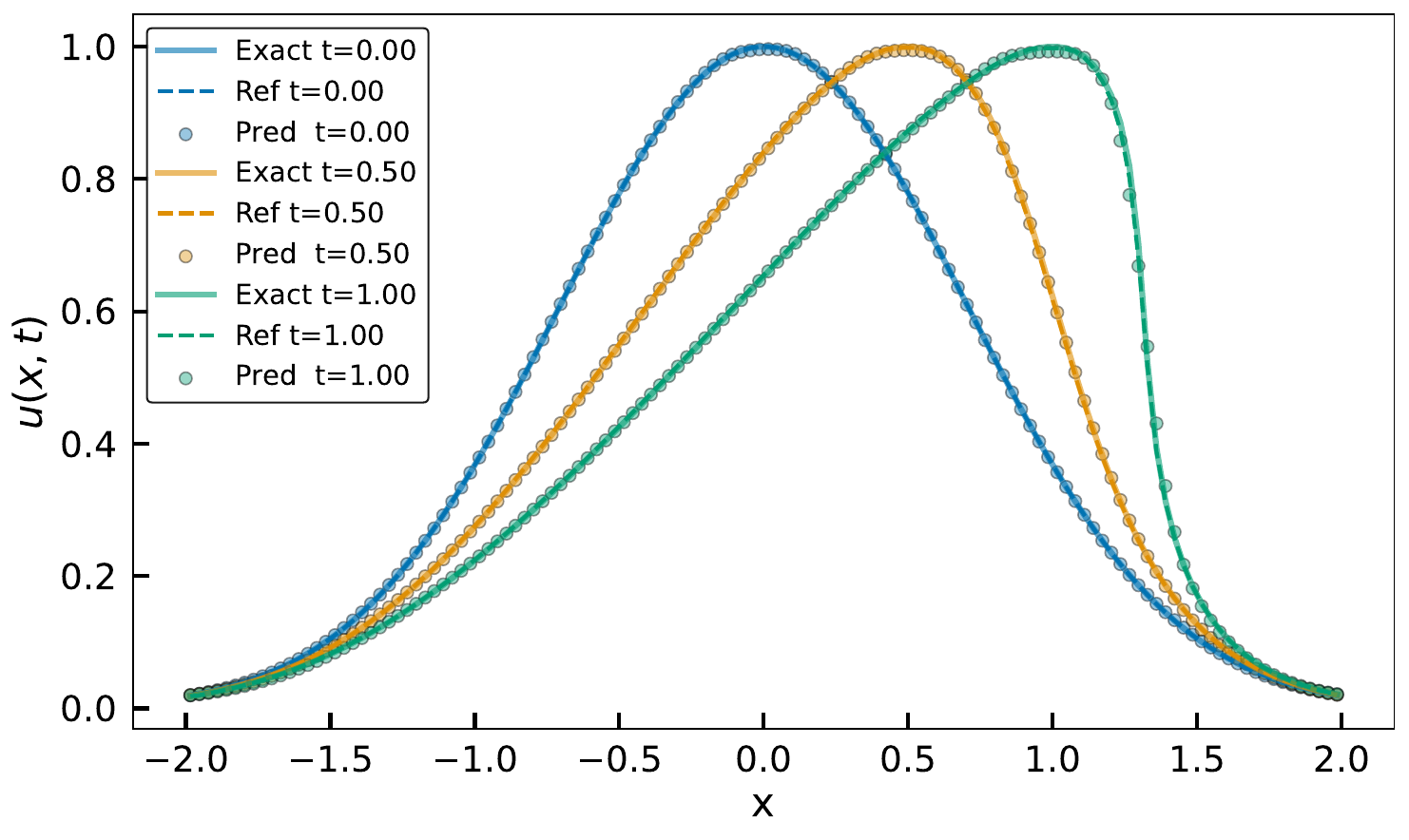}
      \subcaption{Gauss: 3 tasks}
    \end{minipage}
    \hspace{0.02\textwidth}
    \begin{minipage}[b]{0.31\textwidth}
      \centering
      \includegraphics[width=\linewidth]{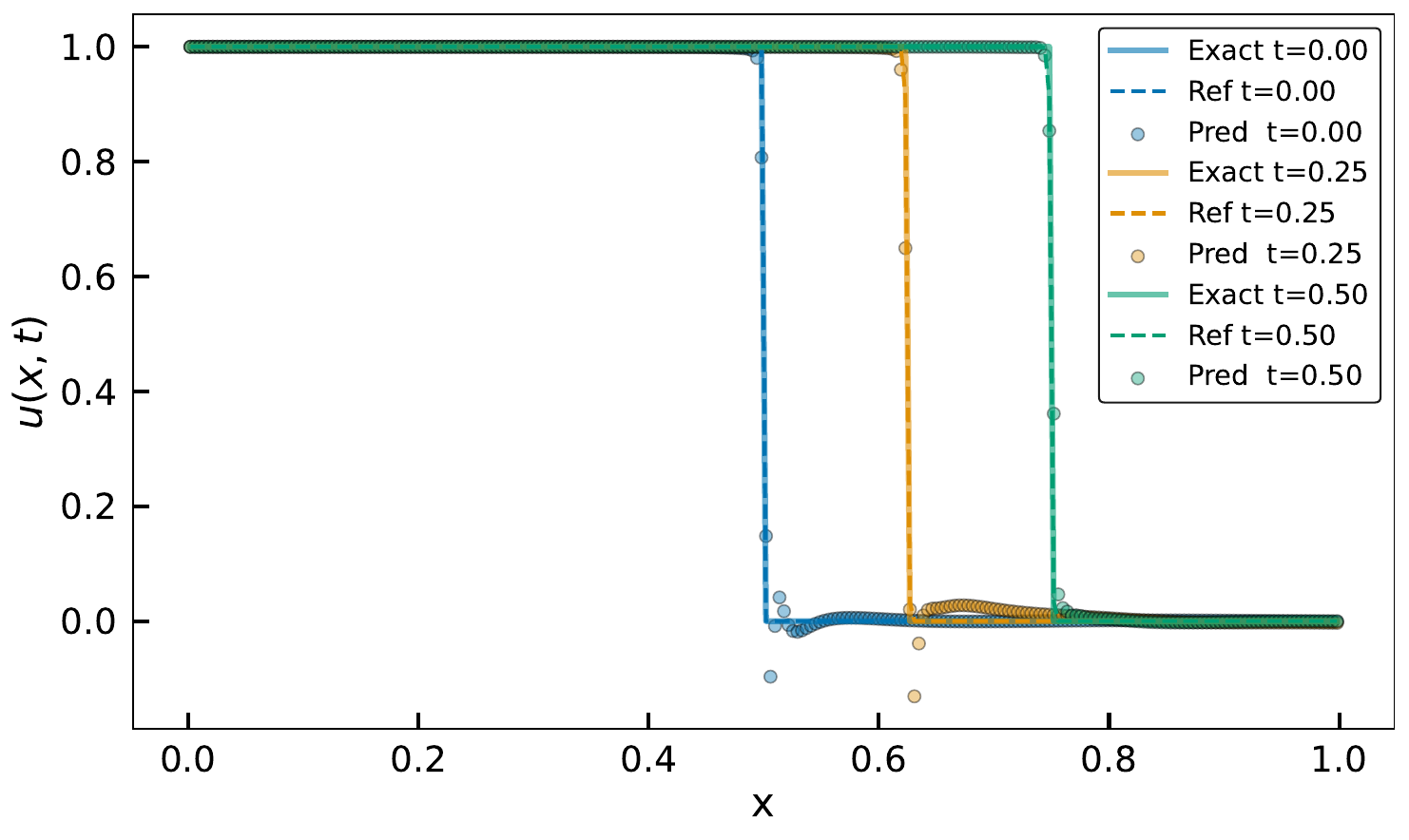}
      \subcaption{Riemann: 3 tasks}
    \end{minipage}
    \hspace{0.02\textwidth}
    \begin{minipage}[b]{0.31\textwidth}
      \centering
      \includegraphics[width=\linewidth]{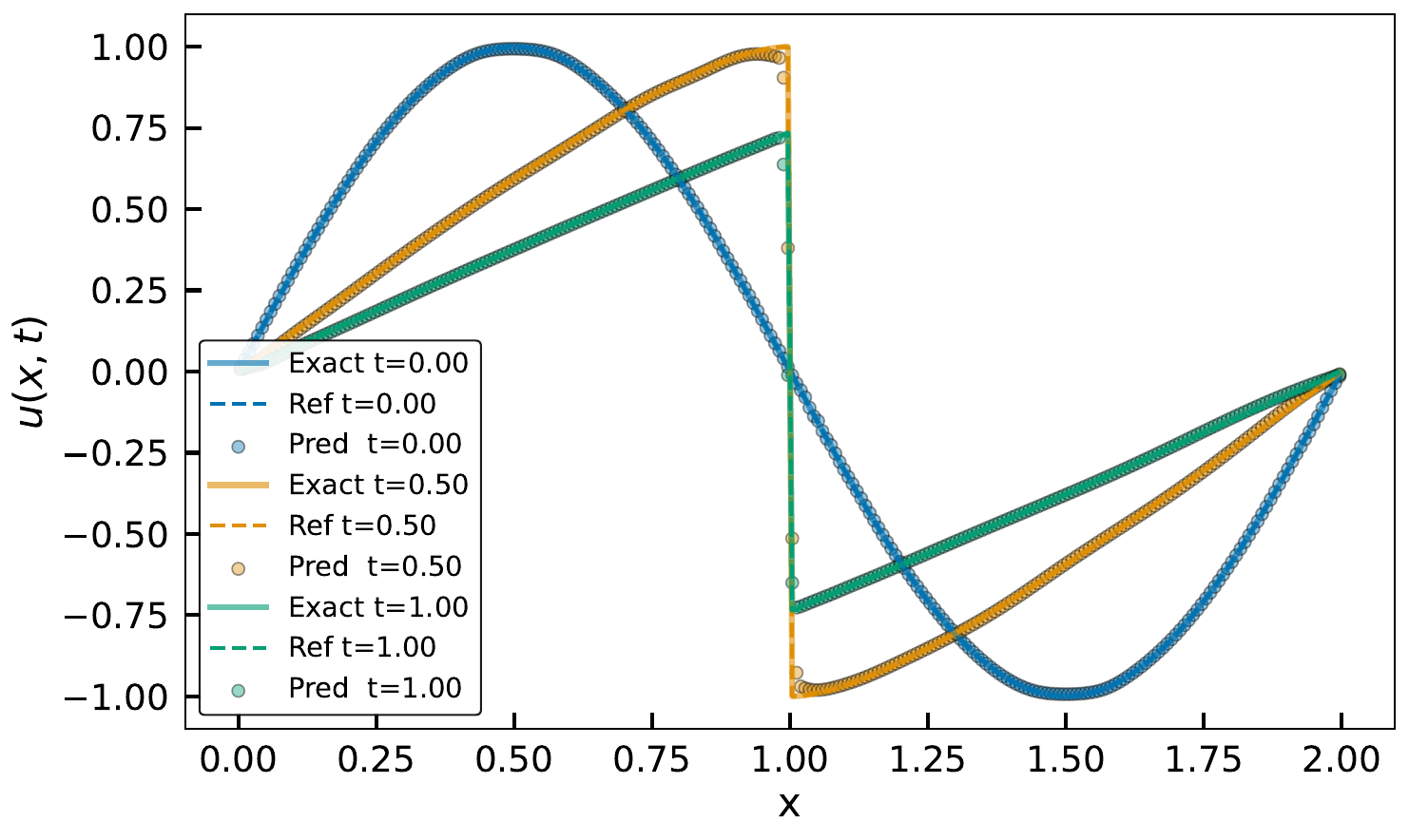}
      \subcaption{Sine: 3 tasks}
    \end{minipage}
  \end{subfigure}  

  \begin{subfigure}[b]{\textwidth}
    \centering
    \begin{minipage}[b]{0.31\textwidth}
      \centering
      \includegraphics[width=\linewidth]{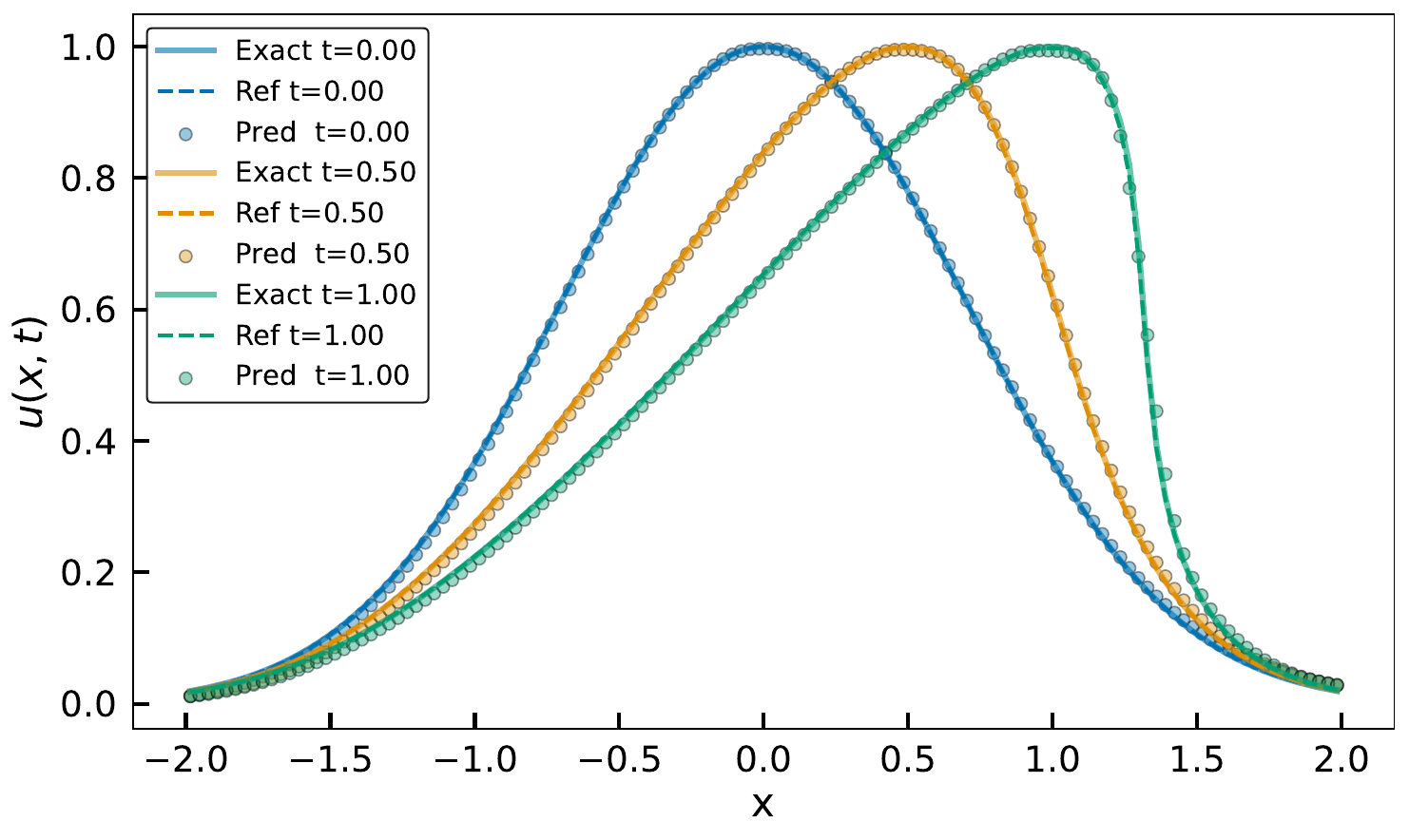}
      \subcaption{Gauss: 4 tasks}
    \end{minipage}
    \hspace{0.02\textwidth}
    \begin{minipage}[b]{0.31\textwidth}
      \centering
      \includegraphics[width=\linewidth]{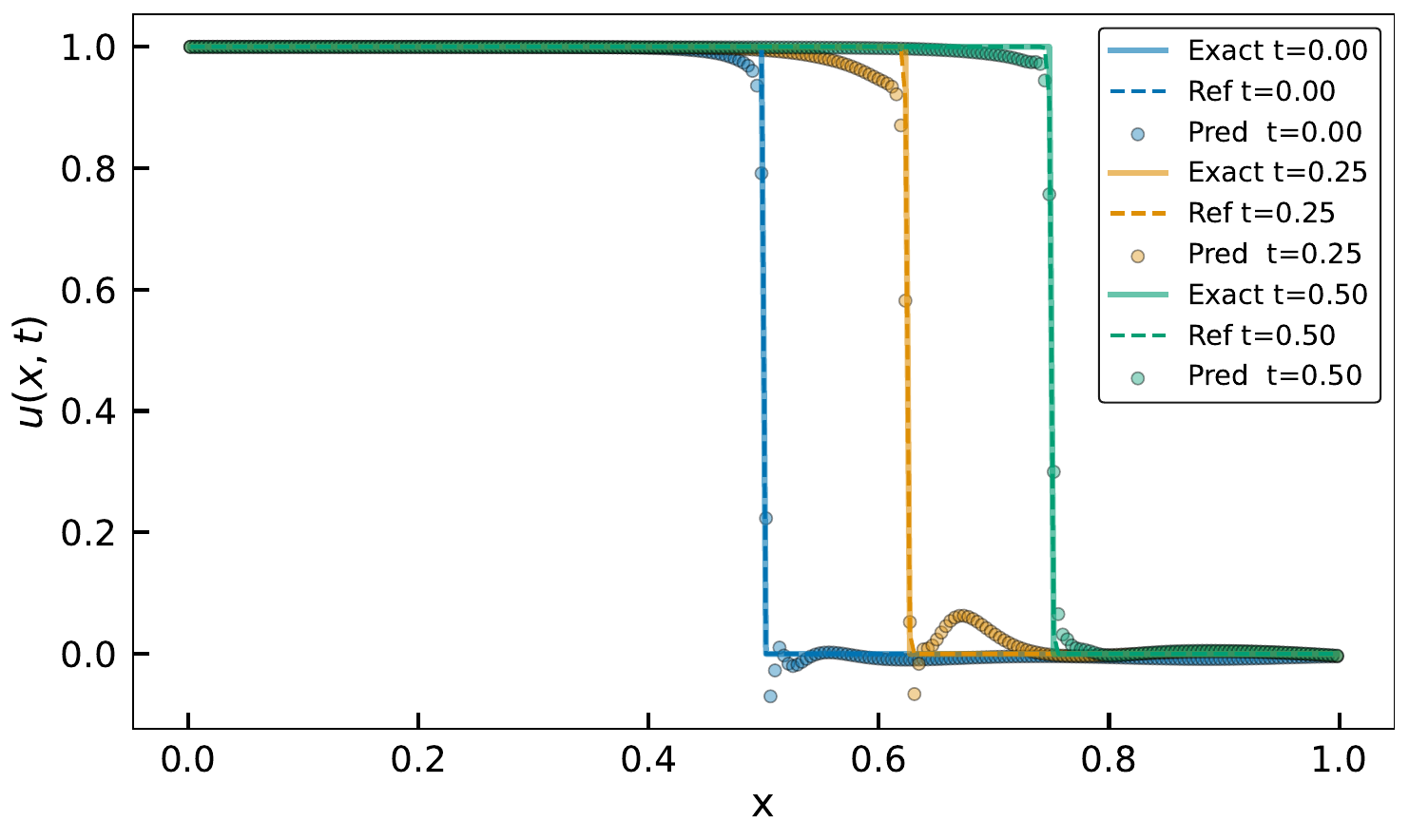}
      \subcaption{Riemann: 4 tasks}
    \end{minipage}
    \hspace{0.02\textwidth}
    \begin{minipage}[b]{0.31\textwidth}
      \centering
      \includegraphics[width=\linewidth]{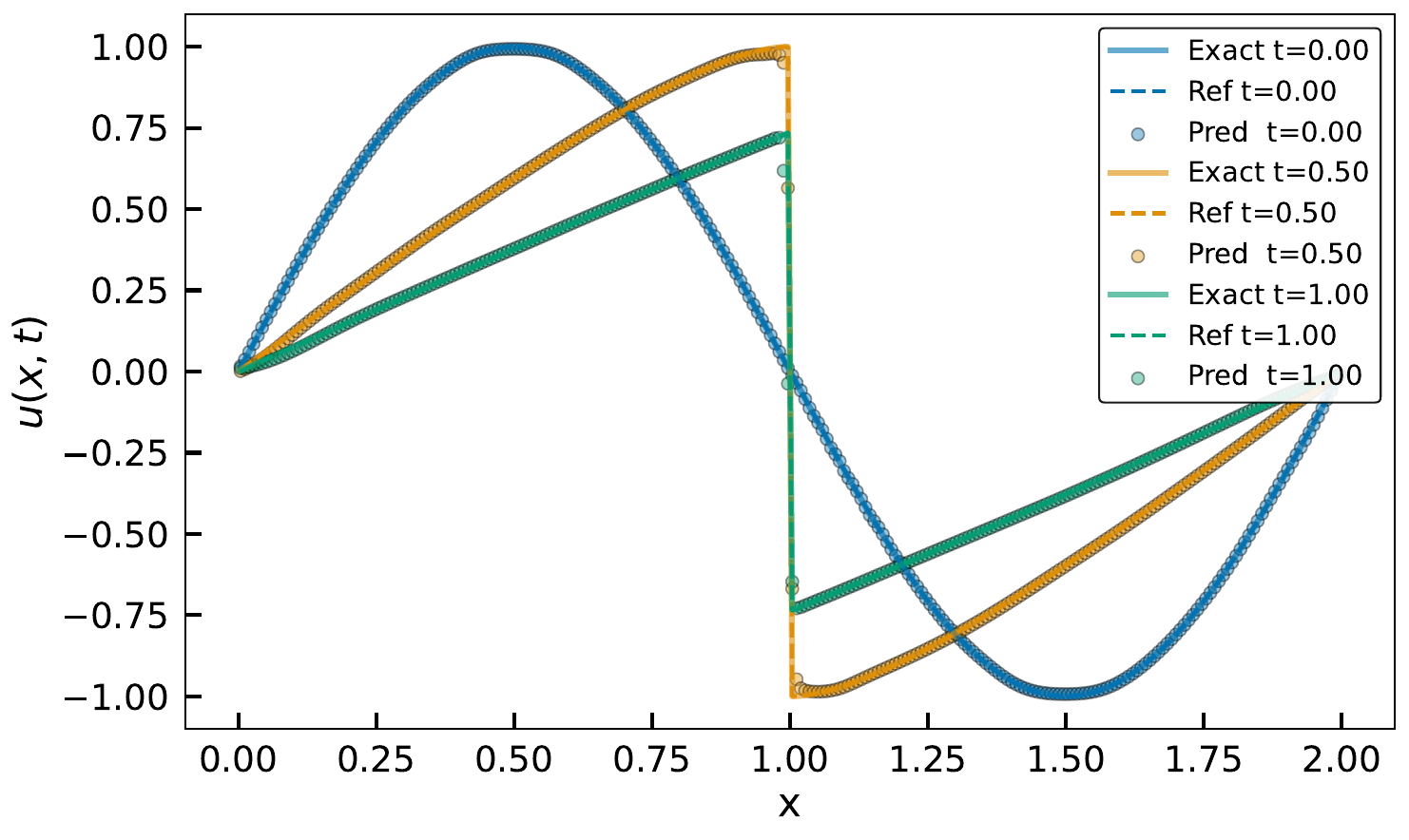}
      \subcaption{Sine: 4 tasks}
    \end{minipage}
  \end{subfigure}

  \caption{Numerical solutions of the Burgers equation under three initial conditions (columns) and four task numbers (rows).}
  \label{fig:burgers_tasks}
\end{figure}

Figure \ref{fig:burgers_tasks} provides representative predictions at different time steps under various task settings. Overall, these results demonstrate that the proposed progressive framework offers a flexible and effective approach to solving nonlinear scalar conservation laws, with competitive accuracy and clear benefits from both temporal and spatial decomposition.

\begin{figure}[!htbp]
  \centering
  
  \begin{subfigure}[b]{\textwidth}
    \centering
    \begin{minipage}[b]{0.31\textwidth}
      \centering
      \includegraphics[width=\linewidth]{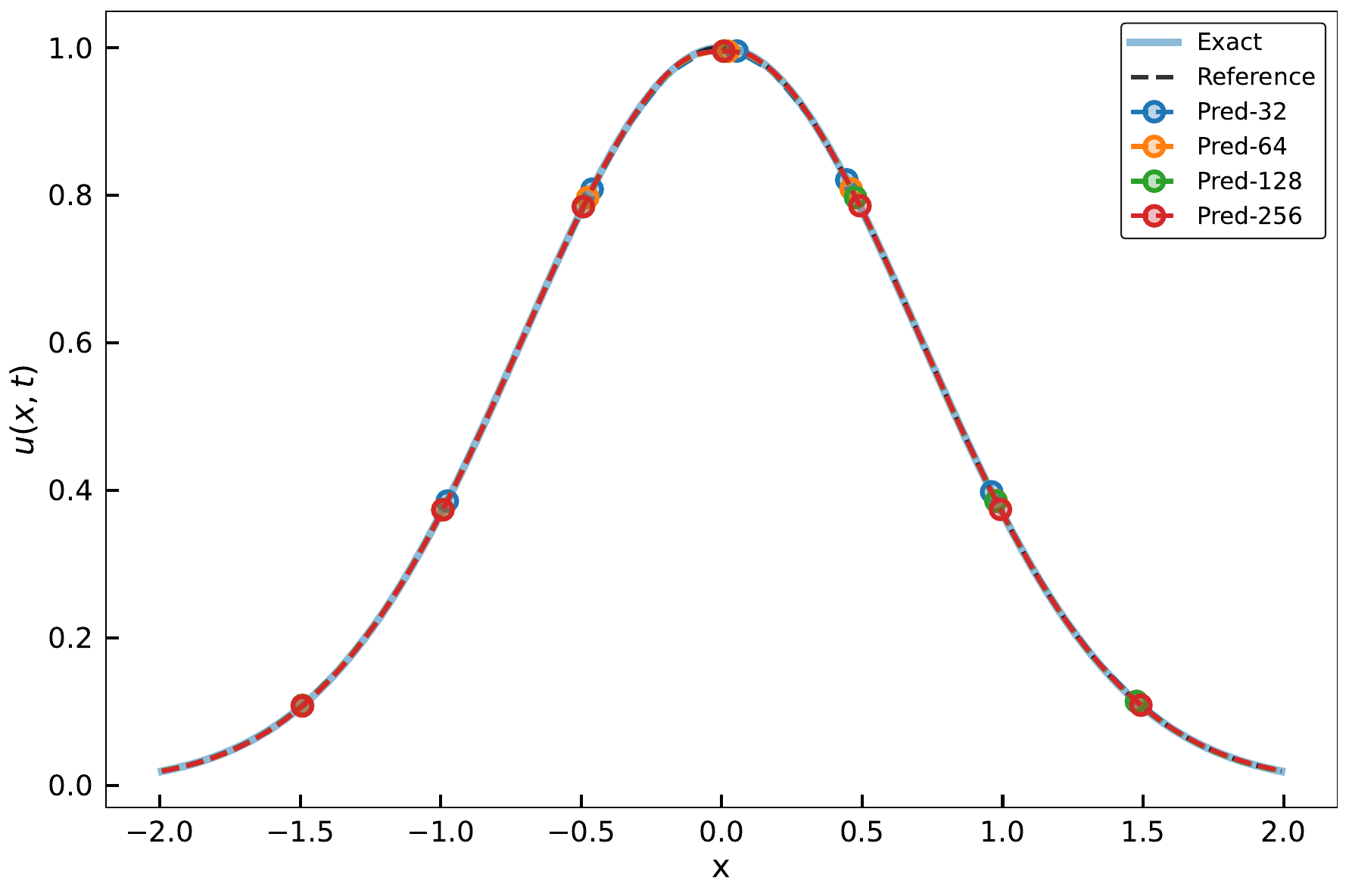}
      \caption*{Gauss case: initial condition}
    \end{minipage}
    \hspace{0.02\textwidth}
    \begin{minipage}[b]{0.31\textwidth}
      \centering
      \includegraphics[width=\linewidth]{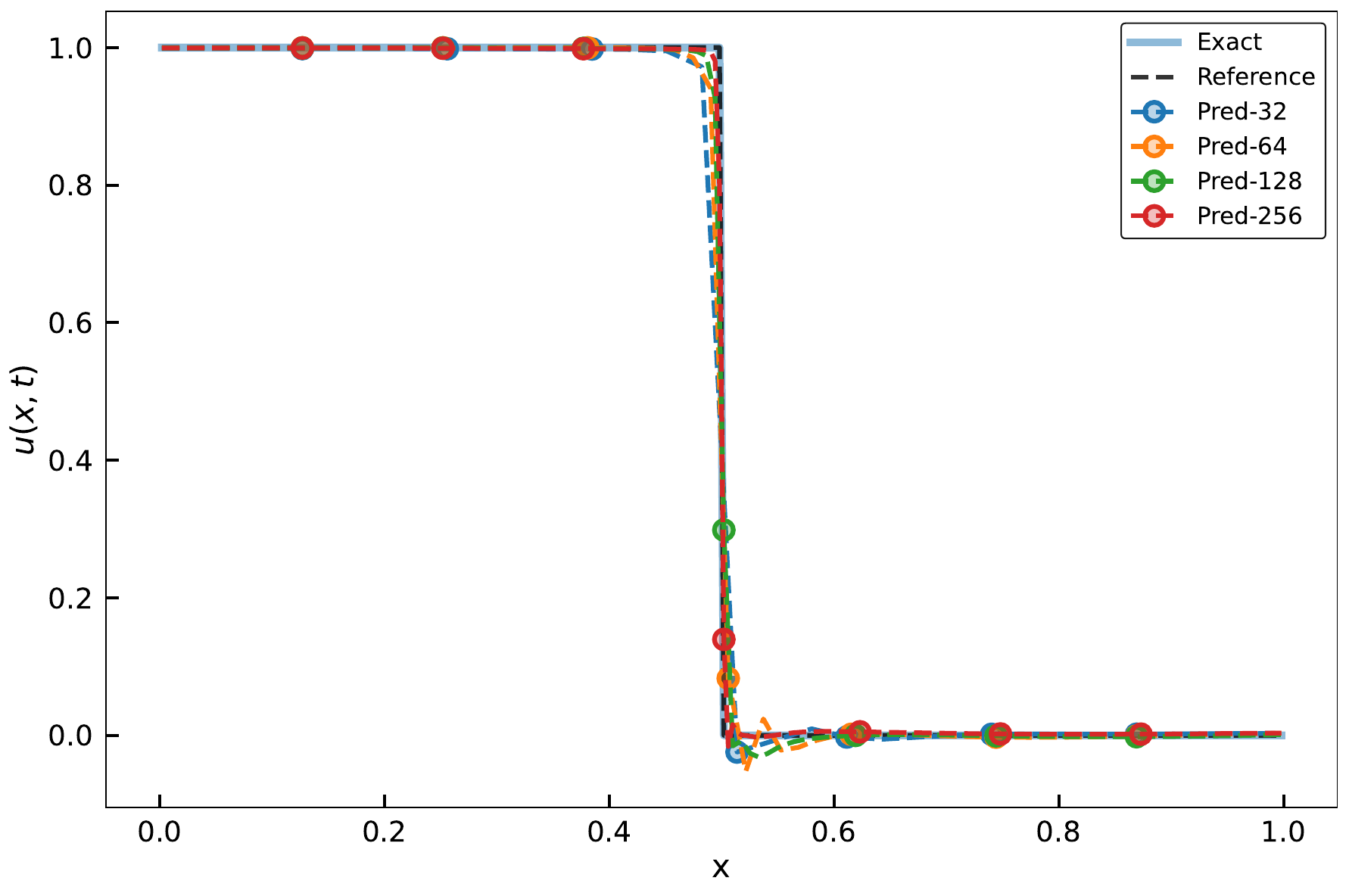}
      \caption*{Riemann case: initial condition}
    \end{minipage}
    \hspace{0.02\textwidth}
    \begin{minipage}[b]{0.31\textwidth}
      \centering
      \includegraphics[width=\linewidth]{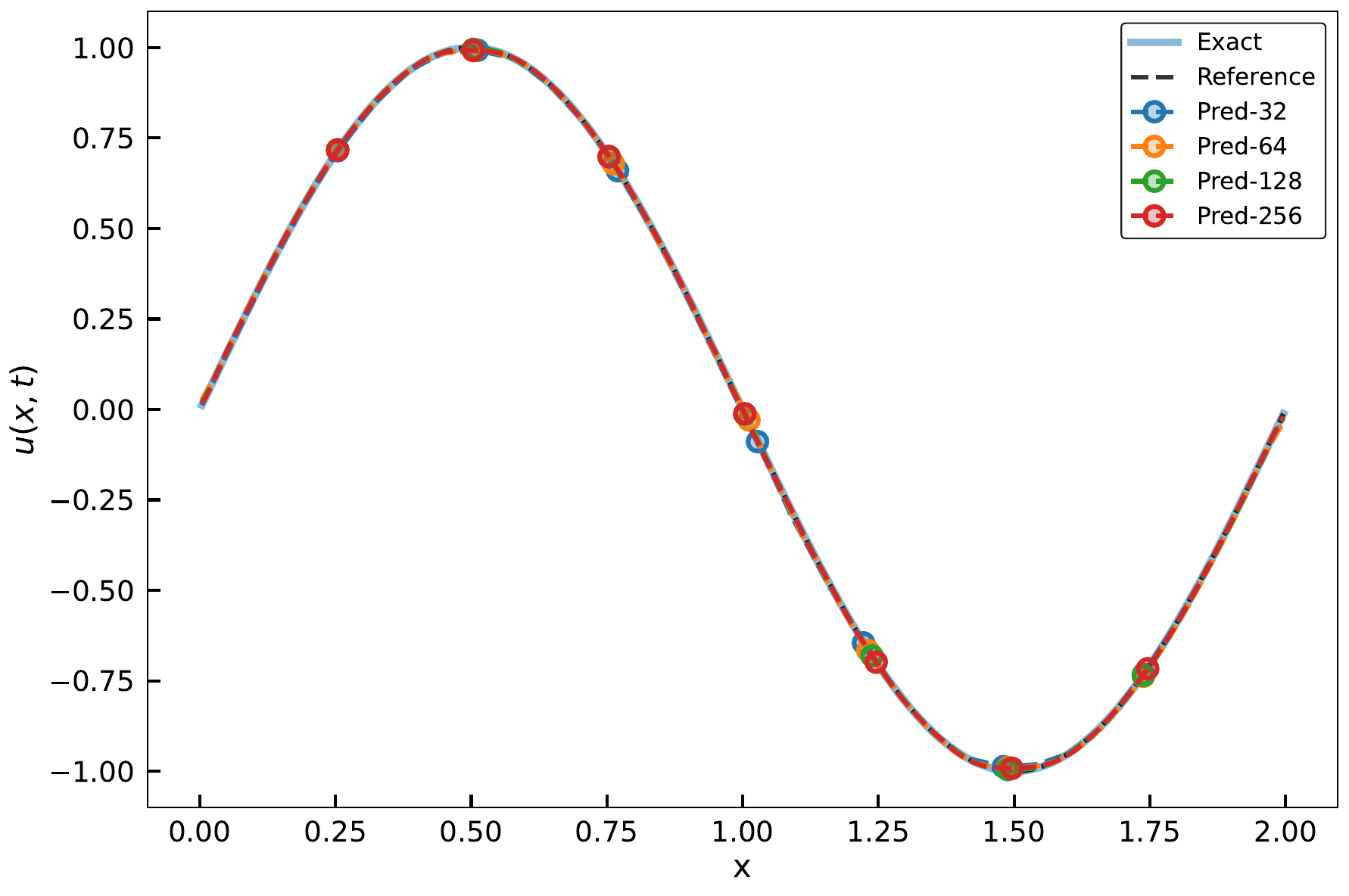}
      \caption*{Sine case: initial condition}
    \end{minipage}
    \caption{Predicted solutions at the initial time under different mesh sizes (2 tasks).}
  \end{subfigure}

  \begin{subfigure}[b]{\textwidth}
    \centering
    \begin{minipage}[b]{0.31\textwidth}
      \centering
      \includegraphics[width=\linewidth]{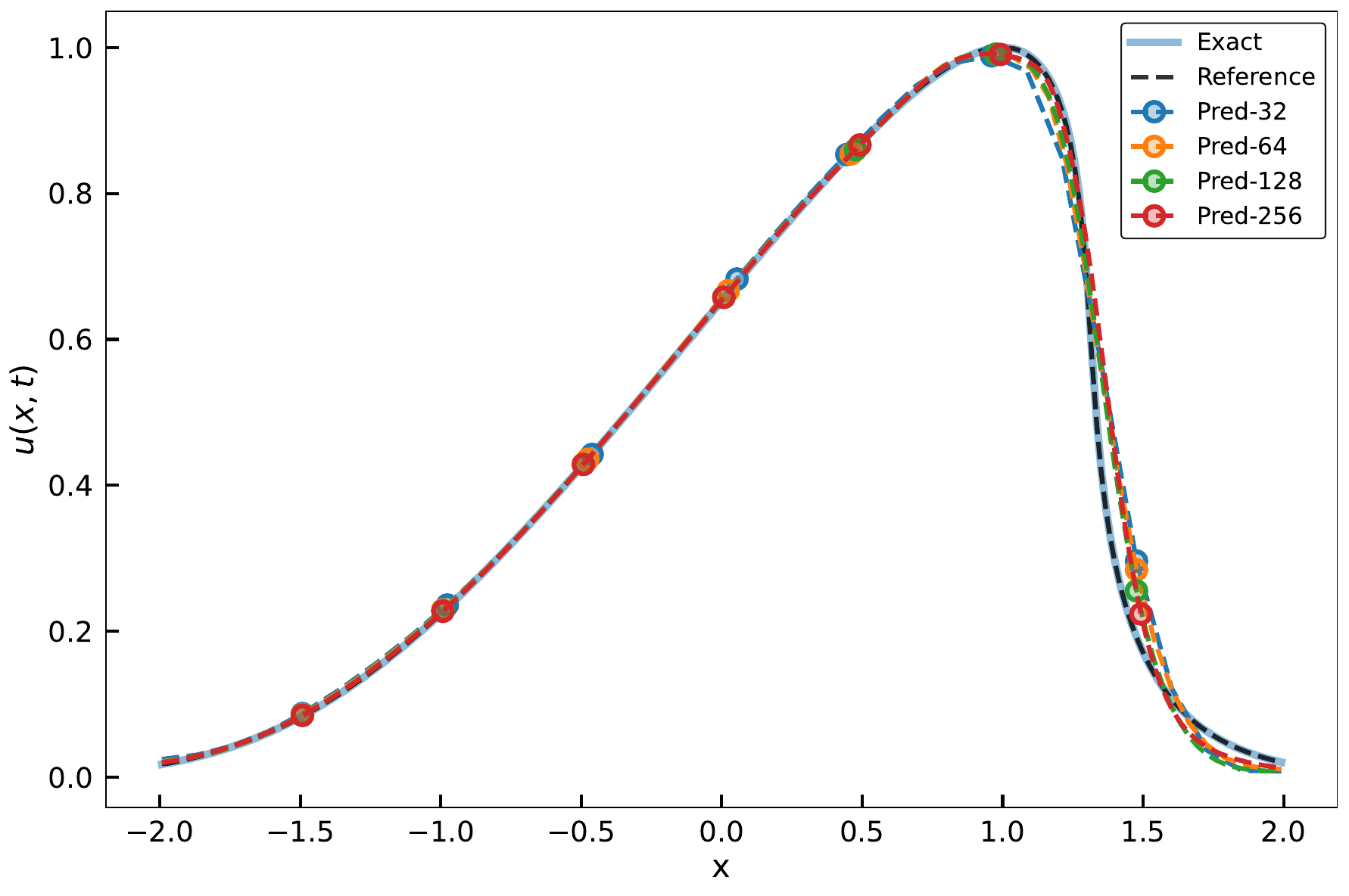}
      \caption*{Gauss case: final time}
    \end{minipage}
    \hspace{0.02\textwidth}
    \begin{minipage}[b]{0.31\textwidth}
      \centering
      \includegraphics[width=\linewidth]{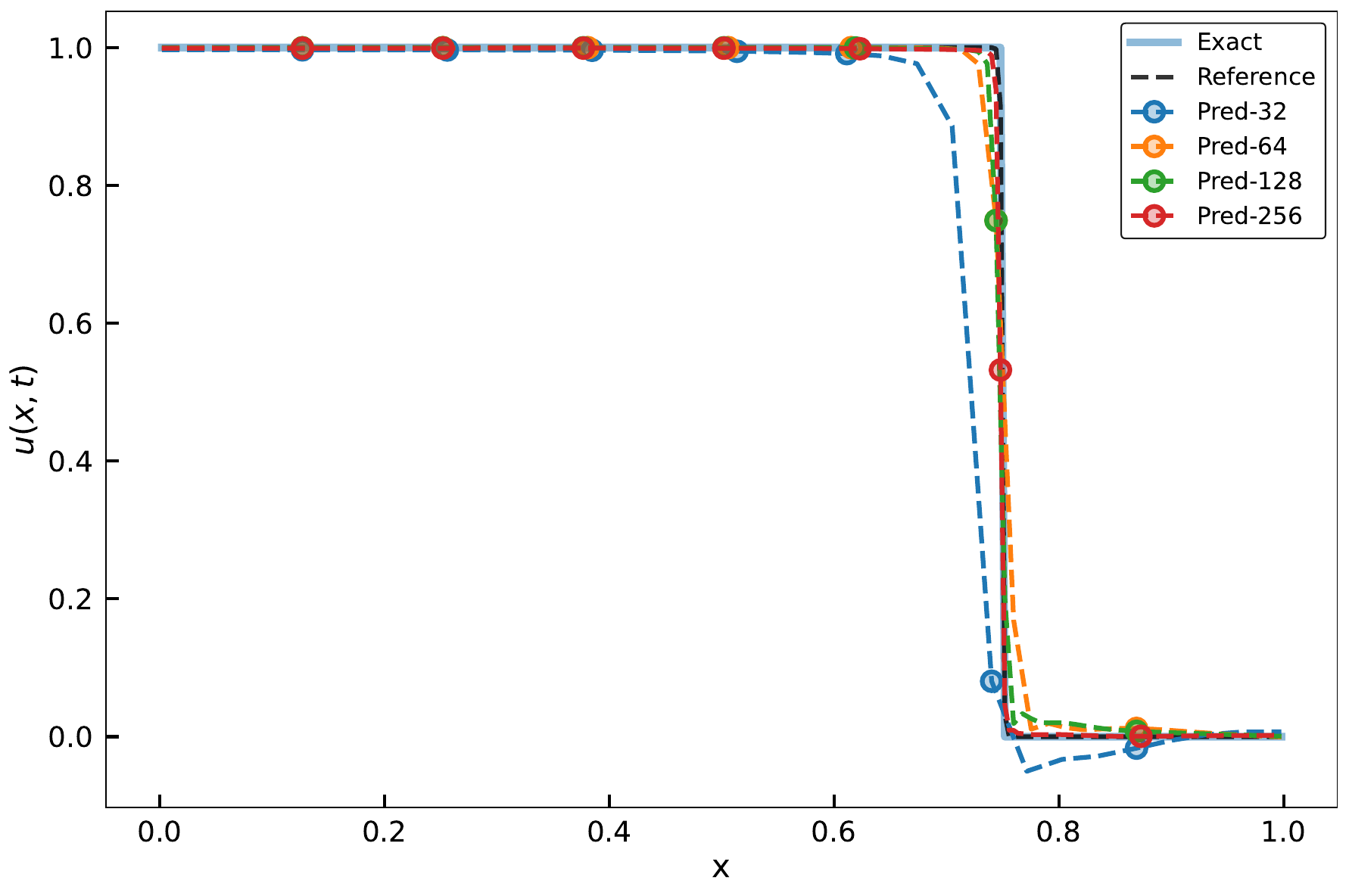}
      \caption*{Riemann case: final time}
    \end{minipage}
    \hspace{0.02\textwidth}
    \begin{minipage}[b]{0.31\textwidth}
      \centering
      \includegraphics[width=\linewidth]{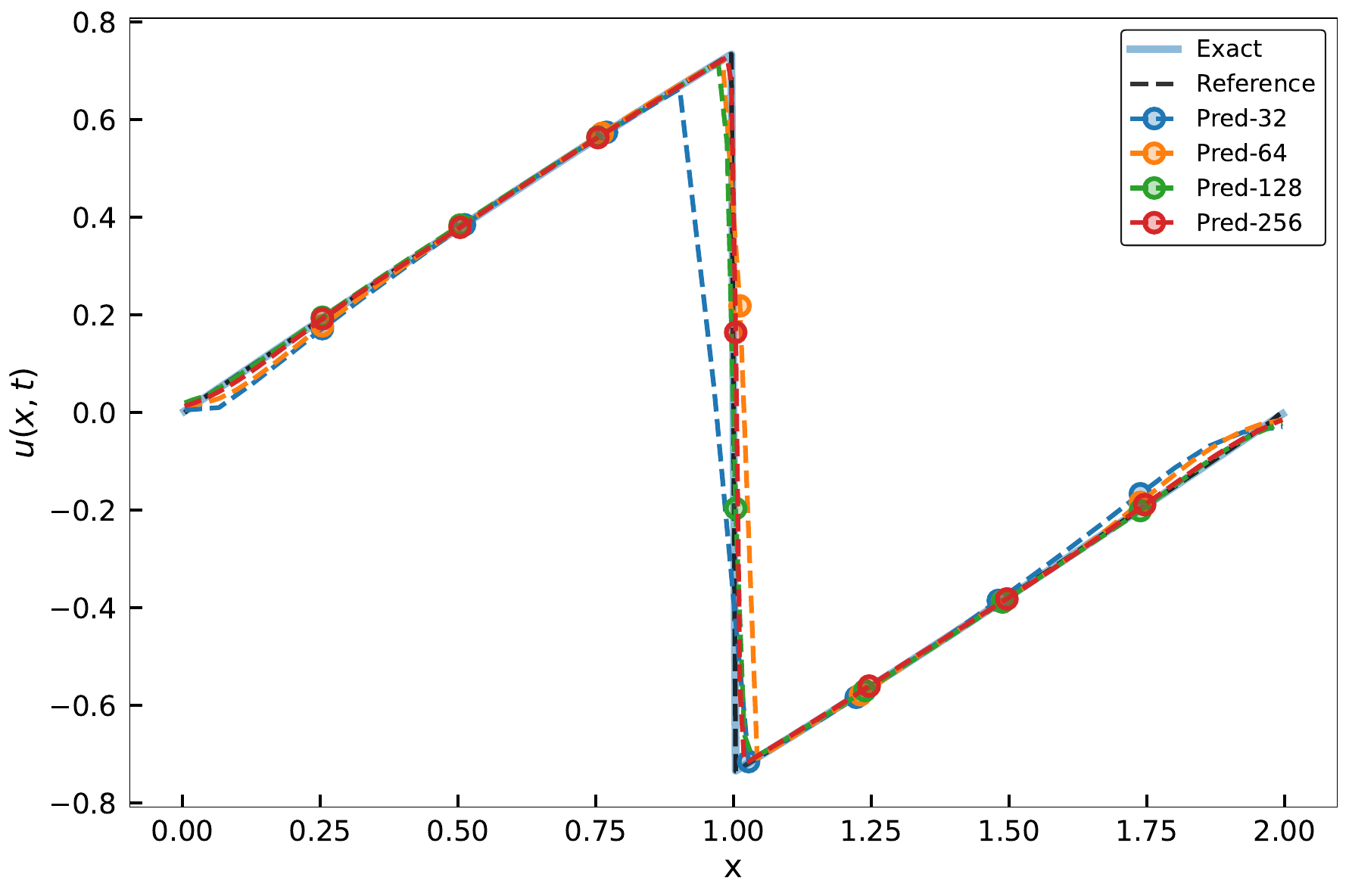}
      \caption*{Sine case: final time}
    \end{minipage}
    \caption{Predicted solutions at the final time under different mesh sizes (2 tasks).}
  \end{subfigure}

  \caption{Burgers equation predictions at the initial (top) and final (bottom) times under varying spatial resolutions for three initial conditions. All results are obtained with a two-task model.}
  \label{fig:burgers_meshes}
\end{figure}

\begin{table}[htbp]
  \centering
  \caption{{Relative $L^2$ errors of the Burgers equation for different initial conditions under varying spatial resolutions (2-task model).}}
  \label{tab:burgers-mesh}
  \sisetup{scientific-notation = true}
  \begin{tabular}{lcccc}
    \toprule
    \textbf{Initial Condition} & \textbf{32 Cells} & \textbf{64 Cells} & \textbf{128 Cells} & \textbf{256 Cells} \\
    \midrule
    Gaussian & \num{1.99e-2} & \num{1.04e-2} & \num{6.16e-3} & \textbf{\num{5.59e-3}} \\
    Riemann  & \num{1.16e-1} & \num{2.45e-2} & \textbf{\num{2.29e-2}} & \num{2.33e-2} \\
    Sine     & \num{2.51e-1} & \num{1.38e-1} & \num{7.71e-2} & \textbf{\num{3.64e-2}} \\
    \bottomrule
  \end{tabular}
\end{table}
We then fix the number of tasks to two and vary the number of spatial cells to study the effect of spatial resolution. The results are shown in Table \ref{tab:burgers-mesh} and Figure \ref{fig:burgers_meshes}. As anticipated, finer meshes lead to more accurate predictions. This trend is consistent across all initial conditions, although the improvement is more pronounced for smooth cases such as the sine wave. For discontinuous inputs like the Riemann problem, the gain from mesh refinement is more moderate, indicating the challenge of resolving shocks even with high resolution.

Building upon previous work \cite{chen2023DG}, it is observed that further mesh refinement generally leads to improved solution accuracy. This trend is also confirmed by the results presented in the figures and tables above. 

\begin{figure}[!htbp]
  \centering
  \begin{subfigure}[b]{0.31\textwidth}
    \includegraphics[width=\linewidth]{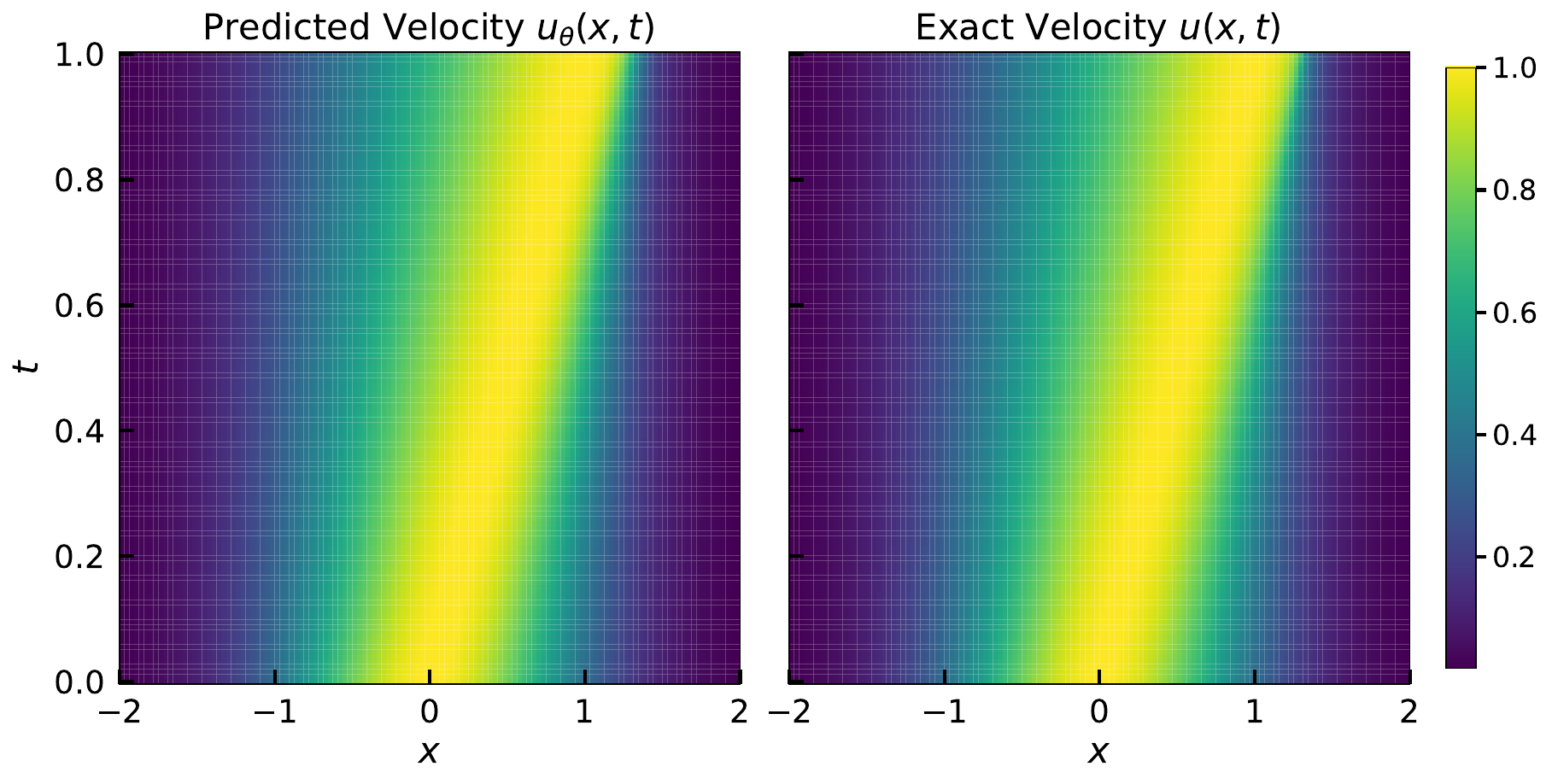}
    \caption{Gaussian Initial Condition}
  \end{subfigure}
  \hspace{0.02\textwidth}
  \begin{subfigure}[b]{0.31\textwidth}
    \includegraphics[width=\linewidth]{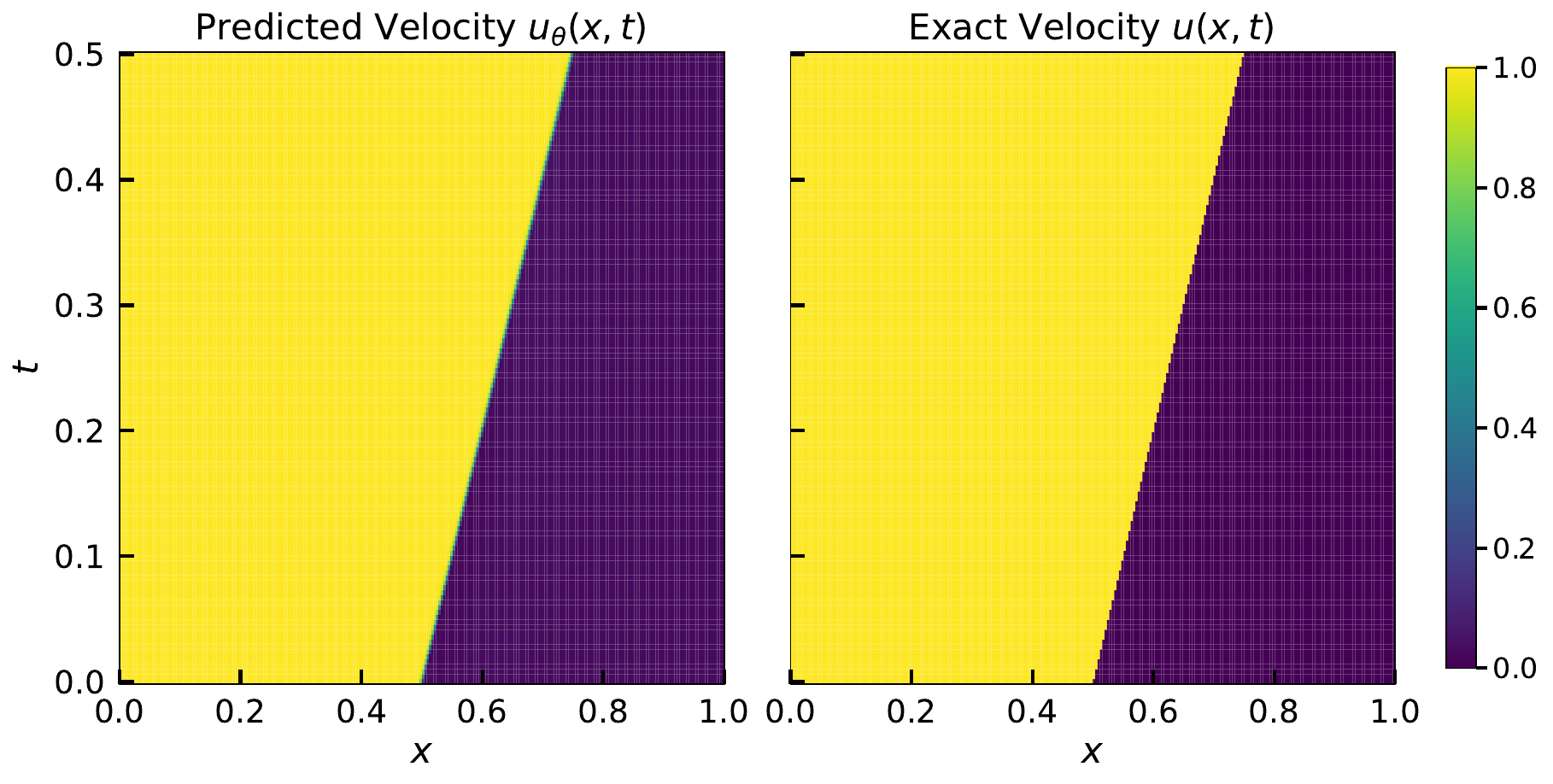}
    \caption{Riemann Initial Condition}
  \end{subfigure}
  \hspace{0.02\textwidth}
  \begin{subfigure}[b]{0.31\textwidth}
    \includegraphics[width=\linewidth]{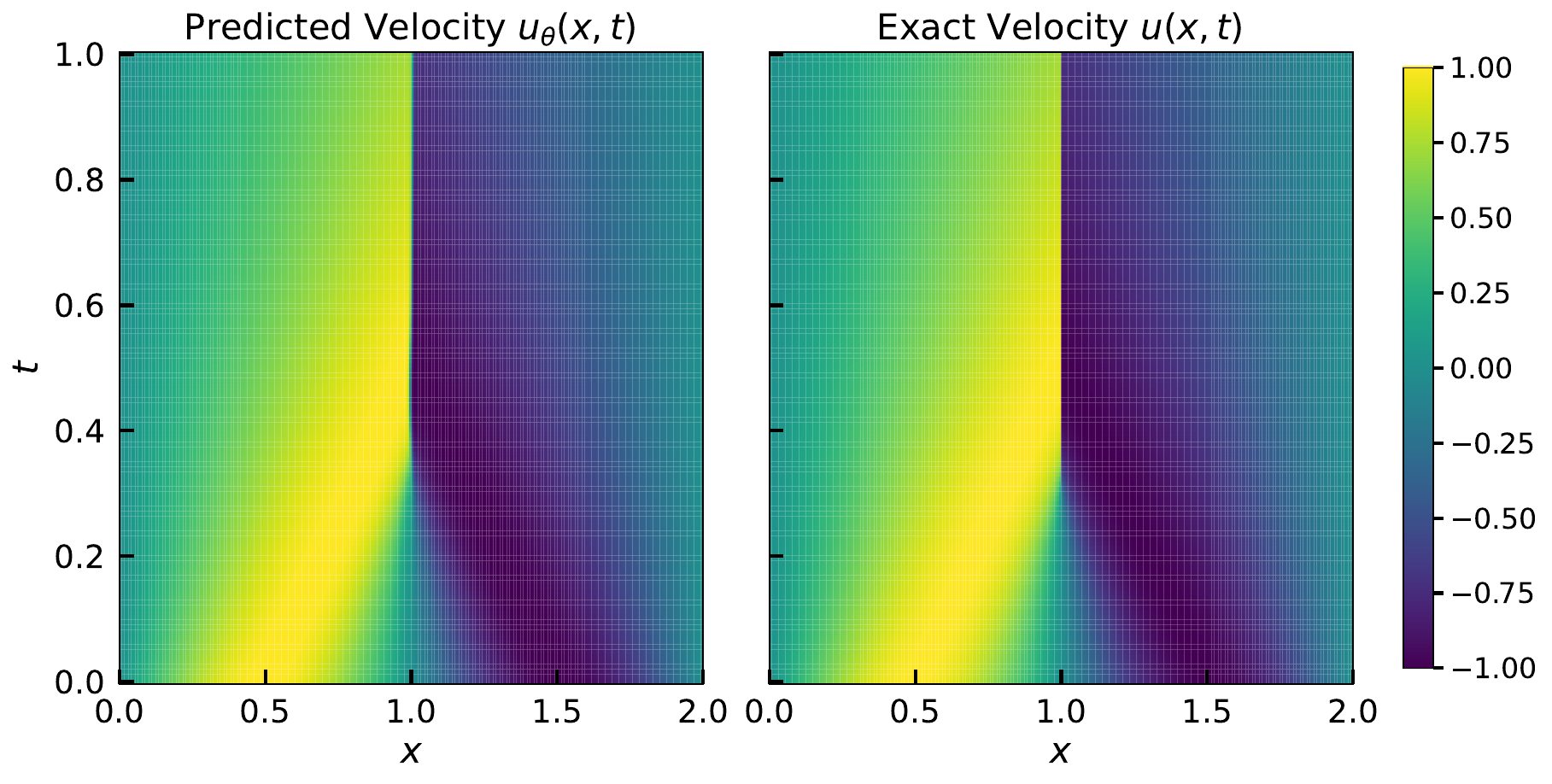}
    \caption{Sine Initial Condition}
  \end{subfigure}
  \caption{Spatio-temporal profiles of the predicted solution $u_\theta(x,t)$ for the Burgers equation under three representative initial conditions, computed with a two-task temporal decomposition.}
  \label{fig:burgers comparasion}
\end{figure}

Figure \ref{fig:burgers comparasion} depicts the spatio-temporal solution profiles computed using the proposed framework with a temporal decomposition into two tasks (task = 2). The results confirm that, even under this moderate temporal segmentation, the method reliably captures essential solution features across distinct initial conditions, achieving accurate and stable predictions throughout the simulation horizon. Similar qualitative behavior and accuracy are observed for other task decompositions, which are omitted here for brevity.

\subsection{Euler Equation}
We now turn to the more complex compressible Euler equation, using the shock tube \eqref{Sod} as an example. The model architecture and training setup is the same as the Burgers experiment, with four hidden layers per task and identical loss function weights. To better capture the richer dynamics of the Euler system, we use more collocation points internally: 16,384 points, as well as 512 boundary points, 1,024 initial points, and 4,096 pseudo-labeled points.

\begin{figure}[!htbp]
  \centering

  \begin{subfigure}[b]{\textwidth}
    \centering
    \begin{minipage}[b]{0.48\textwidth}
      \centering
      \includegraphics[width=\linewidth]{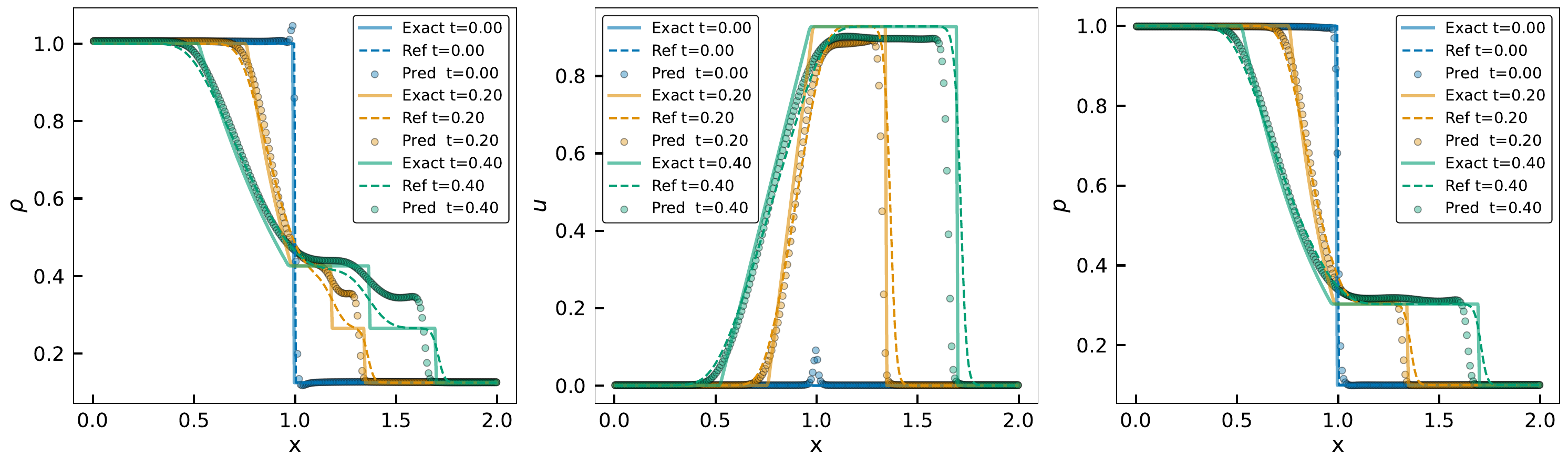}
      \caption*{Sod case: 1 task}
    \end{minipage}
    \hspace{0.02\textwidth}
    \begin{minipage}[b]{0.48\textwidth}
      \centering
      \includegraphics[width=\linewidth]{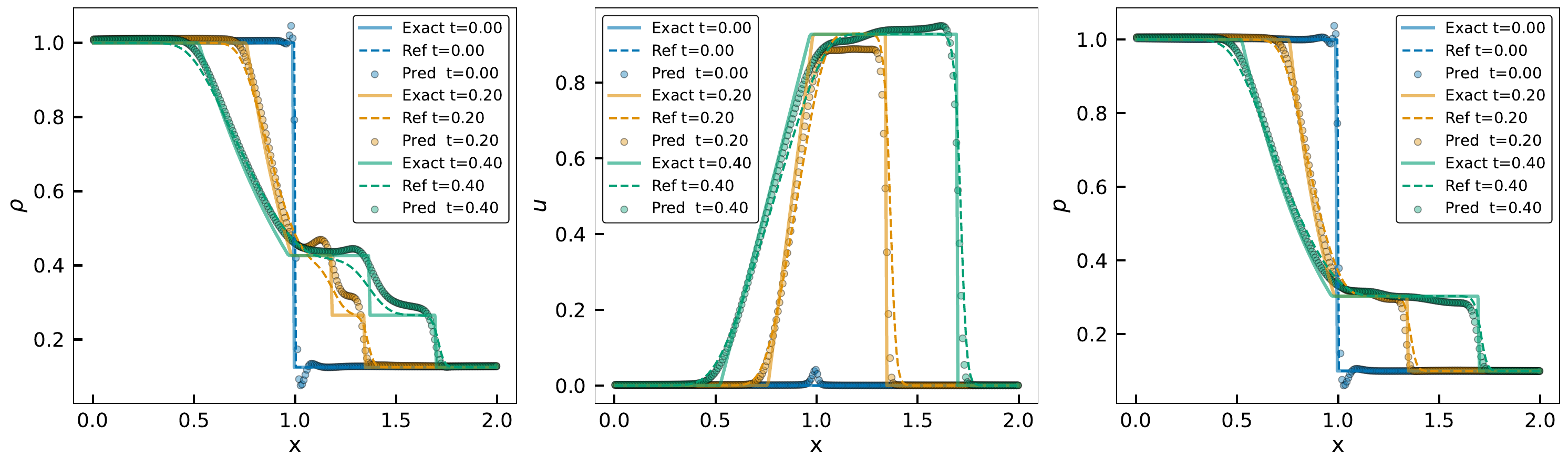}
      \caption*{Sod case: 2 tasks}
    \end{minipage}
  \end{subfigure}

  \begin{subfigure}[b]{\textwidth}
    \centering
    \begin{minipage}[b]{0.48\textwidth}
      \centering
      \includegraphics[width=\linewidth]{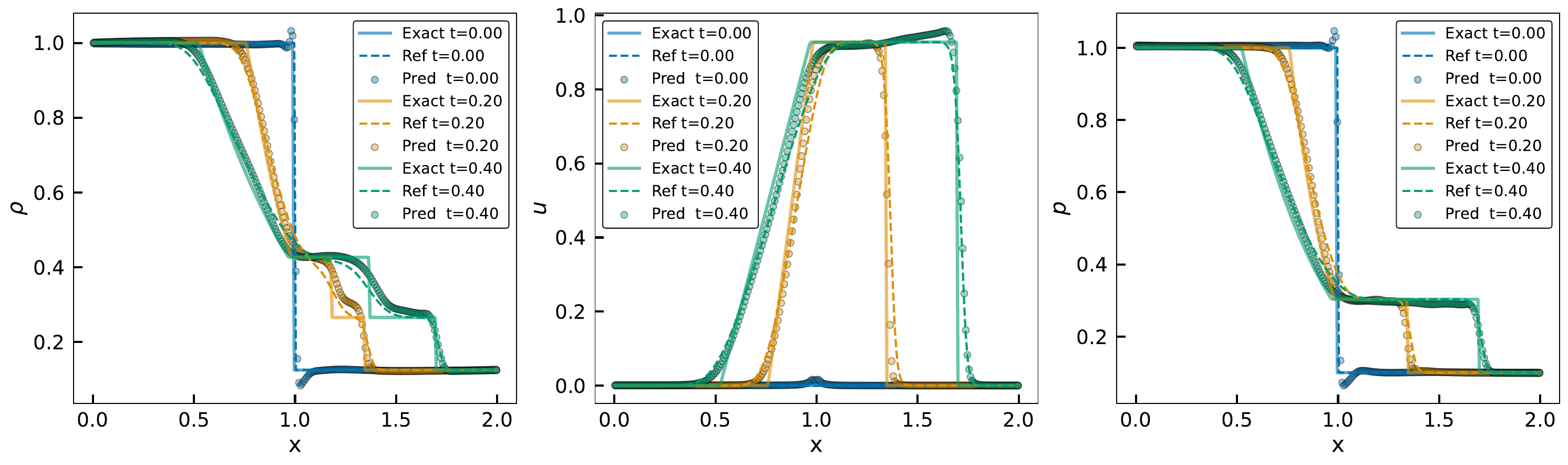}
      \caption*{Sod case: 3 tasks}
    \end{minipage}
    \hspace{0.02\textwidth}
    \begin{minipage}[b]{0.48\textwidth}
      \centering
      \includegraphics[width=\linewidth]{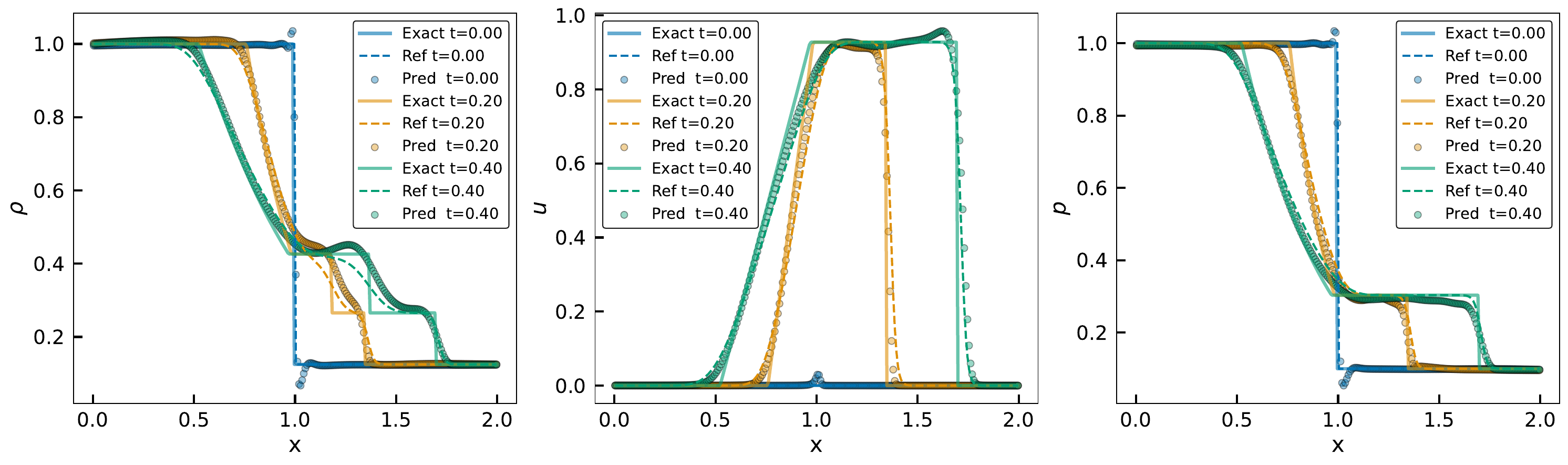}
      \caption*{Sod case: 4 tasks}
    \end{minipage}
  \end{subfigure}

  \caption{Temporal evolution of the predicted solutions for the Sod shock tube problem under different numbers of temporal tasks.}
  \label{fig:euler_tasks}
\end{figure}

\begin{figure}[!htbp]
  \centering

  \begin{subfigure}[b]{\textwidth}
    \centering
    \begin{minipage}[b]{0.48\textwidth}
      \centering
      \includegraphics[width=\linewidth]{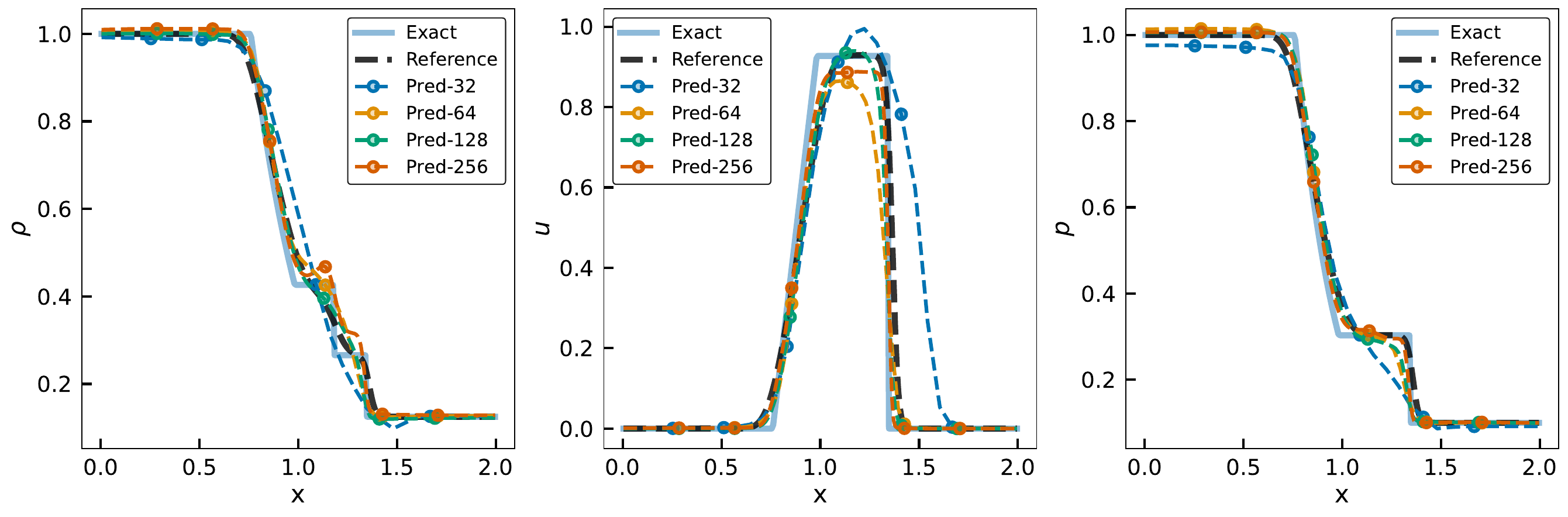}
      \caption*{Intermediate time}
    \end{minipage}
    \hspace{0.02\textwidth}
    \begin{minipage}[b]{0.48\textwidth}
      \centering
      \includegraphics[width=\linewidth]{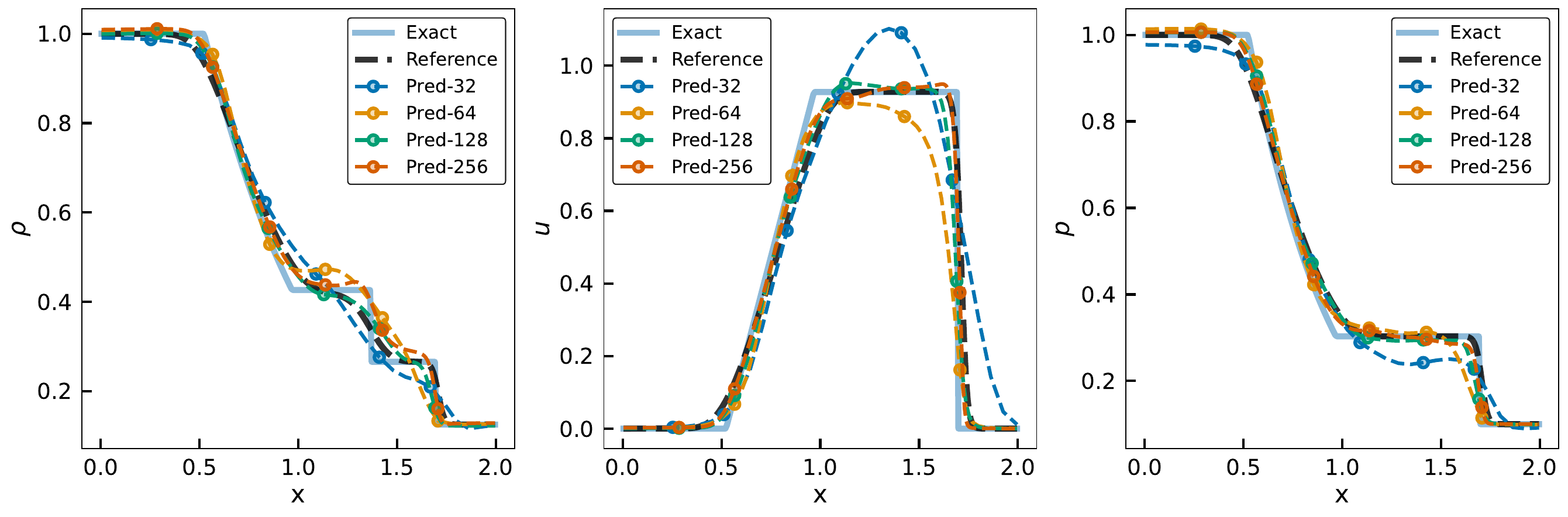}
      \caption*{Final time}
    \end{minipage}
  \end{subfigure}

  \caption{Predicted solutions of the Sod problem at intermediate and final times, computed using a two-task model with varying mesh resolutions.}
  \label{fig:euler_meshes}
\end{figure}

Table \ref{tab:euler-task} demonstrates that incorporating multiple sequential tasks notably enhances the solution accuracy for the Euler equations with the Sod initial condition. 
The progressive learning framework effectively improves the approximation of key physical variables compared to single-task training. 
However, the results also indicate that increasing the number of tasks does not always lead to further gains, suggesting an optimal range for task division that balances learning capacity and model generalization.

\begin{table}
  \centering
  \caption{{Relative $L^2$ errors of the Euler equations with Sod initial condition under varying numbers of temporal tasks (256 cells).}}
  \label{tab:euler-task}
  \sisetup{scientific-notation = true}
  \begin{tabular}{llcccc}
    \toprule
    \textbf{Initial Condition} & \textbf{Variable} & \textbf{1 Task} & \textbf{2 Tasks} & \textbf{3 Tasks} & \textbf{4 Tasks} \\
    \midrule
    \multirow{3}{*}{Sod}
      & $\rho$ & \num{3.73e-2} & \num{3.25e-2} & \textbf{\num{2.35e-2}} & \num{2.86e-2} \\
      & $u$    & \num{2.32e-1} & \num{1.12e-1} & \num{8.14e-2} & \textbf{\num{6.66e-2}} \\
      & $p$    & \num{3.54e-2} & \textbf{\num{2.50e-2}} & \num{2.46e-2} & \num{2.61e-2} \\
    \bottomrule
  \end{tabular}
\end{table}

With the number of tasks fixed at two, we examine the impact of spatial resolution on solution accuracy. Table \ref{tab:euler-mesh} demonstrates consistent improvement across all flow variables and initial conditions with increasing mesh density. 
The model's capability to capture key flow features is further demonstrated in Figure \ref{fig:euler_meshes}, which presents predicted solutions at multiple time instances under different task configurations and mesh resolutions. These results clearly show the method's robustness in accurately representing shock waves and contact discontinuities throughout the temporal evolution.

A comprehensive benchmark comparison against three baseline approaches, the standard PINN, first-order DG method, and PINNs-WE, is presented in Table \ref{tab:euler-comparison}. Under strictly controlled conditions where all neural network models employ identical architectures (four hidden layers with 64 neurons each) and matching training parameters (iterations, optimizer settings, and sampling points), our method achieves superior accuracy in velocity and pressure predictions for the Sod problem while maintaining competitive performance in density estimation. These results collectively demonstrate that the proposed framework successfully combines the strengths of DG discretizations and progressive neural networks, yielding robust and accurate solutions for challenging hyperbolic conservation problems.

\begin{table}[htbp]
  \centering
  \caption{{Relative $L^2$ errors of the Euler equations with Sod initial condition under varying spatial resolutions (2-task model).}}
  \label{tab:euler-mesh}
  \sisetup{scientific-notation = true}
  \begin{tabular}{llcccc}
    \toprule
    \textbf{Initial Condition} & \textbf{Variable} & \textbf{32 Cells} & \textbf{64 Cells} & \textbf{128 Cells} & \textbf{256 Cells} \\
    \midrule
    \multirow{3}{*}{Sod}
      & $\rho$ & \num{5.98e-2} & \num{4.60e-2} & \num{2.85e-2} & \textbf{\num{2.27e-2}} \\
      & $u$    & \num{2.82e-1} & \num{2.75e-1} & \num{1.28e-1} & \textbf{\num{5.84e-2}} \\
      & $p$    & \num{8.02e-2} & \num{6.22e-2} & \num{3.70e-2} & \textbf{\num{3.23e-2}} \\
    \bottomrule
  \end{tabular}
\end{table}

\begin{table}[htbp]
  \centering
  \caption{{Relative $L^2$ errors for the Euler equations under the Sod initial condition, using our method, vanilla PINN, PINNs-WE, and first-order DG (256 cells, 2-task setting).}}
  \label{tab:euler-comparison}
  \sisetup{scientific-notation = true}
  \begin{tabular}{llcccc}
    \toprule
    \textbf{Initial Condition} & \textbf{Variable} & \textbf{Ours} & \textbf{PINN} & \textbf{DG (k=1)} & \textbf{PINNs-WE} \\
    \midrule
    \multirow{3}{*}{Sod} 
      & $\rho$ & \num{4.05e-2} & \num{9.73e-2} & \num{3.94e-2} & {\num{3.60e-2}} \\
      & $u$    & {\num{1.25e-1}} & \num{4.76e-1} & \num{1.81e-1} & \num{1.52e-1} \\
      & $p$    & {\num{3.21e-2}} & \num{6.21e-2} & \num{4.48e-2} & \num{3.35e-2} \\
    \bottomrule
  \end{tabular}
\end{table}

\begin{figure}[!htbp]
  \centering
      \includegraphics[width=0.8\linewidth]{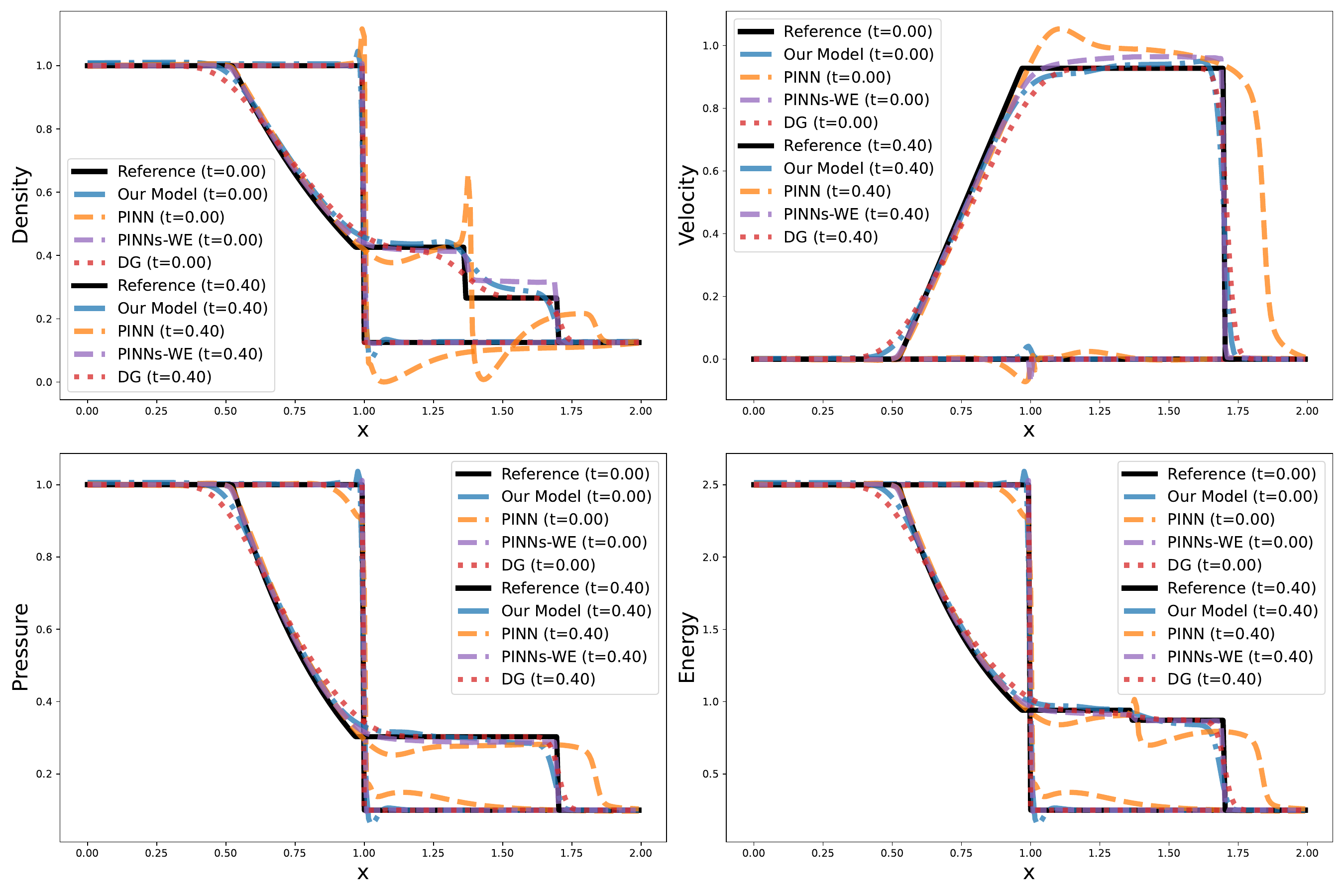}
  \caption{Comparison of Euler equation predictions at multiple time instances for density, velocity, pressure, and energy, obtained using our method, PINN, PINNs-WE, and the first-order DG method. Reference solutions are shown in black.}
  \label{fig:comparison}
\end{figure}

\section{Conclusion}
\label{sec:conclusion}

We present a hybrid framework that combines discontinuous Galerkin (DG) discretizations with temporally progressive neural networks for solving one-dimensional hyperbolic conservation laws. The method preserves the conservative structure of DG schemes while overcoming key limitations of standard physics-informed neural networks (PINNs), particularly in addressing temporal inconsistency and shock resolution challenges.

By employing progressive training across sequential time intervals augmented with physics-guided loss terms, including DG residuals and Rankine-Hugoniot conditions, our framework achieves enhanced stability and accuracy, especially for solutions containing steep gradients and discontinuities. Comprehensive numerical experiments on both Burgers and Euler equations confirm the method's superior performance compared to first-order DG solvers and baseline PINN approaches.

Future research directions include:
\begin{itemize}
    \item \textbf{High-dimensional extension:} Generalization to 2D/3D problems, including unstructured grid applications.
    \item \textbf{Multi-physics systems:} Adaptation to coupled systems like magnetohydrodynamics or reactive flows.
    \item \textbf{High-order accuracy:} Investigation of higher-degree polynomial bases for improved resolution in smooth regions.
    \item \textbf{Adaptive decomposition:} Development of error-based temporal partitioning strategies.
    \item \textbf{Parallel computing:} Leveraging the method's inherent modularity for large-scale parallel simulations.
\end{itemize}

This work demonstrates the significant potential of merging structure-preserving numerical methods with learning-based approaches for hyperbolic conservation laws, establishing a foundation for future developments in scientific machine learning.

\newpage
\appendix
\section{Details about the proofs}
\label{app:proof}
This appendix presents detailed proofs of Theorem \ref{thm:nn_dg} and supporting lemmas concerning the approximation of discontinuous Galerkin (DG) solutions by $\tanh$-activated neural networks. Key practical considerations are noted:

\begin{itemize}
\item In the progressive solver implementation, derivatives are evaluated via DG discretization. Training exclusively minimizes the scalar loss functional $\mathcal{L}^{(k)}$ over sampled points. Crucially, while the proofs require smoothness and boundedness of network maps and their derivatives, they do \emph{not} assume automatic differentiation for PDE terms during training.

\item All constants (e.g., $C_{\mathrm{coer}}^{(k)}$, $\alpha_k$, $L_k$, $C_{\mathrm{grow}}$) are explicitly defined in our assumptions. These depend on: (i) the PDE flux properties, (ii) DG discretization parameters, and (iii) network parameter norm bounds where specified.
\end{itemize}

\subsection{Preliminary lemma: boundedness of tanh networks and their derivatives}

\begin{lemma}[Boundedness of $\tanh$ networks and derivatives]
\label{lem:tanh_bounds}
Let $u_\theta:\Omega\to\mathbb{R}^m$ be a feedforward network with $\tanh$ activations in all hidden layers and parameter vector $\theta$ (weights and biases). Assume the parameter vector is constrained by $\|\theta\|_{\infty}\le B$ for some $B>0$. Then for every multi-index $\beta$ with $|\beta|\le 2$ (up to second derivatives) there exists a constant $C_{\beta}(B,{architecture})$ depending only on $B$ and the network architecture (depth, widths) such that
\begin{equation}
    \sup_{x\in\Omega} \left| \partial^\beta u_\theta(x) \right| \le C_{\beta}(B,{architecture}).
\end{equation}
Moreover the constants are finite because $\tanh$ and all its derivatives are bounded on $\mathbb{R}$ and composition with bounded-affine layers preserves boundedness under the parameter constraint.
\end{lemma}

\begin{proof}
The activation $\tanh(\cdot)$ and all its derivatives are uniformly bounded on $\mathbb R$. A feedforward layer computes an affine map followed by $\tanh$, and repeated composition yields that every partial derivative up to any fixed order is a finite combination (sum and product) of bounded functions and bounded parameters. Under the uniform parameter bound $\|\theta\|_\infty\le B$ and for a fixed architecture (finite depth and widths), repeated application of the chain rule yields finite uniform bounds for derivatives up to any fixed order. Quantitative expressions for $C_\beta$ can be obtained by iteratively applying the chain rule and bounding each factor by its uniform bound; we omit the explicit (but elementary) combinatorial formula.
\end{proof}

\subsection{Lemma A (Approximation of DG functions by tanh networks)}

\begin{lemma}[Tanh network approximation of DG functions ]
\label{lem:approx_consistent}
Let $\mathcal V_h$ denote the DG finite element space of piecewise polynomials of degree $p$ on mesh size $h$. For any $v_h\in\mathcal V_h$ and any $\varepsilon>0$ there exists a $\tanh$ network $u_\theta$ (with parameters bounded as $\|\theta\|_\infty\le B$ for sufficiently large $B$ and a finite architecture) such that
\begin{equation}
    \|u_\theta - v_h\|_{L^2(\Omega)} \le \varepsilon.
\end{equation}
Consequently we define the approximation error
\begin{equation}
    \delta_{\mathrm{approx}}^{(k)} := \inf_{\theta\in\Theta_{\mathrm{NN}}} \|u_\theta - u_{\mathrm{DG}}^{(k)}\|_{L^2(\mathcal T^{(k)}\times\Omega)},
\end{equation}
which can be made arbitrarily small by increasing network capacity.
\end{lemma}

\begin{proof}
The proof follows the two-step strategy of mollification and uniform approximation, adapted to the notation and assumptions of the main text.

\textbf{(i) Local smoothing.} Let ${K_i}$ denote the mesh elements, and let $P_i$ be the polynomial on $K_i$ such that $v_h|_{K_i} = P_i$. For $\delta \in (0, h/4)$, construct a smooth partition of unity ${\varphi_i} \subset C_c^\infty(K_i)$ satisfying $\sum_i \varphi_i \equiv 1$ on $\Omega^\delta := \bigcup_i K_i^\delta$, where $K_i^\delta$ is the $\delta$-interior of $K_i$. Define the smoothed function
\begin{equation}
    v_h^\delta(x) := \sum_i \varphi_i(x) P_i(x).
\end{equation}
As in the main text, the $L^2$ difference $\|v_h-v_h^\delta\|_{L^2(\Omega)}$ is supported only on the boundary layers of total measure $\mathcal O(\delta)$, hence for sufficiently small $\delta$ we obtain
\begin{equation}
    \|v_h-v_h^\delta\|_{L^2(\Omega)} \le \varepsilon/2.
\end{equation}

\textbf{(ii) Uniform approximation of smooth functions by tanh networks.} The smooth function $v_h^\delta$ is uniformly continuous on $\overline\Omega$. By classical universal approximation theorems for sigmoidal activations (see \cite{Hornik1991}) and quantitative $L^\infty$-approximation results for $\tanh$ networks, there exists a network $u_\theta$ with $\tanh$ activations such that
\begin{equation}
    \|u_\theta - v_h^\delta\|_{L^\infty(\Omega)} \le \frac{\varepsilon}{2\sqrt{|\Omega|}}.
\end{equation}
Consequently,
\begin{equation}
    \|u_\theta - v_h^\delta\|_{L^2(\Omega)} \le \varepsilon/2.
\end{equation}

Combining the two steps via the triangle inequality yields the desired bound:
\begin{equation}
    \|u_\theta - v_h\|_{L^2(\Omega)} \le \varepsilon.
\end{equation}

Finally, note that Lemma \ref{lem:tanh_bounds} ensures the approximating network can be constructed with bounded parameters and derivatives, a property crucial for subsequent coercivity and stability analyses.
\end{proof}

\subsection[Lemma B (Residual coercivity: loss controls L2 error)]{Lemma B (Residual coercivity: loss controls $L^2$ error)}

\begin{lemma}[Residual-to-error coercivity (local version)]
\label{lem:coercivity_consistent}
Let $u_{\mathrm{DG}}^{(k)}$ be the DG solution on $\mathcal T^{(k)}$. Define the total loss
\begin{equation}
    \mathcal L^{(k)}(\theta) := \mathcal L_{\mathrm{DG}}^{(k)}(\theta) + \mathcal L_{\mathrm{IC}}^{(k)}(\theta)
    + \mathcal L_{\mathrm{bdy}}^{(k)}(\theta) + \mathcal L_{\mathrm{RH}}^{(k)}(\theta),
\end{equation}
where each term is an $L^2$-type squared residual evaluated via the DG quadrature/weak form (as implemented in the algorithm). Assume:
\begin{enumerate}[(i)]
  \item The residual operator $F:\mathcal X\to\mathcal Y$ (discrete DG residual mapping) is Fr\'echet differentiable in a neighborhood of $u_{\mathrm{DG}}^{(k)}$.
  \item The linearization $F'(u_{\mathrm{DG}}^{(k)})$ satisfies the coercivity estimate
  \begin{equation}
      \langle F'(u_{\mathrm{DG}}^{(k)}) w, w \rangle_{\mathcal Y} \ge \alpha_k \|w\|_{\mathcal X}^2,
  \end{equation}
  for some $\alpha_k>0$ and all admissible $w$.
  \item The nonlinear remainder satisfies the quadratic bound
  \begin{equation}
      \|F(u_{\mathrm{DG}}^{(k)}+w) - F(u_{\mathrm{DG}}^{(k)}) - F'(u_{\mathrm{DG}}^{(k)}) w\|_{\mathcal Y} \le L_k \|w\|_{\mathcal X}^2,
  \end{equation}
  for $\|w\|_{\mathcal X}\le r_k$.
  \item The loss is equivalent to the residual norm up to a constant weight:
  \begin{equation}
      \|F(u_\theta)\|_{\mathcal Y} \le C_{\mathrm{loss}} \, \mathcal L^{(k)}(\theta)^{1/2}.
  \end{equation}
\end{enumerate}
Then there exist $\rho_k\in(0,r_k]$ and a constant $C_{\mathrm{coer}}^{(k)}>0$ (for example $C_{\mathrm{coer}}^{(k)}=\tfrac{2}{\alpha_k}C_{\mathrm{loss}}$) such that whenever
\begin{equation}
    \|u_\theta - u_{\mathrm{DG}}^{(k)}\|_{\mathcal X} \le \rho_k,
\end{equation}
we have
\begin{equation}
    \|u_\theta - u_{\mathrm{DG}}^{(k)}\|_{\mathcal X} \le C_{\mathrm{coer}}^{(k)} \, \mathcal L^{(k)}(\theta)^{1/2}.
\end{equation}
\end{lemma}

\begin{proof}
Set $e := u_\theta - u_{\mathrm{DG}}^{(k)}$. Since $F(u_{\mathrm{DG}}^{(k)})=0$, Taylor expansion yields
\begin{equation}
    F(u_\theta) = F(u_{\mathrm{DG}}^{(k)} + e) = F'(u_{\mathrm{DG}}^{(k)}) e + R(e),
\end{equation}
with the remainder $\|R(e)\|_{\mathcal Y}\le L_k\|e\|_{\mathcal X}^2$ by assumption. Take the duality pairing with $e$ and apply coercivity to obtain
\begin{equation}
    \alpha_k \|e\|_{\mathcal X}^2 \le \langle F'(u_{\mathrm{DG}}^{(k)}) e, e\rangle_{\mathcal Y}
    = \langle F(u_\theta), e\rangle_{\mathcal Y} - \langle R(e), e\rangle_{\mathcal Y}.
\end{equation}
Using Cauchy-Schwarz and the remainder bound,
\begin{equation}
    \alpha_k \|e\|_{\mathcal X}^2 \le \|F(u_\theta)\|_{\mathcal Y}\|e\|_{\mathcal X} + L_k\|e\|_{\mathcal X}^3.
\end{equation}
Divide by $\|e\|_{\mathcal X}$ (for $e\neq 0$) to get
\begin{equation}
    \alpha_k \|e\|_{\mathcal X} \le \|F(u_\theta)\|_{\mathcal Y} + L_k\|e\|_{\mathcal X}^2.
\end{equation}
Choose $\rho_k>0$ so that $L_k \rho_k \le \alpha_k/2$. If $\|e\|_{\mathcal X}\le \rho_k$ then
\begin{equation}
    \frac{\alpha_k}{2}\|e\|_{\mathcal X} \le \|F(u_\theta)\|_{\mathcal Y}.
\end{equation}
Hence
\begin{equation}
    \|e\|_{\mathcal X} \le \frac{2}{\alpha_k}\|F(u_\theta)\|_{\mathcal Y}
    \le \frac{2}{\alpha_k} C_{\mathrm{loss}} \, \mathcal L^{(k)}(\theta)^{1/2},
\end{equation}
which proves the claim with $C_{\mathrm{coer}}^{(k)} = \tfrac{2}{\alpha_k} C_{\mathrm{loss}}$.
\end{proof}

\begin{remark}[DG computation of PDE derivatives]
In practice, the residual $F(u_\theta)$ is evaluated using the DG weak form and numerical quadrature, not through symbolic differentiation of $u_\theta$. The coercivity assumption concerns the mapping $u\mapsto F(u)$ as evaluated by the DG discretization. While Lemma \ref{lem:tanh_bounds} ensures $u_\theta$ is sufficiently smooth for well-defined pointwise evaluations and quadrature-based residuals, the implementation does not require automatic differentiation of $u_\theta$. Instead, the DG operator computes residuals using its discrete flux differences and elementwise integrals.
\end{remark}

\subsection{Lemma C (Interface stability and temporal error propagation)}

\begin{lemma}[Temporal stability and error propagation]
\label{lem:stability_consistent}
Let $\mathcal E_k$ denote the DG evolution operator mapping initial data at $t_{k-1}$ to the discrete solution at $t_k$ (on interval $\mathcal T^{(k)}$). Suppose there exists $C_{\mathrm{grow}}>0$ such that for any two admissible initial data $u_0,\tilde u_0$ with $\|u_0-\tilde u_0\|_{L^2}\le \rho$ the following Lipschitz stability estimate holds:
\begin{equation}
    \|\mathcal E_k(u_0) - \mathcal E_k(\tilde u_0)\|_{L^2} \le S_k \|u_0-\tilde u_0\|_{L^2},
    \qquad S_k := e^{C_{\mathrm{grow}}\Delta t_k},
\end{equation}
with $\Delta t_k = t_k-t_{k-1}$. Assuming the stability factors $S_k$ are uniformly bounded for the class of admissible solutions, the global NN-to-exact error satisfies
\begin{equation}
    \|u_{\mathrm{NN}} - u_{\mathrm{exact}}\|_{L^2(\Omega\times[0,T])}
    \le \sum_{k=1}^M \Big(\prod_{j=k+1}^M S_j\Big)\|u_{\mathrm{NN}}^{(k)} - u_{\mathrm{DG}}^{(k)}\|_{L^2(\mathcal T^{(k)}\times\Omega)}
    + \|u_{\mathrm{DG}} - u_{\mathrm{exact}}\|_{L^2(\Omega\times[0,T])}.
\end{equation}
\end{lemma}

\begin{proof}
The proof follows from standard energy estimates. Define the error $e^{(k)}(t) := u^{(k)}_{\mathrm{DG}}(t)-\tilde{u}^{(k)}(t)$ as the difference between the DG evolution starting from accurate initial data and the evolution starting from perturbed initial data (due to neural network approximation at the time interface). Applying elementwise energy estimates and Gronwall's inequality yields, for $t\in\mathcal{T}^{(k)}$,
\begin{equation}
    \|e^{(k)}(t_k)\|_{L^2} \le e^{C_{\mathrm{grow}}\Delta t_k} \|e^{(k)}(t_{k-1})\|_{L^2}.
\end{equation}
Iterating these bounds from $k=1$ to $M$ and summing the propagated contributions of the local NN-to-DG errors on each interval produces the stated global bound. The DG-to-exact term $\|u_{\mathrm{DG}}-u_{\mathrm{exact}}\|_{L^2}$ accounts for the spatial discretization error described in Theorem \ref{thm:dg_exact}.
\end{proof}

\subsection{Proof of Theorem \ref{thm:nn_dg} (NN-to-DG error bound)}

\begin{proof}[Proof of Theorem \ref{thm:nn_dg}]
Combine Lemma \ref{lem:approx_consistent} and Lemma \ref{lem:coercivity_consistent}. Let $\theta^{(k)}$ denote the parameters obtained by training on interval $\mathcal T^{(k)}$, and define the optimization error
\begin{equation}
    \epsilon_{\mathrm{opt}}^{(k)} := \mathcal L^{(k)}(\theta^{(k)}) - \inf_{\theta}\mathcal L^{(k)}(\theta).
\end{equation}
By Lemma \ref{lem:coercivity_consistent}, whenever $\|u_{\theta^{(k)}} - u_{\mathrm{DG}}^{(k)}\|_{\mathcal X}\le \rho_k$ we have
\begin{equation}
    \|u_{\theta^{(k)}} - u_{\mathrm{DG}}^{(k)}\|_{\mathcal X}
    \le C_{\mathrm{coer}}^{(k)} \, \mathcal L^{(k)}(\theta^{(k)})^{1/2}.
\end{equation}
Write $\mathcal L^{(k)}(\theta^{(k)}) = \epsilon_{\mathrm{opt}}^{(k)} + \inf_\theta \mathcal L^{(k)}(\theta)$. The infimum is attained (or approximated) by any parameter realizing the approximation of $u_{\mathrm{DG}}^{(k)}$ up to $\delta_{\mathrm{approx}}^{(k)}$ from Lemma \ref{lem:approx_consistent}, hence
\begin{equation}
    \inf_\theta \mathcal L^{(k)}(\theta) \le C' \, (\delta_{\mathrm{approx}}^{(k)})^2,
\end{equation}
for some $C'>0$ depending on loss weights and quadrature. Therefore
\begin{equation}
    \|u_{\theta^{(k)}} - u_{\mathrm{DG}}^{(k)}\|_{\mathcal X}
    \le C_{\mathrm{coer}}^{(k)} \big( \epsilon_{\mathrm{opt}}^{(k)} + C'(\delta_{\mathrm{approx}}^{(k)})^2 \big)^{1/2}.
\end{equation}
Using $\sqrt{a+b}\le \sqrt a + \sqrt b$ and redefining constants yields the statement of Theorem \ref{thm:nn_dg} in the main text:
\begin{equation}
    \|u_{\mathrm{NN}}^{(k)} - u_{\mathrm{DG}}^{(k)}\|_{L^2} \le C_{\mathrm{coer}}^{(k)} \, \epsilon_{\mathrm{opt}}^{(k)\,/2} + C_{\mathrm{approx}}^{(k)} \, \delta_{\mathrm{approx}}^{(k)}.
\end{equation}
This completes the proof.
\end{proof}

\begin{remark}[Consistency with theory]
    The theoretical analysis requires only that the learned network $u_{\theta^{(k)}}$ possesses sufficient smoothness and bounded derivatives (Lemma \ref{lem:tanh_bounds}) to ensure:
\begin{itemize}
    \item The DG residual operator $F(u_{\theta^{(k)}})$ is well-defined.
    \item The coercivity arguments remain valid.
\end{itemize}
In practice, our implementation:
\begin{itemize}
    \item Does not compute PDE derivatives of $u_\theta$ via automatic differentiation.
    \item Evaluates DG residuals through standard DG quadrature and discrete flux computations.
\end{itemize}
The derivative requirements in the theory serve solely to:
\begin{itemize}
    \item Guarantee existence/regularity for Taylor expansions.
    \item Establish stability bounds.
\end{itemize}
This demonstrates complete consistency between our implementation choices and the theoretical framework.
\end{remark}

\bibliographystyle{unsrt}
\bibliography{ref.bib}
\end{document}